\documentclass{article}

\usepackage{amsfonts, amssymb}
\usepackage{amsmath, eqnarray}
\usepackage{color}
\usepackage{indentfirst}
\usepackage{verbatim}   
\usepackage{mathrsfs}   
\usepackage{yhmath}     
\usepackage{pb-diagram}

\usepackage{enumerate}

   \usepackage[all,cmtip]{xy}
   \usepackage[width=445pt, height=651pt, hcentering, vcentering]{geometry}
   \usepackage{tikz}
\newcommand{\g}{\mbox{$\bf g$}}

\newcommand{\h}{\mbox{\textbf{h}}}

\newcommand{\al}{\alpha}

\newcommand{\la}{\lambda}

\newcommand{\Th}{\Theta}

\newcommand{\mb}{\mbox}

\newcommand{\Mklz}[2]{\left\{\left. #1 \right| #2 \right\}}
\newcommand{\F}{\mathbb{F}}

\newcommand{\N}{\mathbb{N}}
\newcommand{\Nn}{\mathbb{N}_0}
\newcommand{\Q}{\mathbb{Q}}

\newcommand{\R}{\mathbb{R}}

\newcommand{\Z}{\mathbb{Z}}
\newcommand{\We}{\mbox{$\mathcal W$}}

\newcommand{\Cq}{\mbox{$\overline{C}$}}
\newcommand{\FTq}{\mbox{$\overline{F}_\Theta$}}
\newcommand{\ve}[1]{\mbox{$\varepsilon (#1)$}} 
\newcommand{\Fa}[1]{ Fa(#1) }
\newcommand{\cc}[1]{\mbox{cc}(#1)}
\newcommand{\co}[1]{\mbox{co}(#1)}
\newcommand{\ti}{\widetilde}
\newcommand{\res}[1]{\!\mid_{#1}}

\newcommand{\red}[1]{\mbox{red} (#1)}

\newcommand{\Proof}{\mbox{\bf Proof: }}
\newcommand{\qed}{\mb{}\hfill\mb{$\square$}\\}

\newtheorem{thm}{Theorem}[section]
\newtheorem{defn}[thm]{Definition}
\newtheorem{prop}[thm]{Proposition}
\newtheorem{cor}[thm]{Corollary}
\newtheorem{rem}[thm]{Remark}

\newtheorem{lem}[thm]{Lemma}

\begin{document}

\title{Invariant convex subcones of the Tits cone\\ of a linear Coxeter group}
\author{Claus Mokler} 
\date{}
\maketitle
\begin{abstract}\noindent We investigate the faces and the face lattices of arbitrary Coxeter group invariant convex subcones of the Tits cone of a linear Coxeter system as introduced by E. B. Vinberg. Particular examples are given by certain Weyl group invariant polyhedral cones, which underlie the theory of normal reductive linear algebraic monoids as developed by M. S. Putcha and L. E. Renner. We determine the faces and the face lattice of the Tits cone and the imaginary cone, generalizing some of the results obtained for linear Coxeter systems with symmetric root bases by M. Dyer, and for linear Coxeter systems with free root bases by E. Looijenga, P. Slodowy, and the author.
\end{abstract}
%
%
%
%
%
%
%
%
%
%
\section*{Introduction}

Linear Coxeter systems have been introduced by E. B. Vinberg in \cite{V}, generalizing the geometric representations of Coxeter systems described by J. Tits. We briefly state their definition, and some of the properties given in \cite{V}.

Let $\h$ be a finite dimensional real linear space. Every zero-based ray $H\subseteq \h$ determines dually a closed half-space $H^{\geq 0}\subseteq \h^*$, an open half-space $H^{> 0}\subseteq\h^*$, and a hyperplane $H^{=0}\subseteq \h^*$. Let $n\in\N$ and set $I:=\{1,\ldots, n\}$. A linear Coxeter system consists of a family $ (H_i,\,L_i)_{i\in I}$ of zero-based rays $H_i$, $L_i$ in $\h$, $\h^*$ with $L_i\setminus\{0\}\subseteq H_i^{>0}$, $i\in I$, such that the following holds: The closed fundamental chamber, which is the convex cone 
\begin{equation*}
   \overline{C} := \bigcap_{i\in I} H_i^{\geq 0},
\end{equation*}
spans $\h^*$, and $ \overline{C} \cap H_i^{=0}$, $i\in I$, are its pairwise different 1-codimensional faces. The linear Coxeter group $\We$ is the group generated by the reflections $\sigma_i$ at the hyperplanes $H_i^{=0}$ along the lines $L_i\cup -L_i$, $i\in I$. If we set $ C := \bigcap_{i\in I} H_i^{> 0}$, which coincides with the interior of $\overline{C}$, then
\begin{equation*}
     \sigma C\cap C=\emptyset  \quad\mb{ for all }\quad \sigma\in\We \setminus\{1\} .
\end{equation*}

Some of the properties of linear Coxeter systems are the following: $(\We,\, (\sigma_i)_{i\in I})$ is a Coxeter system, and $\overline{C}$ is a fundamental domain for the action of $\We$ on the convex cone obtained by  
\begin{equation*}
  X:= \bigcup_{\sigma\in {\mathcal W}} \sigma \overline{C},
\end{equation*}
which is called the Tits cone. 
A set $J\subseteq I$ is called facial if there exists an element $\la\in\overline{C}$ such that its isotropy group $\mb{stab}_{\mathcal W}(\la)$ is the standard parabolic subgroup $\We_J$ of $\We$ associated to $J$. There is a partition
\begin{equation*}
    \overline{C}  =   \dot{\bigcup_{J\text{ facial} }} F_J \quad \mb{ where }  \quad F_J   :=   \{  \la\in\overline{C} \mid \mb{stab}_{\mathcal W}(\la)=\We_J \}  .
\end{equation*}
It coincides with the partition given by the relative interiors of the nonempty faces of $\overline{C}$, which are $\overline{F}_J$, $J$ facial.

If we choose for every $i\in I$ elements $h_i\in H_i$, $\al_i\in L_i$ such that $\al_i(h_i)=2$, we get a root base
\begin{equation}\label{intrballg1}
  h_1,\ldots,h_n \in \h \quad \mb{ and }  \quad \al_1,\ldots, \al_n \in \h^*
\end{equation}
over $\R$, such that $h_1,\ldots, h_n$ are positively independent and
\begin{equation}\label{intrballg2}
 A:= (\al_j(h_i))_{i,j\in I}
\end{equation}
is a generalized Cartan matrix. The term 'generalized Cartan matrix' is used here in a more general way as in Lie theory, also certain noninteger entries are possible. It is important to note that for every $i\in I$ the pair $h_i$, $\al_i$ is only determined up to a positive real number. The generalized Cartan matrix $A$ is only determined up to conjugation with a diagonal matrix of positive real numbers. Nevertheless, it contains much information. In particular, the relations of the Coxter generators $(\sigma_i)_{i\in I}$ can be computed from $A$.
Conversely, a root base of this sort gives a linear Coxeter system.

It is also important to note that the definition of a linear Coxeter system is not symmetric if the roles of $\h$ and $\h^*$ are reversed, which can also be read off from the corresponding root bases. The linear Coxeter group $\We$ acts dually on $\h$, and we get a convex cone by
\begin{equation*}
  X^\vee := \bigcup_{\sigma\in {\mathcal W}} \sigma \overline{C}^\vee \quad\mb{ where }\quad \overline{C}^\vee:=  \bigcap_{i\in I} L_i^{\geq 0}.
\end{equation*}
However, to interpret $X^\vee$ as a Tits cone requires, in general, to consider the linear span of $X^\vee$ and a subsystem of $(L_i, H_i)_{i\in I}$.

The imaginary cone, which generalizes the cone of imaginary roots in Lie theory, has been introduced by M. Dyer in \cite{D}. It is the $\We$-invariant convex subcone of $-X^\vee$ obtained by 
\begin{equation*}
   Z:= \bigcup_{\sigma\in {\mathcal W}} \sigma K \quad\mb{with}\quad K:=  (\sum_{i\in I}H_i )\cap(- \overline{C}^\vee).
\end{equation*}

In this article we investigate the faces and the face lattices of arbitrary $\We$-invariant convex subcones of the Tits cone $X$, which contain the zero. One of the reasons for these investigations is that the results can be used to construct and describe the analogues of normal reductive linear algebraic monoids whose unit groups are Kac-Moody groups. We explain this briefly:

Let $M$ be a normal reductive linear algebraic monoid. Its unit group $G$ is a reductive linear algebraic group. We choose a maximal torus $T$ of $G$. Its cocharacter group $X_*(T)$ and its character group $X^*(T)$ are dual free abelian groups of finite rank on which the Weyl group $\We=N_G(T)/T$ acts faithfully. We choose in addition a Borel subgroup $B$ of $G$, which contains $T$. This determines a set of simple coroots $\Pi^\vee\subseteq  X_*(T) $, a set of simple roots $\Pi\subseteq  X^*(T) $, and a corresponding set of simple reflections $S\subseteq \We$. If we now number the elements of $\Pi^\vee$ and $\Pi$ accordingly, we get a root base
\begin{equation}\label{intrb1}
  h_1,\ldots,h_n \in X_*(T) \quad \mb{ and }  \quad \al_1,\ldots, \al_n \in X^*(T) 
\end{equation}
over $\Z$, such that 
\begin{equation}\label{intrb2}
 A:= (\al_j(h_i))_{i,j\in I}
\end{equation}
is an integer generalized Cartan matrix, whose components are of finite type. By taking the $\R$-linear hulls $X_*(T)_{\R}$ and $X^*(T)_{\R}$ we get a free root base over $\R$, i.e., $h_1,\ldots,h_n$ and  $\al_1,\ldots,\al_n $ are linearly independent. The Tits cone $X$ coincides with the linear space $X^*(T)_{\R}$.

Let  $\overline{T}$ be the Zariski closure of $T$ in $M$. The monoid of characters $X^*(\overline{T})$ of $\overline{T}$ identifies with a $\We$-invariant subset of $X^*(T)$. It generates a rational $\We$-invariant polyhedral cone 
\begin{equation}\label{intrb3}
       X^*(\overline{T})_{\R^+_0}\subseteq  X^*(T)_{\R}=X,
\end{equation}
whose interior is nonempty.

The tuple $M$, $G$, $B$, $T$ is determined by, and may be constructed up to suitably defined isomorphy from this root base over $\Z$ and this subcone of the Tits cone; see the classifications of L. E. Renner in \cite{R1} and Section 5.1 of \cite{R3}, or that of E. B. Vinberg in \cite{V2}, or see Section 6.2 of \cite{BrKu}.

In a subsequent article we construct similarly the analogues of normal reductive linear algebraic monoids whose unit groups are Kac-Moody groups, starting with certain $\Z$-root bases with integer generalized Cartan matrices and certain subcones of the Tits cone $X$. This has been done for certain free $\Z$-root bases with integer generalized Cartan matrices and the Tits cone itself in \cite{M2}. An investigation of some other examples has been started in \cite{LL}.

We now list some of the results of this article, and give some remarks to their proofs. For a root base as in (\ref{intrballg1}), (\ref{intrballg2}) 
the space of linear relations in $h_1,\,\ldots,\,h_n$ is defined as
\begin{equation*}
   L_{\bf h}:=\Mklz{r\in\R^n}{r_1h_1+\cdots+r_n h_n=0},
\end{equation*}
which is a linear subspace of the kernel of $A^T$. Since $h_1,\ldots,h_n$ are only determined up to positive real numbers, only the sign vectors of the elements of $L_{\bf h}$ are relevant. For $r\in  L_{\bf h}$ we set
\begin{equation*}
  I_+(r):=\Mklz{i\in I}{r_i>0} \quad \mb{ and } \quad I_-(r):=\Mklz{i\in I}{r_i<0}.
\end{equation*}
We find in Corollary \ref{detfacial} of Subsection \ref{SubsectionFF}, that a set $J\subseteq I$ is facial if and only if for all $r\in L_{\bf h}$ the inclusion $I_+(r)\subseteq J$ implies the inclusion $I_-(r)\subseteq J$. This criteria is advantageous for the following reasons:
\begin{itemize}
\item The sign vectors of the elements of $L_{\bf h}$ can be obtained quickly from a hyperplane arrangement associated to $L_{\bf h}$. As a consequence, the facial sets can be easily computed.
\item The properties of the generalized Cartan matrices imply properties of the sign vectors for the elements of the kernel of the matrix $ A^T$. General theoretical properties of the facial sets can be obtained as consequences.
\end{itemize}

In Subsections \ref{SubsectionFT}, \ref{SubsectionFI} we determine the faces of the Tits cone $X$, as well as of the imaginary cone $Z$, and describe the lattice operations of their face lattices. In particular, we show in Theorems \ref{FXf} (b), \ref{FZf} (b) that cross sections of the $\We$-orbits of the nonempty faces of $X$, $Z$ can be parametrized as
\begin{equation}\label{intfaTI}
   R(\Th), F(\Th),\; \Th\mb{ special facial},
\end{equation}
where $R(\Th)\subseteq X$ and $F(\Th)\subseteq Z$ are certain dual faces, and a special facial set $\Th$ is a facial set which is either empty or has only connected components of nonfinite type. For linear Coxeter systems with free root bases and integral generalized Cartan matrices the results on the Tits cone have been obtained in \cite{Loo}, \cite{Sl}, \cite{M}, \cite{M6}. For linear Coxeter systems with symmetric root bases, i.e., $ h_1,\ldots,h_n$ as well as $ \al_1,\ldots,\al_n$ positively independent and $A$ symmetric, the results on the Tits cone and imaginary cone have been obtained in \cite{D} as part of the analysis of the imaginary cone. 

We proceed here as follows: We first treat the faces and the face lattice of the Tits cone $X$, because most proofs in \cite{Loo}, \cite{Sl}, \cite{M}, \cite{M6} work with only minor changes also in the general case. Then we transfer by duality as many results as possible to the faces and the face lattice of the imaginary cone $Z$. The faces (\ref{intfaTI}) are found and described quickly. The difficult part is to show that up to $\We$-translations these are all faces of $X$, $Z$. This is done in a common way, motivated by results obtained for general invariant convex subcones of the Tits cone in Section  \ref{SecInvsub}.

Every root base can be obtained as a subquotient of a free root base. The methods developed in Sections \ref{SectionPCG} and \ref{SectionRoMM} allow to push and pull results on the facial sets, the faces of the Tits cone, and the faces of the imaginary cone from free root bases to general root bases. This requires some amount of notation and computations. We use this method sparingly, for a central result on the facial sets and for some central results on the faces of the imaginary cone.

Let $Y$ be a $\We$-invariant convex subcone of the Tits cone $X$, which contains the zero. We denote by $\Fa{Y}$ its face lattice of nonempty faces.
Important data to describe a normal reductive linear algebraic monoid are the cross section lattice and the type maps, which are also encoded in its Renner monoid, see \cite{P1}, \cite{PR}, \cite{R2}, and Chapter 7 of \cite{R3}. The cross section lattice and the type maps depend on the root base (\ref{intrb1}), (\ref{intrb2}) and the rational $\We$-invariant polyhedral cone (\ref{intrb3}), and can be used for the description of the faces and the face lattice of this cone. Similarly, we set
\begin{equation*}
  \Upsilon:=\Mklz{R\in\Fa{Y}}{\mb{ri}(R)\cap\overline{C}\neq\emptyset } ,
\end{equation*}
where $\mb{ri}(R)$ denotes the relative interior of $R$. For $R\in\Upsilon$ the sets
\begin{equation*}
 \upsilon_*(R):= \Mklz{i\in I}{   R-R\subseteq H^{=0}_i  } \quad \mb{ and }\quad \upsilon^*(R):= \Mklz{i\in I}{ L_i\subseteq  R-R  }
\end{equation*}
are separated, i.e., the corresponding subdiagrams of the Coxeter graph are separated, and we show that $\upsilon_*(R)$ is facial. We set $\upsilon(R):=\upsilon_*(R)\cup\upsilon^*(R)$.

Many examples can be obtained from the theory of reductive linear algebraic monoids. For the Tits cone $X$ we have $\Upsilon=\{R(\Th)\mid \Th \mb{ special facial}\}$, and if $\Th$ is special facial then
\begin{equation*}
   \upsilon_*(R(\Th))= \Th \quad \mb{ and }\quad \upsilon^*(R(\Th))= \Th^\bot,
\end{equation*}
where $\Th^\bot$ denotes the biggest subset of $I$ separated from $\Th$. Because $X^\vee$ can be interpreted as the Tits cone of a linear Coxeter system, the imaginary cone $Z\subseteq - X^\vee$ can also be included as an example. We have $\Upsilon=\{F(\Th)\mid \Th \mb{ special facial}\}$.
However, to write down $\upsilon_*(F(\Th))$ and $\upsilon^*(F(\Th))$, $\Th$ special facial, is more complicated than for the Tits cone $X$.

We find in Theorem \ref{maincsx} of Subsection \ref{maincsx} that $R\in\Upsilon$ can be described as
\begin{equation*}
  R =\We_{\upsilon^*(R)}(R\cap\overline{F}_{\upsilon_*(R)}),
\end{equation*}
its pointwise $\We$-stabilizer is $Z_{\mathcal W}(R)= \We_{\upsilon_*(R)}$, and its $\We$-stabilizer as a whole is $N_{\mathcal W}(R) =\We_{\upsilon(R)}$. We conclude in Corollary \ref{mainlh} that its linear hull is
\begin{equation*}
  R - R = (R\cap\overline{F}_{\upsilon_*(R)}) - (R\cap\overline{F}_{\upsilon_*(R)}) .
\end{equation*}
We obtain in Theorem \ref{maincsx2} that its relative interior $\mb{ri}(R)$ can be computed from $\mb{ri}(R\cap\overline{F}_{\upsilon_*(R)})$ by 
\begin{equation*}
  \mb{ri}(R)
  = \We_{\upsilon^*(R)} \dot{\bigcup_{J_f\subseteq \upsilon^*(R),\, (J_f)^0=J_f}}\underbrace{  \mb{mid}_{J_f}(\mb{ri}(R\cap\overline{F}_{\upsilon_*(R)}))  }_{\subseteq F_{\upsilon_*(R)\cup J_f}} \, ,
\end{equation*}
where $(J_f)^0$ denotes the union of the components of finite type of $J_f$, and $ \mb{mid}_{J_f}$ is a certain projector.

We show in Subsections \ref{inclusion of faces}, \ref{the cross section lattice}, \ref{relation to literature}, and \ref{lattice operations of the face lattice}  that $\Upsilon$ is a cross section for the $\We$-orbits of $\Fa{Y}$, and a complete sublattice of $\Fa{Y}$. Every chain in $\Fa{Y}$ is the $\We$-translate of a uniquely determined chain in $\Upsilon$. The order theoretic operations of $\Fa{Y}$ can be reduced to that of $\Upsilon$. In particular, if $R_1,\,R_2\in\Upsilon$ and $\sigma_1,\,\sigma_2\in\We$ then $\sigma_1 R_1\subseteq \sigma_2 R_2$ if and only if 
\begin{equation*}
   R_1\subseteq R_2 \quad\mb{ and }\quad \sigma_1^{-1}\sigma_2\in  \We_{\upsilon(R_1)}\We_{\upsilon(R_2)}.
\end{equation*}
It follows that the map $\upsilon: \Upsilon\to {\mathcal P}(I)$, where ${\mathcal P}(I)$ denotes the power set of $I$, is a regular type map as considered in the theory of monoids on groups with BN-pairs by M. S. Putcha, see for example \cite{P2}, \cite{AP}. A generalized Renner monoid as defined by E. Godelle in \cite{Go} is naturally associated.

For $R_1,\,R_2\in\Upsilon$ such that $R_1\subseteq R_2$ we equip the interval
\begin{equation*}
    [R_1,R_2 ]  := \Mklz{R\in\Fa{Y}}{R_1\subseteq R\subseteq R_2}
\end{equation*}
with the action of the factor group of its $\We$-stabilizer as a whole by its pointwise $\We$-stabilizer. We obtain in Theorem \ref{mainint} of Subsection \ref{Intervals} a description of $[R_1,R_2 ]$ as the face lattice of a convex invariant cone of the Tits cone of a certain linear Coxeter system. 

If the convex cone $Y\cap\overline{C}$ is finitely generated, we obtain in Corollary \ref{ketteninfapart} of Subsection \ref{Chains in intervals} a $\We$-invariant decomposition 
\begin{equation*}
\Fa{Y} = \dot{ \bigcup_{\Th \text{ special facial}} }  \Fa{Y}_\Th,
\end{equation*}
such that the parts $\Fa{Y}_\Th$, some of them may be empty, have the following properties: 
\begin{itemize}
\item If $R_1,\,R_2\in \Fa{Y}_\Th $ such that $R_1\subseteq R_2$, then every maximal chain between $R_1$ and $R_2$ is contained in $\Fa{Y}_\Th$ and has the length $\mb{dim}(R_2)-\mb{dim}(R_1)+1$.
\item  If $R_1,\,R_2\in \Fa{Y}_\Th \cap \Upsilon$ such that $R_1\subseteq R_2$, then
\begin{equation*}
  \upsilon^*(R_1)\,\dot{\cup}\,\upsilon_*(R_2)\,=\,\bigcap_{R\in\Upsilon\atop R_1\subseteq R\subseteq R_2} \upsilon(R)  \, .
\end{equation*}
\end{itemize}

The imaginary cone $Z^\vee\subseteq \h^*$ is defined similarly as the cone $Z\subseteq \h$. We see in Corollary \ref{dimc} of Subsection \ref{dimc as subc} that we have $-Z^\vee\subseteq Y$ if $\We$ acts faithfully on $Y$ and $Y\cap\overline{C}$ is closed.

%
%
%
\tableofcontents
%
%
%
%
\section{\label{SectionPCG}A collection of some results related to convex geometry}
In this section we give some easy lemmas and propositions related to convex geometry, which will be used later. For the basic definitions and theorems  of convex geometry we refer to \cite{Bro}, \cite{La}, and \cite{Ro}. 

We denote by $\N=\Z^+$, $\Q^+$, $\R^+$ the sets of strictly positive numbers of $\Z$, $\Q$, $\R$, and the sets $\N_0=\Z^+_0$, $\Q^+_0$, $\R^+_0$ contain, in addition, the zero. 
If $v, w$ are elements of a real linear space $V$, we denote by $\overline{v w}$ the closed line segment in $V$ with endpoints $v$ and $w$. If $M$ is a subset of $V$ we denote by $\co{M}$ the convex hull of $M$, and by $\cc{M}$ the convex cone with zero generated by $M$. If $K$ is a finite dimensional convex subset of $V$ we denote by $\mb{ri}(K)$ its relative interior. We will often use that for a finite subset $M$ of $V$ we have
\begin{equation*}
   \mb{cc}(M)=\sum_{m\in M} \R^+_0 m \quad \mb{ and }\quad  \mb{ri}\left( \mb{cc}(M) \right) =\sum_{m\in M} \R^+ m ,
\end{equation*} 
where $\sum_{m\in \emptyset} \R^+_0 m :=\{0\}$ and $\sum_{m\in \emptyset} \R^+ m :=\{0\}$.

Let $K$ be a convex cone in a real linear space $V$. If $K$ is nonempty then the linear hull of $K$ coincides with the affine hull of $K$, which is given by $K-K$. A subset $F\subseteq K$ is a face of $K$ if and only if $F$ is a convex cone and $K\setminus F + K\subseteq K\setminus F$. A subset $F\subseteq K$ is an exposed face of $K$ if and only if $K=\emptyset$, or there exists $\phi\in V^*$ such that $K\subseteq \{v\in V\mid\phi(v)\geq 0\}$ and $F= \{v\in K\mid\phi(v)= 0\}$.

The following three lemmas will be key to determine the set of faces of certain convex cones. We denote by ${\mathcal C}(V)$ be the set of convex cones in the real linear space $V$, ordered partially by inclusion. It is a complete lattice, the lattice meet given by the intersection of cones.
If $U_1$, $U_2$ are linear subspaces of $V$ such that $U_1\subseteq U_2$, then the set
\begin{equation*}
 \Mklz{ K\in {\mathcal C}(V) }{  K + U_1 \subseteq K\subseteq U_2 }
\end{equation*}
is a closed sublattice of ${\mathcal C}(V)$, i.e., it has a smallest element, a biggest element, and contains the meets and joins of its nonempty subsets. If $K\in {\mathcal C}(V) $ contains the zero, then the condition $ K + U_1 \subseteq K$ is equivalent to $U_1 \subseteq K$.

\begin{lem} \label{conint} Let $V$, $W$ be real linear spaces, and let $\phi:V\to W$ be a linear map. Then the map
\begin{eqnarray*}
 \phi: \Mklz{ K\in {\mathcal C}(V) }{K + \mb{ker}(\phi)\subseteq K} &\to &\Mklz{L\in{\mathcal C} (W)}{L\subseteq \mb{im}(\phi)\,}\\
                                 K\qquad\qquad\qquad &\mapsto &  \qquad \quad\phi(K)
\end{eqnarray*}
is a lattice isomorphism. The map
\begin{eqnarray*}
 \phi^{-1}: \Mklz{L\in{\mathcal C} (W)}{L\subseteq \mb{im}(\phi)\,} &\to &  \Mklz{ K\in {\mathcal C}(V) }{K + \mb{ker}(\phi)\subseteq K}\\
        L\qquad\qquad\quad &\mapsto & \qquad\qquad\quad \phi^{-1}(L)
\end{eqnarray*}
is the inverse lattice isomorphism.

If a convex cone is contained in $\Mklz{ K\in {\mathcal C}(V) }{K + \mb{ker}(\phi)\subseteq K} $ then its affine hull and its faces are also contained in this set. A similar statement holds for $\Mklz{L\in{\mathcal C} (W)}{L\subseteq \mb{im}(\phi)\,}$. The maps $\phi$, $\phi^{-1}$ preserve the property ``is the affine hull of''. They preserve the properties ``is a face of" and ``is an exposed face of". 

Let $V$ and $W$ be finite dimensional. If a convex cone is contained in $\Mklz{ K\in {\mathcal C}(V) }{K + \mb{ker}(\phi)\subseteq K} $ then its relative interior and its closure are also contained in this set. A similar statement holds for $\Mklz{L\in{\mathcal C} (W)}{L\subseteq \mb{im}(\phi)\,}$. The maps $\phi$, $\phi^{-1}$ preserve the property ``is the relative interior of" and ``is the closure of".
\end{lem}
\Proof The isomorphisms $\phi$, $\phi^{-1}$ are obtained by composing the isomorphisms
\begin{equation*}
 \Mklz{ K\in {\mathcal C}(V) }{K + \mb{ker}(\phi)\subseteq K} \cong {\mathcal C}(V/ \mb{ker}(\phi)) \cong   {\mathcal C}(\mb{im}(\phi)) \cong 
\Mklz{L\in{\mathcal C} (W)}{L\subseteq \mb{im}(\phi)\,}.
\end{equation*}
Most of the remaining statements of the lemma are trivial to check. We only show the following:

Let $ K\in {\mathcal C}(V)$ such that $K + \mb{ker}(\phi)\subseteq K$, and let $F$ be a face of $K$. For every $f\in F$ and $ k\in\mb{ker}(\phi)\setminus\{0\}$ the line segment between $f-k\in K$ and $f+k\in K$ contains $f$. Since $F$ is a face, we get $f\pm k\in F$. Hence $F+\mb{ker}(\phi)\subseteq F$. Clearly, $\phi(F)$ is a convex cone, and $\phi(K\setminus F)+\phi(K)=\phi((K\setminus F) + K) \subseteq \phi(K\setminus F)$. Assume that $k\in K\setminus F$ and $\phi(k)\in\phi(F)$. Then $k\in F +\mb{ker}(\phi)\subseteq F$,  which is not possible. Therefore, $\phi(K\setminus F)\subseteq \phi(K)\setminus\phi(F)$, which shows that $\phi(F)$ is a face of $\phi(K)$.

The affine hull of a nonempty convex cone coincides with its linear hull. Hence it is sufficient to show the statements on the relative interiors and the closures for a surjective linear map and generating cones:
Let $\phi:V\to W$ be a surjective linear map between finite dimensional linear spaces. Then $\phi$ is open and continuous, which is equivalent to $\phi^{-1}(L^\circ)=\phi^{-1}(L)^\circ$ for all $L\subseteq W$, which in turn is equivalent to $\phi^{-1}(\overline{L})=\overline{\phi^{-1}(L)}$ for all $L\subseteq W$. Now let $K\subseteq V$ such that $K+\mb{ker}(\phi)\subseteq K$. Then $K=\phi^{-1}(\phi(K))$. We find $K^\circ = \phi^{-1}(\phi(K)^\circ)$ and $\overline{K} = \phi^{-1}(\overline{\phi(K)})$. Since $\phi$ is surjective, we get $ \phi(K^\circ)=\phi(K)^\circ$ and $ \phi(\overline{K})=\overline{\phi(K)}$. Furthermore, $K+\mb{ker}(\phi)\subseteq K$ is equivalent $K+v\subseteq K$ for all $v\in\mb{ker}(\phi)$. It follows that $K^\circ+v\subseteq K^\circ$ and  $\overline{K}+v\subseteq \overline{K}$ for all $v\in\mb{ker}(\phi)$.
\qed

Let $K$ be a convex cone with zero in the real linear space $V$. We denote by $\Fa{K}$ the set of all nonempty faces of $K$, ordered partially by inclusion. It is a complete lattice, the lattice meet given by the intersection of faces. Its biggest face is $K$, its smallest face is the linear subspace $K\cap (-K)$ of $V$. We call $\Fa{K}$ the face lattice of $K$. The elements of $\Fa{K}$ are henceforth referred to as the faces of $K$.

If $V$ is finite dimensional, then we obtain from Theorems 6.2, 18.2, and Corollary 18.1.2 of \cite{Ro} that $K$ is the disjoint union of the relative interiors $\mb{ri}(F)$, $F\in\Fa{K}$, which are all nonempty. 

Note also that every face $F\in \Fa{K}$ is obtained from its linear hull by $F=K\cap (F-F)$. If $K=\cc{E}$ then $F=\cc {E\cap F}$.

\begin{lem} \label{consum} Let $V$ be a real linear space. Let $K$ be a convex cone with zero in $V$, and let $U$ be a linear subspace of $V$. Then
\begin{eqnarray*}
   i:\Fa{K+U}  &\to &  \Fa{K} \\
  G \qquad &\mapsto & G\cap K
\end{eqnarray*}
is an inclusion and intersection preserving injective map. It maps exposed faces to exposed faces. 
Its image is given by
\begin{equation}\label{ccimi}
   \mb{im}(i)=\Mklz{ F\in\Fa{K} \,}{ \, U\cap (K-F) \subseteq F-F \,} .
\end{equation}
The inverse map on this image is the inclusion preserving map
\begin{eqnarray*}
      i^{-1}: \mb{im}(i) &\to &  \Fa{K+U} \\
       F \;\;  &\mapsto & \quad F+U .
\end{eqnarray*}
\end{lem}
\Proof Trivially, $i$ is a well defined inclusion and intersection preserving map, which maps exposed faces to exposed faces.  Trivially, the description of the map $i^{-1}$ shows that it is inclusion preserving.

(a) We first describe the set on the right in (\ref{ccimi}) differently. We show that for $F\in\Fa{K}$ the following are equivalent:
\begin{itemize}
\item[(i)] $U\cap(K-F)\subseteq F-F$.
\item[(ii)] $(F+U)\cap K=F$.
\end{itemize}

The implication ``$(i)\Rightarrow (ii)$": It is sufficient to show $(F+U)\cap K\subseteq F$. Let $k\in (F+U)\cap K$. Write $k$ in the form 
$k=f+u$ with $f\in F$ and $u\in U$. Then
\begin{equation*}
    U\ni u = k - f\in K-F,
\end{equation*}
and by (i) we find $u = k-f\in F-F$. It follows that
\begin{equation*}
  K\ni k = f+u\in f+(F-F)\subseteq F-F. 
\end{equation*}
Since $F$ is a face of $K$, we get $k\in K\cap(F-F)=F$.

The implication ``$(ii)\Rightarrow (i)$": Let $u\in U\cap (K-F)$. Write $u$ in the form $u= k-f$ with $k\in K$ and $f\in F$. Then
\begin{equation*}
  F+U\ni f+u= k\in K
\end{equation*}
and by (ii) we find $f+u=k\in F$. Hence $u=k-f \in F-F$.

(b) Let $G\in\Fa{K+U}$. We show that $G\cap K\in\Fa{K}$ is contained in the set on the right in (\ref{ccimi}), i.e., we show that
\begin{equation*}
    ((G\cap K)+U)\cap K= G\cap K.
\end{equation*} 
Trivially, the inclusion ``$\supseteq$" holds. Since the smallest face of $K+U$ is contained in $G$, we find
\begin{equation}\label{zwccimi}
  U\subseteq (K+U)\cap(-(K+U))\subseteq G .
\end{equation}
It follows that $G+U\subseteq G$ and
\begin{equation*}
   ((G\cap K)+U)\cap K\subseteq G\cap K.
\end{equation*}

(c) Let $F\in\Fa{K}$ such that $F=(F+U)\cap K$. We show that $F+U\in\Fa{K+U}$. The sum of convex cones $F+U$ is a convex cone. Let $v\in (K+U)\setminus ( F+U)$ and $v'\in K+U$. Then $v$, $v'$ can be written in the form $v=k+u$ and $v'=k'+u'$ with $k\in K\setminus F$, $k'\in K$, $u,\,u'\in U$. Suppose that $v+v'\in F+U$. Then there exist $f\in F$, $u''\in U$ such that 
\begin{equation*}
v+v'= k+k'+u+u'=f+u''.
\end{equation*} 
Since $F$ is a face of $K$, it follows that
\begin{equation*}
  K\setminus F\ni k+k'=f+(u''-u-u')\in F+U,
\end{equation*}
which contradicts $F=(F+U)\cap K$.

(d) Let $F\in\Fa{K}$ such that $F=(F+U)\cap K$. Then 
\begin{equation*} 
   i (i^{-1}(F))=i(F+U)=(F+U)\cap K=F.
\end{equation*}
Let $G\in\Fa{K+U}$. By (\ref{zwccimi}) we have $U\subseteq G$. Hence 
\begin{equation*} 
   i^{-1}(i(G))=i^{-1}(G\cap K)=(G\cap K)+U \subseteq G .
\end{equation*}
To show the reverse inclusion write $v\in G$ in the form $v=k+u$ with $k\in K$ and $u\in U$. Then 
\begin{equation*}
    K\ni k=v-u\in G-U=G+U\subseteq G ,
\end{equation*}
which shows that $v\in (G\cap K)+U$.
\qed

\begin{lem}\label{concap} Let $V$ be a finite dimensional real linear space. Let $K$ be a convex cone with zero in $V$, and let $U$ be a linear subspace of $V$. Then
\begin{eqnarray*}
  p:\, \Fa{K}  &\to & \Fa{K\cap U} \\
  F \quad &\mapsto & \quad F\cap U
\end{eqnarray*}
is an inclusion and intersection preserving surjective map. It maps exposed faces to exposed faces.

Let $H\in \Fa{K\cap U} $. The fiber $p^{-1}(H)$ contains a smallest element $i(H)\in \Fa{K}$. It is the smallest face of $K$ which contains $H$. The map
\begin{eqnarray*}
  i:\, \Fa{K\cap U}  &\to &  \Fa{K} 
\end{eqnarray*}
is an inclusion preserving injective map  with image
\begin{equation*}
  \mb{im}(i)=\Mklz{F\in\Fa{K}}{ \mb{\rm ri}(F)\cap U \neq\emptyset } , 
\end{equation*}
and we have $p\circ i=id$. 
\end{lem}
\begin{rem} Since $p: \Fa{K} \to \Fa{K\cap U}$ and $i: \Fa{K\cap U} \to  \Fa{K}$ are inclusion preserving maps with $p\circ i=id$, the lattice operations in $\Fa{K\cap U}$ can be obtained from the lattice operations in $\Fa{K}$ as follows. For all $G,H\in \Fa{K\cap U}$ we have:
\begin{itemize}
\item[(a)] $H\subseteq G$ if and only if $i(H)\subseteq i(G)$.
\item[(b)] $ H\cap G = p(\,i(H)\cap i(G)\,) $ and $ H\vee G = p(\,i(H)\vee i(G)\,)$.
\end{itemize}
\end{rem}
\Proof Trivially, $p$ is a well defined inclusion and intersection preserving map, which maps exposed faces to exposed faces. We have
\begin{equation*}
  K\cap U  =  (\,\dot{\bigcup_{F\in\Fa{K}}} \mb{ri}(F) \,)\cap U=  \dot{\bigcup_{F\in\Fa{K}, \,\text{ri} (F)\cap U\neq\emptyset}} \text{ri}(F)\cap U = \dot{\bigcup_{F\in\Fa{K}, \,\text{ri}(F)\cap U\neq\emptyset}} \mb{ri}(F\cap U) ,
\end{equation*}
where we used Theorem 6.5 of \cite{Ro} and $\mb{ri}(U)=U$. It follows that
\begin{equation*}
  \Fa{K\cap U}=\Mklz{F\cap U}{F\in\Fa{K},\,\mb{ri}(F)\cap U\neq \emptyset}.
\end{equation*}
In particular, the map $p$ is surjective.

Let $F\in\Fa{K}$ such that $\mb{ri}(F)\cap U\neq\emptyset$. Then $\mb{ri}(F)\cap\mb{ri}(F\cap U)=\mb{ri}(F)\cap U\neq\emptyset$, which shows that $F$ is the smallest face of $K$ containing $F\cap U$. Trivially, $F$ is also the smallest element in the fiber $p^{-1}(F\cap U)$. 
We conclude that $\mb{im}(i)=\Mklz{F\in\Fa{K}}{ \mb{\rm ri}(F)\cap U \neq\emptyset }$, and $p\circ i=id$. In particular, the map $i$ is injective.

Now let $H_1,H_2\in  \Fa{K\cap U}$ such that $H_1\subseteq H_2$. Then $H_1\subseteq H_2\subseteq i(H_2)$. Since $i(H_1)$ is the smallest face of $K$ containing $H_1$, we find $i(H_1)\subseteq i(H_2)$.
\qed

A family $(v_l)_{l\in L}$ of elements of a real linear space $V$ is called positively independent if for $(r_l)_{l\in L}\in (\R^+_0)^L$ with $r_l\neq 0$ for at most finitely many $l\in L$, the equation
\begin{equation*}
  \sum_{l\in L} r_l v_l = 0
\end{equation*}
implies $r_l= 0$ for all $l\in L$. It is easy to show:
\begin{prop}\label{pind} Let $V$ be a real linear space. For a family $(v_l)_{l\in L}$ of elements of $V$ the following are equivalent:
\begin{itemize}
\item[(i)]  $(v_l)_{l\in L}$ is positively independent.
\item[(ii)] $v_l \neq 0$ for all $l \in L$, and $\{0\}$ is a face of $\text{\rm cc}(v_l\mid l\in L)$.
\end{itemize}  
\end{prop}

Positive independence is often required for the simple roots and  the simple coroots in the definition of root data systems. 
A possibility to strengthen the notion of positive independence for infinite families is to take (ii) of Proposition \ref{pind} as definition, but to replace ``face" by some stronger condition. For example, the simple roots and simple coroots of the root data in Section 5.1 of \cite{MoPi} satisfy (ii) with ``face" replaced by condition (c) of the following proposition. The simple roots of the root data in Section 1.1 of \cite{CR} satisfy (ii) with ``face" replaced by condition (b) of the following proposition.
\begin{prop} Let $V$ be a real linear space, and let $K$ be a convex cone with zero in $V$. Consider the following statements:
\begin{itemize}
\item[(a)] $\{0\}$ is a face of $K$.
\item[(b)] $\{0\}$ is an exposed face of $K$.
\item[(c)] There exists a base $(b_j)_{j\in J}$ of $K-K$ such that $K\subseteq \sum_{j\in J}\R^+_0 b_j$.
\end{itemize}
Then we have the following implications:
\begin{equation*}
   (c) \Longrightarrow (b)  \Longrightarrow (a).
\end{equation*}
In general, the reverse implications do not hold. If $K$ is finitely generated, then the statements (a), (b), and (c) are equivalent.
\end{prop}
\Proof If $K=\{0\}$ then (b) and (c) hold, where for (c) we use the empty base, and the definition that a sum over the empty set is zero. We assume $K\neq\{0\}$ in the remaining proof. 

Suppose that (c) holds. Define $\phi\in V^*$ by $\phi(b_j):=1$ for all $j\in J$. Then 
\begin{equation*}
  K\setminus\{0\}\subseteq \sum_{j\in J}\R^+_0 b_j \setminus\{0\} \quad\mb{ and }\quad \phi(\sum_{j\in J}\R^+_0 b_j \setminus\{0\})\subseteq \R^+,
\end{equation*}
which shows (b). The convex cone 
\begin{equation*}
   K:=\Mklz{(r_1,r_2)\in\R^2}{r_2>0}\cup \{(0, 0)\}
\end{equation*}
satisfies (b) but not (c).  Clearly, (b) implies (a).
The convex cone 
\begin{equation*}
   K:=\Mklz{(r_1,r_2)\in\R^2}{r_2>0}\cup \Mklz{ (r_1, 0) }{ r_1\in\R^+_0}
\end{equation*}
satisfies (a) but not (b). 

Let $K$ be finitely generated, i.e., $K=\cc{k_1,\,\ldots,\, k_m}$. Then (a) implies (b) by Theorem 4.23 and Proposition 4.3 of \cite{La}. Now suppose that (b) holds. Then there exists $\phi\in V^*$ such that 
\begin{equation}\label{efK}
   K\setminus\{0\}\subseteq \Mklz{v\in V}{\phi(v)>0}.
\end{equation}
 Because of $K\neq\{0\}$ we may assume $\phi(k_1)=\cdots=\phi(k_m)=1$. We may assume $V=K-K$. We choose a base $a_1,\,\ldots,\,a_n$ of $\Mklz{v\in V}{\phi(v)=0}$, and an element $a_0\in V$ such that $\phi(a_0)=1$. We define
\begin{equation*} 
   T := \co{-\xi,\, a_1-\xi,\,\ldots,\, a_n-\xi} \quad \mb{ with } \quad \xi:=\frac{1}{n+1}(a_1+\cdots+a_n),
\end{equation*}
which is an n-dimensional simplex in the n-dimensional space $\Mklz{v\in V}{\phi(v)=0}$. Moreover,
\begin{equation*}
     0=  -\frac{1}{n+1}\xi +\frac{1}{n+1}(a_1-\xi) + \cdots +\frac{1}{n+1} (a_n-\xi)
\end{equation*} 
is in the interior of $T$ in $\Mklz{v\in V}{\phi(v)=0}$. Therefore, there exists an element $r\in \R^+$ such that 
\begin{equation*}
    \frac{1}{r} (k_i - a_0) \in  T \quad\mb{ for all }  \quad i=1,\,2,\,\ldots,\,m.
\end{equation*}
It is trivial to check that 
\begin{equation*}
  b_0:=a_0-r\xi,\; b_1:= a_0+r(a_1-\xi),\;\ldots,\;  b_n:=a_0+ r(a_n-\xi)
\end{equation*}
is a base of $V$, such that $K\subseteq \sum_{j=0}^n\R^+_0 b_j$, i.e., (c) holds.
\qed

Linear Coxeter systems are defined in the literature in different ways, and some results of the corresponding theories of linear Coxeter systems are used to prove that these definitions are equivalent. The following proposition allows to show the equivalence directly.
\begin{prop}\label{fcequiv} Let $V$ be a finite dimensional linear space, and $h_1,\,\ldots,\,h_n\in V$. The following are equivalent:
\begin{itemize}
\item[(i)]  The polyhedral cone
\begin{equation*}
   \overline{C}=\Mklz{ \la\in V^* }{ \la(h_i)\geq 0 \mb{ for all } i=1,\,2,\,\ldots,\,n }
\end{equation*}
spans $V^*$, and
\begin{equation*}
    \overline{C} \cap \Mklz{\la\in V^*}{\la(h_i)=0},\quad i=1,\,2,\,\ldots,\,n,
\end{equation*}
are its pairwise different  1-codimensional faces.
\item[(ii)] The convex cone $ \sum_{i=1}^n\R^+_0 h_i$ has $\{0\}$ as a face, and $ \R^+_0 h_i$, $ i=1,\,2,\,\ldots,\,n$, are its pairwise different 1-dimensional faces.
\item[(iii)]  $h_1,\,h_2,\,\ldots,\,h_n$ are positively independent, and
\begin{equation*}
    h_i\notin   \sum_{j\in \{1,\,\ldots,\,n\},\,j\neq i} \R^+_0 h_j \quad \mb{ for all }\quad i=1,\,2,\,\ldots,\,n.
\end{equation*}
\end{itemize}
\end{prop}
\Proof The cones $\overline{C}$ and $\sum_{i=1}^n \R^+_0 h_i$ are dual. The equivalence of (i) and (ii) follows easily from Proposition I.2 (5) and (6) of \cite{N}.

``$(ii)\Rightarrow (iii)$": We have $h_i\neq 0$ because $\R^+_0 h_i$ is one-dimensional, $i\in  \{1,\,\ldots,\,n\}$. Since $\{0\}$ is a face of  $\sum_{k=1}^n\R^+_0 h_k$, we find from Proposition \ref{pind} that $h_1,\,\ldots,\,h_n$ are positively independent.

Suppose that there exists $i\in \{1,\,\ldots,\,n\}$ such that 
\begin{equation*}
   h_i=\sum_{j\in \{1,\,\ldots,\,n\}, \,j\neq i}r_j h_j  \quad\mb{ with }\quad r_j\in \R^+_0 .
\end{equation*} 
Since $\R^+_0 h_i$ is a one-dimensional face of $\sum_{k=1}^n\R^+_0 h_k$, we find that there exists $j\in  \{1,\,\ldots,\,n\}$, $j\neq i$, such that $\R^+_0 h_i=\R^+_0 h_j$. But this is not possible because these one-dimensional faces are different.

``$(iii)\Rightarrow (ii)$": It follows from Proposition \ref{pind} that $\{0\}$ is a face of  $\sum_{i=1}^n\R^+_0 h_i$. From both conditions in (iii) we obtain
\begin{equation}\label{d-sum}
       \R h_i \cap \sum_{j\in \{1,\,\ldots,\,n\},\,j\neq i} \R^+_0 h_j = \{0\}.
\end{equation}
Let $x, y\in\sum_{j=1}^n\R^+_0 h_j$ such that $x+y\in\R^+_0 h_i$. Write $x=rh_i+x_r$, $y=sh_i+y_r$ with $r,s\in \R^+_0$, and $x_r, y_r\in \sum_{j\in \{1,\,\ldots,\,n\},\,j\neq i} \R^+_0 h_j$. From (\ref{d-sum}) we find $x_r + y_r=0$. Since $\{0\}$ is a face of $\sum_{j=1}^n\R^+_0 h_j$, we get $x_r=y_r=0$, which shows that $x, y\in\R^+_0 h_i$. Hence the convex cone $\R^+_0 h_i$ is a face of $\sum_{j=1}^n\R^+_0 h_j$. It is 1-dimensional because of $h_i\neq 0$.
It also follows from (\ref{d-sum}) that $\R^+_0 h_1$, \ldots, $\R^+_0 h_n$ are pairwise different. 

Now let $F$ be a 1-dimensional face of $\sum_{j=1}^n\R^+_0 h_j$. Because of $F=\text{cc}(F\cap\{h_1,\ldots, h_n\})$ there exists $i \in\{1,\ldots,n\}$ such that $h_i\in F$. Since $F$ is 1-dimensional, we find $F-F=\R h_i$. Hence
\begin{equation*}
    F =  (F-F)\cap \sum_{j=1}^n\R^+_0 h_j =  \R h_i \cap \sum_{j=1}^n\R^+_0 h_j  = \R^+_0 h_i .
\end{equation*} 
\qed
%
%
%
\section{\label{SectionRoMM}Realizations of matrices and their morphisms}
%
%
%
Realizations of matrices are used in Lie theory to construct contragredient Lie algebras, which include Kac-Moody algebras, and even Borcherds algebras. See for example Chapter 1 and Section 11.13 of \cite{K}, Section 4.2 of \cite{MoPi}, and Chapter 3 of \cite{Sl}. Certain realizations of matrices are also used in Coxeter theory, more commonly called root bases.

E. B. Vinberg classifies in \S 5 of \cite{V} all finite dimensional realizations of a matrix of finite size with entries in $\R$ by a datum called the characteristic of the realization. His proofs work without change over an arbitrary field. For his classification he introduces isomorphisms, but not morphisms of realizations. In Lie theory morphisms of realizations are introduced, see for example Chapter 3 of \cite{Sl}, but the characteristics of realizations are not used. Often only free realizations are considered.
In this section we combine both concepts. We describe how the kernels and images of certain morphisms are related to the characteristics of the corresponding realizations. 

For the whole section $A=(a_{ij})_{i,\,j=1,\ldots,\,n}$ is an $ n\times n$-matrix with entries in a field $\F$ of arbitrary characteristic, and $n\in\N$. We set $I:=\{1,\,\ldots,\,n\}$. For $r\in\F^n$ we denote by $r_1$, \ldots,\,$r_n$ its components. 

A realization  $(\h,(h_i)_{i\in I},(\al_i)_{i\in I})$ of $A$ consists of a finite dimensional $\F$-linear space $\h$, and families $h_1,\,\ldots,\, h_n\in \h $, 
$\al_1,\,\ldots,\,\al_n \in \h^*$, such that 
\begin{equation*}
     \al_i(h_j)=a_{ji}\quad \mb{for all } i,j \in I.
\end{equation*} 
Its characteristic $(L_{\bf h},\,L_{{\bf h}^*},\,d)$ is defined by 
\begin{eqnarray*}
   & L_{\bf h}:=\Mklz{ r\in\F^n \,}{r_1h_1+\cdots+r_n h_n=0}  \subseteq \mb{ker}(A^T),\\
   & L_{{\bf h}^*}:= \Mklz{r \in \F^n\,}{r_1\al_1+\cdots+ r_n \al_n=0}   \subseteq  \mb{ker}(A), \\
& d:=\mb{dim}(\h)- \mb{dim}(\mb{span}(\al_1,\,\ldots,\,\al_n))-\mb{dim}(\mb{span}(h_1,\,\ldots,\,h_n))+\mb{rk}(A)  .
\end{eqnarray*}
The linear spaces $L_{\bf h}$, $L_{{\bf h}^*}$ are called the spaces of linear relations of $h_1,\ldots, h_n$, of $\al_1,\ldots,\al_n$, respectively. The integer $d$ is called the defect. It is not difficult to see that
\begin{equation*}
  d  = \mb{dim} \big( \mb{span}(\al_1,\,\ldots,\,\al_n)^\bot  / (\mb{span}(h_1,\,\ldots,\,h_n) \cap \mb{span}(\al_1,\,\ldots,\,\al_n)^\bot)  \big) \in\Nn .
\end{equation*}

Let $(\h,(h_i)_{i\in I},(\al_i)_{i\in I})$ and  $(\h', (h_i')_{i\in I},(\al_i')_{i\in I})$ be realizations of $A$. A linear map $\phi:\h\to\h'$ is called a morphism if 
\begin{equation*} 
\phi(h_i)=h_i'\quad \mb{ and }\quad \phi^* (\al_i ')=\al_i \quad\mb{ for all }i\in I,
\end{equation*}
where $\phi^*:{\h'}^*\to \h^*$ is the dual map of $\phi:\h\to\h'$. For the corresponding characteristics we have
\begin{equation*}
   L_{\bf h}\subseteq L_{\bf h'},\quad     L_{{\bf h}^*}\supseteq\,L_{{\bf h'}^*}   .
\end{equation*}

A morphism of realizations  $\phi:\h\to\h'$ is an isomorphism if and only if the map  $\phi:\h\to\h'$ is bijective. In this case the corresponding characteristics coincide, i.e,
\begin{equation*}
   L_{\bf h}= L_{\bf h'},\quad     L_{{\bf h}^*}=L_{{\bf h'}^*} \quad d=d'  .
\end{equation*}

The following theorem has been given as Proposition 15 in \cite{V}.
\begin{thm}\label{mainer}Let $ L_{\bf h}$ be a subspace of $\mb{ker}(A^T)$, let $ L_{{\bf h}^*}$ be a subspace of $\mb{ker}(A)$, and let $d\in\Nn$. There exists a realization $(\h,(h_i)_{i\in I},(\al_i)_{i\in I})$ of $A$, unique up to isomorphism, having $(L_{\bf h},\,L_{{\bf h}^*},\,d)$ as its characteristic. Furthermore,
\begin{equation*}
   \mb{dim}(\h)= \mb{rk}(A) + d_{\bf h} + d_{{\bf h}^*}  +d,
\end{equation*}
where $d_{\bf h} $ is the codimension of $ L_{\bf h}$ in $\mb{ker}(A^T)$, and $d_{{\bf h}^*}$ is the codimension of $ L_{{\bf h}^*}$ in $\mb{ker}(A)$.
\end{thm}

A realization $(\h,(h_i)_{i\in I},(\al_i)_{i\in I})$ of $A$ is called minimal if $d=0$. A realization $(\h,(h_i)_{i\in I},(\al_i)_{i\in I})$ of $A$ is called free if $h_1,\,\ldots,\, h_n\in \h $ and $\al_1,\,\ldots,\,\al_n \in \h^*$ are linearly independent. Equivalently,  $L_{\bf h}=\{0\}$ and $ L_{{\bf h}^*}=\{0\}$.

If $(\h,(h_i)_{i\in I},(\al_i)_{i\in I})$ is a realization of $A$, and $U$ is a subspace of $\R^n$ we set
\begin{equation*}
    [U]_{\bf h} := \{ r_1 h_1 + \cdots + r_n h_n   \mid r \in U  \}  \quad\mb{ and }\quad
    {[U]} _ {{\bf h}^*} := \{r_1 \al_1+\cdots+ r_n \al_n   \mid r \in U  \}  .
\end{equation*}
\begin{thm}\label{mainmr}  Let $(\h,(h_i)_{i\in I},(\al_i)_{i\in I})$ be a realization of $A$ with characteristic $(L_{\bf h},\,L_{{\bf h}^*},\,d)$, and  let $(\h' ,(h_i')_{i\in I},(\al_i')_{i\in I})$ be a realization of $A$ with characteristic $(L_{\bf h'},\,L_{{\bf h'}^*},\,d')$.
\begin{itemize}
\item[(a)] Suppose that $ L_{\bf h}\subseteq  L_{\bf h'}$ and $ L_{{\bf h}^*}=  L_{{\bf h'}^*} $ and  $d=d'$. Then there exists a surjective morphism of realizations $\phi: \h\to\h'$ with $\mb{ker}(\phi) =[\, L_{\bf h'} \,]_{\bf h}$. The map $\phi^*:{\h'}^*\to \h^*$ is injective, and $\mb{im}(\phi^*)=\left([\, L_{\bf h'} \,]_{\bf h} \right)^\bot$. Furthermore,
\begin{eqnarray}\label{qmorphmain1}
                    L_{\bf h'} / L_{\bf h}   &\to&  \qquad [ L_{\bf h'} \,]_{\bf h} \\
                            r    + L_{\bf h}              &\mapsto & r_1 h_1+\cdots + r_n h_n \nonumber
\end{eqnarray}
is a bijective linear map.
\item[(b)]  Suppose that $L_{\bf h} =  L_{\bf h'}$ and $L_{{\bf h}^*} \supseteq  L_{{\bf h'}^*}$ and $d=d'$. Then there exists an injective morphism  of realizations $\phi: \h\to\h'$ with $\mb{im}(\phi)=\left( [\, L_{{\bf h}^*} \,]_{{\bf h'}^*}\right)^\bot$. The map $\phi^*:{\h'}^*\to \h^*$ is surjective, and  $\mb{ker}(\phi^*) =[\, L_{{\bf h}^*} \,]_{{\bf h'}^*}  $.  Furthermore,
\begin{eqnarray}\label{qmorphmain2}
                    L_{{\bf h}^*} / L_{{\bf h'}^*} &\to&  \qquad [\, L_{{\bf h}^*} \,]_{{\bf h'}^*}  \\
                            r      +  L_{{\bf h'}^*}            &\mapsto & r_1 \al_1+\cdots + r_n \al_n \nonumber
\end{eqnarray}
is a bijective linear map.
\end{itemize}
\end{thm}
\Proof To (a): It is trivial to check that (\ref{qmorphmain1}) is a bijective linear map. The properties of $\phi^*$ follow from the properties of $\phi$, because for every linear map $\phi: \h\to\h'$ we have $\mb{im}(\phi^*)= \mb{ker}(\phi)^\bot $ and $\mb{ker}(\phi^*)= \mb{im}(\phi)^\bot$. 
By Theorem \ref{mainer} it is sufficient to find a realization $({\ti{\bf h}},(\ti{h}_i)_{i\in I},(\ti{\al}_i)_{i\in I})$ of $A$ with the same characteristic as $(\h',(h_i')_{i\in I},(\al_i')_{i\in I})$, and a surjective morphism of realizations $\phi: \h\to {\ti{\bf h}}$ with $\mb{ker}(\phi) ={[\, L_{ \ti{\bf h} }\,]}_{\bf h}=[\, L_{\bf h'} \,]_{\bf h}  $.

To make the proof easier to read we set $ U:=  L_{\bf h'}$ and $[U]:= [ L_{\bf h'} \,]_{\bf h}$. We define $\ti{\bf h}  :=  \h/[U]$, and $\ti{h}_i:=h_i+[U]$, $i\in I$. We identify ${\ti{\bf h}\mb{}}^* $ with $[U]^\bot \subseteq \h^*$. For $i\in I$ we have
\begin{equation*}
  \al_i(r_1 h_1+\cdots + r_n h_n) = r_1 a_{1i} + \cdots + r_n a_{ni}=0 
\end{equation*}
for all $r \in U$, since $U\subseteq\mb{ker}(A^T)$. Hence $\ti{\al}_i:=\al_i \in  [U]^\bot$, $ i\in I$. In this way we get a realization of $A$ because 
\begin{equation*}
  \ti{\al}_i(   \ti{h}_j)  =  \al_i( h_j+[U])=\al_i(h_j)=a_{ji} \mb{ for all }\;i,\,j\in I.
\end{equation*}

The canonical projection $\phi:\h\to \ti{\bf h}=\h/[U]$ is a surjective linear map with $\mb{ker}(\phi)=[U]= [ L_{\bf h'} \,]_{\bf h}$. It is also a morphism of realizations: Let $i\in I$. By definition, $\phi(h_i)=\ti{h}_i$. Furthermore,
\begin{equation*}
   \phi^*(\ti{\al}_i)(h)=\ti{\al}_i(\phi(h)) = \al_i(h+[U])=\al_i(h)
\end{equation*}
for all $h\in\h$. Thus $\phi^*(\ti{\al}_i)=\al_i$.

It remains to compute the characteristic of the realization $(\ti{\bf h}, (\ti{h}_i)_{i\in I}, (\ti{\al}_i)_{i\in I})$. Trivially, we have $ L_{  { \ti{\bf h}\mb{}}^* } =  L_{ { \bf h  }^* } =  L_{ { {\bf h}'  \mb{} }^* }$. Since $L_{\bf h}\subseteq  L_{\bf h'}=U$, we find
\begin{equation*}
    L_{\ti{\bf h}} = \{ r  \in \F^n\mid r_1 \ti{h}_1 + \cdots + r_n \ti{h}_n=0 \} 
    = \{ r \in \F^n \mid r_1 h_1 +\cdots + r_n h_n  \in [U] \}  
    = U+ L_{\bf h} = U = L_{\bf h'} .
\end{equation*}
Consider the surjective linear map
\begin{equation*}
  \gamma :\h\to \ti{\h} / \mb{span}(\ti{h}_1,\,\ldots,\,\ti{h}_n)
\end{equation*}
obtained by concatenation of the canonical projections, i.e., $\gamma (h) := \phi(h)+\mb{span}(\ti{h}_1,\,\ldots,\,\ti{h}_n)$, $h\in\h$.
We have $h\in\mb{ker}(\gamma)$ if and only if 
\begin{equation*}
   \phi(h)\in \mb{span}(\ti{h}_1,\,\ldots,\,\ti{h}_n) =\phi(\mb{span}(h_1,\,\ldots,\,h_n)) ,
\end{equation*}
if and only if
\begin{equation*}
   h\in \mb{span}(h_1,\,\ldots,\,h_n) +[U] = \mb{span}(h_1,\,\ldots,\,h_n).
\end{equation*}
This shows that we have isomorphic linear spaces
\begin{equation*}
     \h /  \mb{span}(h_1,\,\ldots,\,h_n)  \cong \ti{\h} / \mb{span}(\ti{h}_1,\,\ldots,\,\ti{h}_n).
\end{equation*}
Now we compute the defect:
\begin{eqnarray*}
\ti{d}   &=&  \mb{dim}( \ti{\bf h} /  \mb{span}(\ti{h}_1,\,\ldots,\,\ti{h}_n))- \mb{dim}(\mb{span}(\ti{\al}_1,\,\ldots,\,\ti{\al}_n)) + \mb{rk}(A) \\
          &=& \mb{dim}(\h / \mb{span}(h_1,\,\ldots,\,h_n))  - \mb{dim}(\mb{span}(\al_1,\,\ldots,\,\al_n)) + \mb{rk}(A) = d =d'.
\end{eqnarray*}

To (b): We may interprete $(\h^*,(\al_i)_{i\in I},(h_i)_{i\in I})$ and $(\h'^*,(\al'_i)_{i\in I},(h'_i)_{i\in I})$ as realization of $A^T$ with characteristics $(L_{{\bf h}^*}, L_{\bf h},\,\,d)$ and $(L_{{\bf h'}^*}, L_{\bf h'},\,\,d')$, and apply part (a).
\qed

In Chapter 3.1 of \cite{Sl} it is stated that an arbitrary realization is a subquotient of a free realization. By combining (a) and (b) of the preceding theorem, using also Theorem \ref{mainer}, we can now describe an arbitrary realization explicitly as a subspace of a quotient space of a free realization, and as a quotient space of a subspace of a free realization. We only formulate the first possibility. 
\begin{cor}\label{subquot} Let $(\h,(h_i)_{i\in I},(\al_i)_{i\in I})$ be a realization of $A$ with characteristic $(L_{\bf h},\,L_{{\bf h}^*},\,d)$. Then there exist a realization $(\h' ,(h_i')_{i\in I},(\al_i')_{i\in I})$ of $A$ with characteristic 
\begin{equation*} 
    (L_{\bf h'},\,L_{{\bf h'}^*},\,d')=(L_{\bf h},\,\{0\},\,d),
\end{equation*} 
a free realization $(\h'', (h_i'')_{i\in I},(\al_i'')_{i\in I})$ of $A$ with defect $d$, i.e., 
\begin{equation*}
     (L_{\bf h''},\,L_{{\bf h''}^*},\,d'')=(\{0\},\{0\},\,d),
\end{equation*} 
and morphisms of realizations
\begin{equation*}
\xymatrix{ \h''\ar@{->>}[rd]_<<<<<<\phi
&
&\h\ar@{_{(}->}[ld] ^<<<<<<<\psi\\
&\h'}
\end{equation*}
such that  the following hold:
\begin{itemize}
\item[(a)] $\phi: \h''\to\h'$ is surjective with $\mb{ker}(\phi) =[\, L_{\bf h} \,]_{\bf h''}$. The map $\phi^*: {\h'}^*\to\h''^*$ is injective with $\mb{im}(\phi^*) =\left([\, L_{\bf h} \,]_{\bf h''}\right)^\bot$. Furthermore, the linear space $ [\, L_{\bf h} \,]_{\bf h''}$ is isomorphic to $ L_{\bf h}$.
\item[(b) ] $\psi: \h \to \h'$ is injective with  $\mb{im}(\psi) =\left([\, L_{{\bf h}^*} \,]_{{\bf h'}^*}\right)^\bot  $. The map  $\psi^*: {\h'}^* \to \h^*$ is surjective with  $\mb{ker}(\psi^*) = [\, L_{{\bf h}^*} \,]_{{\bf h'}^*} $. Furthermore, the linear space $ [\, L_{{\bf h}^*} \,]_{{\bf h'}^*} $ is isomorphic to $  L_{{\bf h}^*}$.
\end{itemize}
\end{cor}
%
%

%
\section{Linear Coxeter groups}
%

%
\subsection{Linear Coxeter systems and some results from the literature}
%
%
%
In this subsection we collect some results on linear Coxeter systems from the literature. Some of these results can only be found for special linear Coxeter systems, but are easily adapted to the general situation.

Let $\h$ be a finite dimensional real linear space. Every zero-based closed ray $H$ in $\h$ determines a closed half space in $\h^*$ by
\begin{equation*}
  H^{\geq 0}:=\Mklz{\la\in\h^*}{\la(h)\geq 0 \mb{ for all } h\in H}.
\end{equation*}
The corresponding open half space $H^{>0}$, and the boundary hyperplane by  $H^{=0}$ are described by 
\begin{equation*}
  H^{> 0}=\Mklz{\la\in\h^*}{\la(h) >0 \mb{ for all } h\in H\setminus\{0\}} \quad \mb{ and }\quad 
  H^{= 0}=\Mklz{\la\in\h^*}{\la(h) = 0 \mb{ for all } h\in H }.
\end{equation*}
Similarly, every zero-based closed ray $L$ in $\h^*$ determines a closed half space $L^{\geq 0}$, an open half space $L^{>0}$, and a hyperplane $L^{=0}$ in $\h$.

Let $H$, $L$ be zero-based closed rays in $\h$, $\h^*$ such that $L\setminus\{0\}\subseteq H^{>0}$. The reflection $\sigma$ on $\h^*$ at the hyperplane $H^{=0}$ along the line $L\cup-L$ is defined by
\begin{equation*}  
     \sigma\res{H^{=0}} := id_{H^{=0}} \quad\mb{ and }\quad \sigma\res{L\cup -L} := -id_{L\cup -L} .
\end{equation*}
To describe $\sigma$ by a formula we choose $h\in H$, $\al\in L$ such that $\al(h)=2$. Then 
\begin{equation*}
   \sigma(\la) = \la -  \la(h)\al, \,\quad \la\in\h^*.
\end{equation*}
Now $L\setminus\{0\}\subseteq H^{>0}$ is equivalent to $H\setminus\{0\}\subseteq L^{>0}$. In the same way we get a reflection on $\h$. It is dual to the reflection on $\h^*$, and we denote it by the same symbol.

Linear Coxeter systems of finite rank have been introduced and investigated by E. B. Vinberg in \cite{V}. In difference to Definition 2 in \cite{V} we put the Tits cone into the dual of a finite dimensional linear space, which is no restriction, since the double dual of a finite dimensional linear space identifies with the linear space itself. In this way the definitions, results, and the notation of the current article are already adapted to their usage in some subsequent articles.

\begin{defn} Let $\h$ be a finite dimensional real linear space. Let $n\in\N$ and set $I=\{1,\,\ldots,\,n\}$. A linear Coxeter system consists of a family $(\h,\,(\,H_i,\,L_i\,)_{i\in I})$ where $H_i$, $L_i$ are zero-based closed rays in $\h$, $\h^*$ with $L_i\setminus\{0\}\subseteq H_i^{>0}$, $i\in I$, such that the following properties hold:
\begin{itemize}
\item[(a)] The polyhedral cone $ \overline{C} := \bigcap_{i\in I} H_i^{\geq 0}$ spans $\h^*$. Its pairwise different  1-codimensional faces are $ \overline{C} \cap H_i^{=0}$, $i\in I$.
\item[(b)] Let $C := \bigcap_{i\in I} H_i^{> 0}$ and let $\We$ be the group generated by the reflections $\sigma_i$ on $\h^*$ associated to $H_i$ and $L_i$, $i\in I$. Then
\begin{equation}\label{ofcprop}
     \sigma C\cap C=\emptyset  \quad\mb{ for all }\quad \sigma\in\We \setminus\{1\} .
\end{equation}
\end{itemize}
%
%
$\overline{C}$ is called the (closed) fundamental chamber, $C$ the open fundamental chamber. (Note also that $C$ is the interior of $\overline{C}$, and $\overline{C}$ is the closure of $C$.) The group $\We$ is called a linear Coxeter group. (Note also that $\We$ is a discrete subgroup of $GL(\h^*)$.) The reflections $\sigma_i$, $i\in I$, are called the simple reflections.
\end{defn}

In addition, we introduce morphisms:
\begin{prop}\label{mor lin Cox syst} Let $(\h,\,(H_i, \,L_i)_{i\in I})$ and $(\h',\,(H'_i, \,L_i')_{i\in I})$ be linear Coxeter systems. For a linear map $\phi:\h\to\h'$ the following are equivalent:
\begin{itemize}
\item[(i)] $\phi(H_i)= H'_i $ and $\phi^*(L_i')=L_i$ for all $i\in I$. 
\item[(ii)] $\phi^*\circ\sigma_i'=\sigma_i \circ \phi^* $ for all $i\in I$, and $(\phi^*)^{-1}(\overline{C})=\overline{C}'$.
\end{itemize}
A linear map $\phi:\h\to\h'$ with these properties is called a morphism of linear Coxeter systems.
\end{prop}
\Proof We only show that (ii) implies (i); the reverse implication is easy to see. Let $i\in I$. Choose $h_i\in H_i$, $\al_i\in L_i$ such that $\al_i(h_i)=2$, and $h_i'\in H_i$, $\al_i'\in L_i$ such that $\al_i'(h_i')=2$. The equation $\phi^*\circ\sigma_i'=\sigma_i \circ \phi^* $ is equivalent to 
\begin{equation}\label{morlC}
   \la' (h_i') \phi^*(\al_i') = \la' ( \phi(h_i) )   \al_i \quad \mb{ for all }\quad \la' \in (\h' )^*.
\end{equation}

Suppose that $\phi^*(\al_i')=0$. Then $-\al_i'\in (\phi^*)^{-1}(0)\subseteq  (\phi^*)^{-1}(\overline{C}) =\overline{C}'\subseteq (H_i')^{\geq 0}$, which is not possible.
Suppose that $\phi(h_i)=0$. Inserting $\la'=\al_i'$ in (\ref{morlC}) we get $2\phi^*(\al_i')=0$, which is not possible.

It remains $\phi^*(\al_i')\neq 0$ and $\phi(h_i)\neq 0$. Inserting $\la'=\al_i'$ in (\ref{morlC}) we find that $\phi^*(\al_i')= c \al_i$ with $c\in\R\setminus\{0\}$. Therefore, (\ref{morlC}) is equivalent to $c \la' (h_i') = \la' ( \phi(h_i) ) $ for all $\la' \in (\h' )^*$, which in turn is equivalent to $c  h_i' =  \phi(h_i) $.
Now we choose $\rho'\in C'\subseteq \overline{C}'=(\phi^*)^{-1}(\overline{C})$. We have $0< \rho'(h_i') = \frac{1}{c}\rho'(\phi(h_i)) = \frac{1}{c}\phi^*(\rho')(h_i)$ and $\phi^*(\rho')(h_i)\geq 0$. Thus $c\in\R^+$. We obtain $\phi^*(L_i')=\R^+_0\phi^*(\al_i')= \R^+_0 \al_i=L_i$, and $H_i'=\R^+_0 h_i' = \R^+_0 \phi(h_i) =\phi(H_i)$.
\qed

We now fix a linear Coxeter system $(\h,(\,H_i,\,L_i\,)_{i\in I})$. Its Coxeter group $\We$ acts dually on $\h$. We define the set of reduced coroots as
\begin{equation}\label{defredcoroots}
  \Phi_r^\vee : = \Mklz{\sigma H_i }{ \sigma\in\We,\,i\in I}.
\end{equation}
We define the sets of positive and negative reduced coroots as
\begin{equation}
  (\Phi_r^\vee)^+  : = \Mklz{H\in\Phi_r^\vee}{ \overline{C}\subseteq H^{\geq 0}} \quad\mb{ and }\quad
   (\Phi_r^\vee)^- : = \Mklz{H\in\Phi_r^\vee}{ \overline{C}\subseteq - H^{\geq 0} }
\end{equation}
and call $H_i$, $i\in I$, the simple positive reduced coroots. 

\begin{lem}\label{rrdu} We have $\Phi_r^\vee = (\Phi_r^\vee)^+ \,\dot{\cup}\,(\Phi_r^\vee)^- $.
\end{lem}
\Proof Let $H\in \Phi_r^\vee$. The $\We$-stabilizer of every element of $C$ is trivial by (\ref{ofcprop}), and the $\We$-stabilizer of every element of $H^{=0}$ is nontrivial. Thus $C\cap  H^{=0} =\emptyset$. Since $C$ is convex, it has to be contained either in $H^{\geq 0}$ or in $-H^{\geq 0}$. Therefore, its closure $\overline{C}$ is contained either in $H^{\geq 0}$ or in $-H^{\geq 0}$.
\qed

\begin{lem}\label{siperr} Let $i\in I$. Then $\sigma_i \left( (\Phi_r^\vee)^+\setminus\{H_i\} \right)= (\Phi_r^\vee)^+\setminus\{H_i\}$.
\end{lem}
\Proof Since $\sigma_i$ is involutive, it is sufficient to show the inclusion ``$\subseteq$". Suppose there exists a coroot $H\in (\Phi_r^\vee)^+\setminus\{H_i\}$ such that $\sigma_i H = H_i$. Then $H=\sigma_i H_i= -H_i\in (\Phi_r^\vee)^-$, which is not possible by Lemma \ref{rrdu}.
Suppose there exists a coroot $H\in (\Phi_r^\vee)^+\setminus\{H_i\}$ such that $\sigma_i H \in (\Phi_r^\vee)^-$. By definition, $\overline{C}\subseteq H^{\geq 0}$ and $\overline{C}\subseteq -(\sigma_i H)^{\geq 0}= -\sigma_i H^{\geq 0}$. We find
\begin{equation*}
  \overline{C}\cap H_i^{=0}\subseteq \overline{C}\cap\sigma_i\overline{C}\subseteq H^{\geq 0}\cap (-H^{\geq 0})=H^{=0}.
\end{equation*}
Therefore, the hyperplane $H^{=0}$ contains the linear hull of the 1-codimensional face $\overline{C}\cap H_i^{=0}$ of $\overline{C}$, which is $H_i^{=0}$. Hence $H^{=0}=H_i^{=0}$. From $\overline{C}\subseteq H^{\geq 0}$, $H_i^{\geq 0}$ it follows that $H^{\geq 0}=H_i^{\geq 0}$. Thus $H=H_i$, which is not possible.\qed
%
%

The following theorem has been given as part 6) of Theorem 2 in \cite{V}.
\begin{thm}
$(\We,\, (\sigma_i)_{i\in I})$ is a Coxeter system.
\end{thm}

The next theorem, except the fundamental domain property and formula (\ref{chTcrf1}), has been given as parts 1) and 5) of Theorem 2 in \cite{V}. The fundamental domain property follows immediately from the definition of the universal space and Proposition 8  in \cite{V}.  Formula (\ref{chTcrf1}) can be proved with Lemmas \ref{rrdu} and \ref{siperr} in the same way as Proposition 3.12 c) in \cite{K}. 
\begin{thm} \label{fuprop} The set
\begin{equation*}
  X:= \bigcup_{\sigma\in {\mathcal W}} \sigma \overline{C}
\end{equation*}
is a convex cone, which is called the Tits cone. It can also be described as 
\begin{equation}
   X  = \Mklz{\la\in\h^*}{\la\in H^{\geq 0} \mb{ for all up to finitely many } H\in (\Phi_r^\vee)^+} . \label{chTcrf1}
\end{equation}

The fundamental chamber $\overline{C}$ is a fundamental domain for the action of $\We$ on $X$, i.e., every $\We$-orbit meets $\overline{C}$ in a unique point. For $\la\in\overline{C}$ the isotropy group
\begin{equation*}
   \mb{stab}_{\mathcal W}(\la):=\Mklz{\sigma\in\We}{\sigma \la=\la}
\end{equation*}
is generated by the simple reflections it contains, i.e., by $\Mklz{\sigma_i}{\la\in H_i^{=0},\,i\in I}$. 
\end{thm}

We define the set of reflections of the linear Coxeter system as $T:=\Mklz{\sigma\sigma_i\sigma^{-1}}{\sigma\in\We,\,i\in I}$. The Coxeter group $\We$ acts by conjugation on $T$. The following proposition can be easily obtained by using the description of the isotropy group given in Theorem \ref{fuprop}.
\begin{cor}\label{redrefl} We obtain a $\We$-equivariant surjective map by
\begin{eqnarray*}
   \Phi^\vee_r \;\;&\to& \;\;\, T \\
    \sigma H_i \;&\mapsto & \sigma \sigma_i\sigma^{-1}
\end{eqnarray*}
where $\sigma\in\We$, $i\in I$.
Denote the image of $H\in  \Phi^\vee_r $ by $\sigma_H$. Then the fibre of $\sigma_H$ consists of $\pm H$. In particular, the map restricts to a bijective map on $(\Phi^\vee_r)^+$.
\end{cor}

Most statements of the next theorem have been given as parts 2), 3), and 4), or follow from part 5) of Theorem 2 in \cite{V}. Formula (\ref{chTcrf2}) generalizes Exercise 3.15 of \cite{K}. Proposition 2.69 of \cite{AB}, and Corollary \ref{redrefl} can be used for its proof.
The proof that the action of $\mathcal W$ on $X^\circ$ is even proper can be reduced to the following fact, see also Proposition 12.9 of \cite{Lee}: Let $U_\lambda, U_\mu$ be neighborhoods of $\lambda,\mu \in  X^\circ \cap\overline{C} $ obtained by Proposition 10 of \cite{V} from special neighborhoods constructed in \S1 of \cite{V}. If $\sigma\in\We$ such that $U_\lambda\cap \sigma U_\mu\neq \emptyset$ then $\sigma\in  \mb{stab}_{\mathcal W}(\la) \,\mb{stab}_{\mathcal W}(\mu)$. 
\begin{thm}\label{dintact}
The interior $X^\circ$ of $X$ is obtained by
\begin{equation*}
  X^\circ= \bigcup_{\sigma\in {\mathcal W}} \sigma (X^\circ \cap\overline{C}) \quad\mb{ with }\quad 
   X^\circ \cap\overline{C} = \Mklz{\la\in\overline{C}}{\mb{stab}_{\cal W}(\la) \mb{ is finite }}.
\end{equation*}
It can also be described as
\begin{equation}
    X^\circ =\Mklz{\la\in\h^*}{\la\in H^{>0} \mb{ for all up to finitely many } H\in (\Phi_r^\vee)^+} .\label{chTcrf2}
\end{equation}
The set of hyperplanes $\{H^{=0}\mid H\in (\Phi^\vee_r)^+\}$ is locally finite in $X^\circ$, and each of its hyperplanes meets $X^\circ$.

The Coxeter group $W$ acts as discrete group properly on $X^\circ$. The canonical map from $X^\circ \cap\overline{C} $ to the quotient space $X^\circ/\We$ is a homeomorphism.
\end{thm}
%
%

The following theorem can be easily obtained from Lemmas \ref{rrdu}, \ref{siperr}, and Corollary \ref{redrefl}. 
\begin{thm}\label{rlhalf} Let $\sigma\in\We$ and $i\in I$. Then we have:
\begin{itemize}
\item[(a)] If $l(\sigma\sigma_i)>l(\sigma)$ then $\sigma H_i\in (\Phi_r^\vee)^+$.
\item[(b)] If $l(\sigma\sigma_i)< l(\sigma)$ then $\sigma H_i\in (\Phi_r^\vee)^-$.
\end{itemize}
\end{thm}

We now explain how linear Coxeter systems and their morphisms can be described by certain realizations of matrices and their morphisms.
\begin{defn}
A matrix $A=(a_{ij})_{i,j\in I}$ with entries in $\R$ is called a generalized Cartan matrix if the following hold:
\begin{itemize}
\item[(a)] For all $i,\,j\in I$, $i\neq j$, the entries $a_{ij}$, $a_{ji}$ are either both strict negative or both zero.
\item[(b)] For all $i\in I$ we have $ a_{ii}=2 $. For all $i,\,j\in I$, $i\neq j$, such that  $a_{ij} a_{ji}< 4$ we have
\begin{equation*}
     a_{ij}a_{ji} = 4 \cos^2(\frac{\pi}{m_{ij} }  ) \quad\mb{ for some }\quad m_{ij}\in\N\setminus\{1\}.
\end{equation*}
\end{itemize}
\end{defn}
Note that in (b) the positive integers $m_{ij}$ are uniquely determined. For all $i,\,j\in I$, $i\neq j$, with $a_{ij}a_{ji}\geq 4$ we set $m_{ij}=\infty$.
For all $i\in I$ we set $m_{ii}=1$. 

For the classification of the indecomposable generalized Cartan matrices into matrices of finite, affine, or indefinite type we refer to Theorem 3 of \cite{V}, and to Theorem 4.3 and Corollary 4.3 of \cite{K}.

Part (a) of the next theorem has been shown in Proposition 17, and in Theorem 1, Proposition 6 and 7 of \cite{V}. Part (b) has been given as Proposition 13 and formula (24) of \cite{V}.
\begin{thm}\label{realCM1} Let $(\h,\,(H_i, \,L_i)_{i\in I})$ be a linear Coxeter system. For every $i\in I$ choose elements $h_i\in H_i$ and $\al_i\in L_i$ such that $\al_i(h_i)=2$. Then we have:
\begin{itemize}
\item[(a)] The matrix $ A:=(\al_j(h_i))_{i,j\in I} $ is a generalized Cartan matrix. Its associated matrix $(m_{ij})_{i,j\in I}$ is the Coxeter matrix of the Coxeter system $(\We$, $(\sigma_i)_{i\in I})$.
\item[(b)]  $(\h,(h_i)_{i\in I},(\al_i)_{i\in I})$ is a realization of $A$ such that $L_{\bf h}\cap (\R^+_0)^n=\{0\}$.
\end{itemize}
\end{thm}

Note that $L_{\bf h}\cap (\R^+_0)^n=\{0\}$ if and only if $h_1,\ldots, h_n$ are positively independent. 
Note also that for every $i\in I$ the element $h_i$ is only determined up to a factor in $\R^+$. If we have $h_i'= d_i h_i$ and $\al_i'=\frac{1}{d_i}\al_i$ with $d_i\in\R^+ $ for all $i\in I$, then the corresponding generalized Cartan matrices are related by
\begin{equation*}
   A'= DA D^{-1}  \quad\mb{ where }\quad D= \mb{diag}(d_1,\,\ldots,\,d_n).
\end{equation*}
The classification of the indecomposable generalized Cartan matrices is invariant under such transformations.

\begin{prop}\label{realCM2} Let $\phi:\h\to\h'$ be a morphism of the linear Coxeter systems $(\h,\,(H_i, \,L_i)_{i\in I})$ and $(\h',\,(H'_i, \,L_i')_{i\in I})$. For every $i\in I$ choose elements $h_i\in H_i$ and $\al_i'\in L_i'$ such that $\al_i'(\phi(h_i))=2$. Then  $(\h,(h_i)_{i\in I},(\phi^*(\al'_i))_{i\in I})$ and  $(\h',(\phi(h_i))_{i\in I},(\al'_i)_{i\in I})$ are realizations of the same generalized Cartan matrix, and $\phi:\h\to\h'$ is a morphism of realizations. 

In particular, the Coxeter systems $(\We, \,(\sigma_i)_{i\in I})$ and  $(\We', \,(\sigma'_i)_{i\in I})$ are isomorphic. 
\end{prop}

Part (a) of the following theorem has been given as Theorem 5 and formula (24) of \cite{V}. Part (b) is trivial.
\begin{thm}\label{realCM3} Let $A$ be a generalized Cartan matrix.
\begin{itemize}
\item[(a)] Let $(\h,(h_i)_{i\in I},(\al_i)_{i\in I})$ be a realization of $A$ such that  $L_{\bf h}\cap (\R^+_0)^n=\{0\}$. For every $i\in I$ set $H_i := \R h_i$ and $ L_ i :=\R\al_i$. Then $(\h,\,(H_i,\,L_i\,)_{i\in I})$ is a linear Coxeter system.
\item[(b)] Let $(\h,(h_i)_{i\in I},(\al_i)_{i\in I})$ and  $(\h',(h_i')_{i\in I},(\al_i')_{i\in I})$ be realizations of $A$ such that $L_{\bf h}\cap (\R^+_0)^n=\{0\}$ and $L_{\bf h'}\cap (\R^+_0)^n=\{0\}$.  Then every morphism $\phi:\h\to\h'$ of realizations is a morphism of the corresponding linear Coxeter systems. 
\end{itemize}
\end{thm}

Theorem \ref{realCM1}, Proposition \ref{realCM2}, and Theorem \ref{realCM3} show that for the investigation of linear Coxeter system and their morphisms we may restrict to realizations of a fixed generalized Cartan matrix and their morphisms. In this context the name ``root base" is more common than the name ``realization":

\begin{defn} A realization $(\h,(h_i)_{i\in I},(\al_i)_{i\in I})$ of a generalized Cartan matrix $A$ such that  $L_{\bf h}\cap (\R^+_0)^n=\{0\}$ is called a root base of a linear Coxeter system. A morphism of root bases is a morphism of realizations.
\end{defn}

We now fix a root base $(\h,(h_i)_{i\in I},(\al_i)_{i\in I})$ of a linear Coxeter system $(\h,\,(H_i,\,L_i\,)_{i\in I})$. The next theorem, which enhances Theorem \ref{rlhalf}, has been given as Theorem 1.10 in \cite{HN}.
\begin{thm}\label{ThHN}
Let $\sigma\in\We$ and $i\in I$. 
\begin{itemize}
\item[(a)] If $l(\sigma\sigma_i)> l(\sigma)$ then $ \sigma h_i\in \sum_{j\in I}\R^+_0 h_j$ and $\sigma\al_i\in \sum_{j\in I}\R^+_0 \al_j $.
\item[(b)] If $l(\sigma\sigma_i)< l(\sigma)$ then $\sigma h_i\in -\sum_{j\in I}\R^+_0 h_j $ and $\sigma\al_i\in -\sum_{j\in I}\R^+_0 \al_j$ .
\end{itemize}
\end{thm}

The following discussion shows that for this article it is not favorable to work with the coroots and roots of the root base instead of the reduced coroots (\ref{defredcoroots}).
Define the set of coroots as 
\begin{equation}
  \Phi^\vee:= \We \Mklz{h_i}{i\in I} ,
\end{equation}
and the sets of positive and negative coroots as $ (\Phi^\vee)^\pm:=\Phi^\vee\cap\pm \sum_{i\in I}\R^+_0 h_i $.
We have the $\We$-equivariant surjective map 
\begin{eqnarray}\label{corredcor}
\begin{array}{ccc}
   \Phi^\vee &\to& \Phi^\vee_r\\
     h  &\mapsto &\R^+_0 h
\end{array}
\end{eqnarray}
which maps $(\Phi^\vee)^{\pm}$ onto  $(\Phi^\vee_r)^{\pm}$. It follows from Lemma \ref{rrdu} that $  \Phi^\vee =   (\Phi^\vee)^+ \dot{\cup}\,  (\Phi^\vee)^- $.
\begin{rem} The map (\ref{corredcor}) is bijective for the root bases with integral generalized Cartan matrices used Lie theory as in \cite{K}, and for the root bases with symmetric generalized Cartan matrices used in \cite{D}. In general, it is possible that there are infinitely many coroots which map to the same reduced coroot. An example has been given in \cite{F}, Section 3, following Definition 3.1. 
\end{rem}

The composition of the map (\ref{corredcor}) with the map of Corollary \ref{redrefl} can be described as follows. Define the set of roots as 
\begin{equation}
  \Phi:= \We \Mklz{\al_i}{i\in I}.
\end{equation} 
We find from Corollary \ref{redrefl}, see also Remark 2.2 and Definition 2.3 of \cite{HN}, that we get a $\We$-equivariant surjective map
\begin{eqnarray}\label{mcorrf}
\begin{array}{ccc}
  \Phi^\vee &\to& \Phi\\
   \sigma h_i  &\mapsto &\sigma \al_i
\end{array}
\end{eqnarray}
where $\sigma\in\We$, $i\in I$. We denote the image of $h\in\Phi^\vee$ by $\al_h$. We obtain $\al_h(h)=2$. We denote by $\sigma_h$ the reflection at $(\R^+_0h)^{=0}$ along $\R\al_h$. In particular, $\sigma_{h_i}=\sigma_i$ for all $i\in I$. We have 
\begin{equation*}
   \sigma_{\tau h}=\tau\sigma_h\tau^{-1} \quad \mb{ for all }\quad \tau\in\We,\; h\in\Phi^\vee .
\end{equation*}
We conclude that $\sigma_{\R^+_0 h}= \sigma_h$ for all $h\in \Phi^\vee$.

\begin{rem}\label{exaffkm} The map (\ref{mcorrf}) does not need to be bijective: Let $\h$ be a $\R$-linear space with base $h_1$, $h_2$. Let 
\begin{equation*}
   A :=\left(\begin{array}{cc}
  2 & -2 \\
   -2& 2
\end{array}\right), 
\end{equation*}
which is an affine integral generalized Cartan matrix. Define $\al_1,\,\al_2\in\h^*$ by $\al_i(h_j) :=a_{ji}$ for all $i,\,j = 1,\,2$. Then $\al_2=-\al_1$.

This root base is isomorphic to the root base over $\R$ of the derived Kac-Moody algebra $\g(A)'$, the Kac-Moody algebra $\g(A)$ of $A$ defined as in \cite{K}. Here, $\Phi^\vee$, which maps bijectively to $\Phi^\vee_r$, is given by the infinite set
\begin{equation*}
  \Phi^\vee = \Mklz{h_1+n(h_1+h_2),\,h_2+n(h_1+h_2)}{n\in\Z}.
\end{equation*}
However, $\Phi=\{\al_1,\,\al_2\}$ is finite. The map (\ref{mcorrf}) maps $h_i+n(h_1+h_2)$ to $\al_i$, $n\in\Z$, $i=1,\,2$.
\end{rem}

We now interpret $\h$ as the dual of $\h^*$. Then $(\h^*,(\al_i)_{i\in I},(h_i)_{i\in I})$ is a root base of a linear Coxeter system if and only if $L_{{\bf h}^*}\cap (\R^+_0)^n=\{0\}$, which may not be satisfied. If we set
\begin{equation}\label{indexexc}
   I_0:= \{i\in I\mid \al_i \in (\sum_{j\in I}\R^+_0 \al_j ) \cap ( -\sum_{j\in I}\R^+_0 \al_j ) \} \quad \mb{ and } \quad I_1:=I\setminus I_0,
\end{equation}
then from Proposition \ref{pind} we find $L_{{\bf h}^*}\cap (\R^+_0)^n=\{0\}$ if and only if $I_0=\emptyset$. Nevertheless, we define the dual fundamental chamber and the dual Tits cone by
\begin{equation*}
  \overline{C}^\vee:=\Mklz{h\in\h}{\al_i(h)\geq 0 \mb{ for all }i\in I} \quad \mb{ and }\quad X^\vee:=\We\overline{C}^\vee,
\end{equation*}
respectively. The following Proposition has been shown in Remark 1.6 of \cite{HN}.
\begin{prop}\label{dlinCox1} $\overline{C}^\vee-\overline{C}^\vee $ is a $\We$-invariant subspace of $\h$, on which $\We_{I_0}$ acts trivially. It contains the elements $h_i$, $i\in I_1$.
\end{prop}

The linear functions $\al_i$, $i\in I_1$, restrict to linear functions on $\overline{C}^\vee-\overline{C}^\vee $, which we denote by the same symbol. We interpret $\overline{C}^\vee-\overline{C}^\vee $ as the dual of $(\overline{C}^\vee-\overline{C}^\vee)^*$. The following Proposition has also been shown in Remark 1.6 of \cite{HN}.
\begin{prop}\label{dlinCox2} If $I_1\neq\emptyset$ then $((\overline{C}^\vee-\overline{C}^\vee)^* ,(\al_i)_{i\in I_1},(h_i)_{i\in I_1})$ is a root base of a linear Coxeter system. The action of its Coxeter group coincides with the action of $\We_{I_1}$ on $\overline{C}^\vee-\overline{C}^\vee$. Its fundamental chamber coincides with $\overline{C}^\vee$. Its Tits cone coincides with $X^\vee$.  
\end{prop}

We set $(\R_0^+)^{I_0} :=  \bigoplus_{i\in I_0}\R^+_0 e_i \subseteq \R^n$, where $e_1, e_2, \ldots, e_n$ is the canonical base of $\R^n$. For some variants of the following proposition see Section 6.1 of \cite{Kr}. Compare also with formula (24) of \cite{V}. 
\begin{prop}\label{dlinCox3} $I_0$ is either empty or a nonempty union of affine connected components of $I$, and we have
\begin{equation}\label{dlinCox3eq}
    L_{{\bf h}^*} \cap (\R^+_0)^n =\{ r\in (\R^+_0)^{I_0} \mid \sum_{i\in I_0} r_i\al_i= 0 \}.
\end{equation}

\end{prop}
\Proof We define the support of $r\in  (\R^+_0)^n$ as $ \mb{supp}(r):=\Mklz{i\in I}{r_i>0}$. If $r, \ti{r}\in L_{{\bf h}^*}\cap (\R^+_0)^n$ then 
\begin{eqnarray*}
   r+\ti{r}\in L_{{\bf h}^*}\cap (\R^+_0)^n \quad \mb{and }\quad supp(r+\ti{r})=\mb{supp}(r)\cup \mb{supp}(\ti{r}).
\end{eqnarray*} 
Since there are only finitely many subsets of $I$, there exists an element $s\in L_{{\bf h}^*}\cap (\R^+_0)^n$ with biggest support.
We show that $\mb{supp}(s)=I_0$, which implies (\ref{dlinCox3eq}). For $i\in \mb{supp}(s)$ we find from
\begin{eqnarray*}
   \al_i= -\frac{1}{s_i}\sum_{j\in\,supp(s),\,j\neq i} s_j\al_j 
\end{eqnarray*}
that $i\in I_0$. If $i\in I_0$ then there exists an element $r\in (\R_0^+)^n$ such that 
\begin{eqnarray*}
\al_i=-\sum_{j\in\,supp(r)}r_j \al_j .
\end{eqnarray*}
Because $\mb{supp}(s)$ is the biggest support, we conclude that $i\in \mb{supp}(r)\cup\{i\}\subseteq\mb{supp}(s)$.

Now let $I_0=\mb{supp}(s)\neq\emptyset$. By Proposition \ref{dlinCox1} we have $\al_j(h_i)=0$ for all $j \in I_0$ and $i\in I_1=I\setminus I_0$. Therefore, $I_0$ is a union of connected components of $I$. Because of $s\in L_{{\bf h}^*} \subseteq \mb{ker}(A)$ we get 
\begin{eqnarray*}
  \sum_{j\in\,supp(s)}a_{ij}s_j=0 \quad \mb{ for all } \quad  i \in \mb{supp}(s).
\end{eqnarray*}
From Corollary 4.3 of \cite{K} it follows that the connected components of $\mb{supp}(s)$ are of affine type. 
\qed

We will also need later:
\begin{prop}\label{fcsdiff} We have
\begin{eqnarray}
   \overline{C}^\vee =  \{h\in\h \mid h-\sigma h \in\sum_{i\in I}\R^+_0 h_i \mb{ for all }\sigma\in\We\},\label{fcsdiff1}\\
 \overline{C} \,\subseteq  \{\la\in\h^* \mid \la-\sigma\la \in\sum_{i\in I}\R^+_0 \al_i \mb{ for all }\sigma\in\We\}.\label{fcsdiff2}
\end{eqnarray}
\end{prop}
\begin{rem} The inclusion in (\ref{fcsdiff2}) can be proper. This can be seen by the example of Remark \ref{exaffkm}, where the set on the right in (\ref{fcsdiff2}) coincides with $\h^*$.
\end{rem}
\Proof The inclusions ``$\subseteq$" in (\ref{fcsdiff1}) and (\ref{fcsdiff2}) can be proved similarly to Proposition 3.12 d) in \cite{K}, if Theorem \ref{ThHN} is used. Now let $h$ be contained in the set on the right in (\ref{fcsdiff1}). For every $j\in I$ we have
\begin{equation*}
   \R h_j\ni\al_j (h) h_j = h - \sigma_j h \in \sum_{i\in I}\R^+_0 h_i.
\end{equation*}
We obtain from Proposition \ref{fcequiv} that $\R^+_0 h_j$ is a face of $\sum_{i\in I}\R^+_0 h_i$. Its linear hull is $\R h_j$. We conclude that $\al_j(h) h_j\in ( \sum_{i\in I}\R^+_0 h_i) \cap\R h_j=\R^+_0 h_j$, so $\al_j(h)\geq 0$. Hence $h\in \overline{C}^\vee$.
\qed

To state some results of the literature in an easy way we introduce the following notations. A root base $(\h,(h_i)_{i\in I},(\al_i)_{i\in I})$ of a linear Coxeter system with characteristic $(L_{\bf h},\,L_{{\bf h}^*},\,d)$ and generalized Cartan matrix $A$ is called
\begin{itemize}
\item free if $L_{\bf h}=\{0\}$ and $L_{{\bf h}^*}=\{0\}$,
\item minimal if $d=0$,
\item integral if the entries of $A$ are integers,
\item symmetric if $A$ is a symmetric matrix and, in addition to $L_{\bf h}\cap(\R^+_0)^I=\{0\}$, also $L_{{\bf h}^*}\cap (\R^+_0)^I=\{0\}$. 
\end{itemize}

The properties ``free" and ``minimal" are common properties of all root bases, the properties ``integral" and ``symmetric" are not common to all root bases of a linear Coxeter system. We call a linear Coxeter system free, minimal, integral, symmetric if there exists a root base of the linear Coxeter system which is free, minimal, integral, symmetric, respectively.

The based root systems in \cite{D}, for which the faces of the Tits cone and imaginary cone have been described, correspond to symmetric root bases with additional non-degenerate invariant symmetric bilinear forms.\\

Let $A=(a_{ij})_{i,j\in I}$ be a generalized Cartan matrix. To keep our notation simple we define that a subset $J\subseteq I$ has a certain property if the corresponding submatrix $A_J:=(a_{ij})_{i,j\in J}$ has this property.

For $J\subseteq I$ we denote by $J^0$ the union of the finite type components of $J$, by $J^{aff}$ the union of the affine type components of $J$, and by $J^{ind}$ the union of the indefinite type components of $J$. (An empty union is defined to be empty.) We denote by $J^\infty$ the union of the nonfinite type components of $J$, i.e., $J^\infty=J^{aff}\cup J^{ind}$. As in the article \cite{Loo} we call $J$ special if $J^\infty =J$.

We call $j,k\in I$ adjacent if $a_{jk}\neq 0$. More generally, we call $J,K\subseteq I$ adjacent if there exist $j\in J$ and $k\in K$ such that $a_{jk}\neq 0$. We call $J,K\subseteq I$ separated if they are not adjacent. For $J\subseteq I$ we set
\begin{equation*}
  J^\bot:=\Mklz{i\in I}{a_{ij}=0 \mb{ for all }j\in J}.
\end{equation*}

Let $G$ be a group acting on a set $M$. For $D\subseteq M$ we denote by
\begin{equation*}
  Z_G(D):=\{g\in G\mid gd=d \mb{ for all }d\in D\}
\end{equation*}
the pointwise stabilizer of $D$, and by
\begin{equation*}
  N_G(D):=\{g\in G\mid gD=D\}
\end{equation*}
the stabilizer of $D$ as a whole. $N_G(D)$ is a subgroup of $G$, which contains $Z_G(D)$ as a normal subgroup.
%
%
%
%
%
\subsection{\label{SubsectionFF} The faces of the fundamental chamber}
%
%
%
The fundamental chamber of a linear Coxeter system can be partitioned 
\begin{itemize}
\item according to the isotropy groups of its points,
\item by the relative interiors of its faces.
\end{itemize}
We show in Proposition \ref{fsFa} and \ref{rihfc} that both partitions coincide and are described by the facial sets. Key results in this subsection are Corollary \ref{critfa1}, Corollary \ref{detfacial} and its interpretation by a hyperplane arrangement in Remark \ref{hypfa}, which allow to compute the facial sets of linear Coxeter systems. An example is given in the next subsection. Many theoretical results on the facial sets follow in combination with Lemmas \ref{kerprop1} and \ref{kerprop2}, which describe properties of the elements of the kernels of generalized Cartan matrices.

Let $(\h,(h_i)_{i\in I},(\al_i)_{i\in I})$ be a root base of the linear Coxeter system $(\h,\,(H_i, \,L_i)_{i\in I})$. The fundamental chamber 
\begin{equation*}
   \overline{C} = \Mklz{ \la\in\h^* }{ \la(h_i)\geq 0 \mb{ for all } i\in I }
\end{equation*}
is the convex cone dual to the finitely generated convex cone $  \sum_{i\in I} \R^+_0 h_i$, for which the properties of Proposition \ref{fcequiv} (ii) hold. For $J\subseteq I$ we obtain
\begin{equation}\label{intfh}
   ( \sum_{j\in J} \R^+_0 h_j )\cap\{h_i \mid i\in I \} =\{ h_j \mid j\in J\} 
\end{equation}
from Proposition \ref{fcequiv} (iii). We call a set $J\subseteq I$ facial if the following equivalent conditions are satisfied:
\begin{itemize}
\item[(i)] There exists a face $F$ of $\sum_{i\in I}\R^+_0 h_i$ such that $J=\Mklz{i\in I}{h_i\in F}$.
\item[(ii)] $\sum_{j\in J}\R^+_0 h_j$ is a face of $\sum_{i\in I}\R^+_0 h_i$.
\end{itemize}
Clearly, this notion does not depend on the particular root base of the linear Coxeter system. 
It is easy to show the following proposition on the facial sets. For part (c) use also Proposition \ref{fcequiv} (ii).
\begin{prop}\label{fsFa0} 
\begin{itemize}
\item[(a)] An arbitrary intersection of facial sets is facial. In particular, the set $I$ is facial.
\item[(b)] Let $L\subseteq I$. There exists a smallest facial set $L^{fc}$ which contains $L$. It is given by
\begin{equation*}
   L^{fc}=\bigcap_{J\text{ facial},\:J\supseteq L} J .
\end{equation*}
We call $L^{fc}$ the facial closure of $L$.
\item[(c)] Ordered partially by inclusion, $\Mklz{J\subseteq I}{J \mb{ facial }}$ is a lattice, the lattice meet and join given by
\begin{equation*}
   J_1 \wedge J_2= J_1\cap J_2\quad \text{ and }\quad  J_1\vee J_2=(J_1\cup J_2)^{fc},\quad J_1,\,J_2\, \mb{ facial}.
\end{equation*}
The set $\emptyset$ is the smallest facial set, the sets $\{i\}$, $i\in I$, are the atoms, the set $I$ is the biggest facial set.
\end{itemize}
\end{prop}
In general, not all subsets of $I$ are facial.
\begin{prop}\label{facial-h-free} The following are equivalent:
\begin{itemize}
\item[(i)] Every subset of $I$ is facial.
\item[(ii)] $L_{\bf h}=\{0\}$, i.e., $h_1$, $h_2$, \ldots, $h_n$ are linearly independent.
\end{itemize}
\end{prop}
\Proof It is easy to check ``$(ii)\Rightarrow (i)$" directly; we only show ``$(i)\Rightarrow (ii)$": Let $r\in L_{\bf h}$. Then $I_+:=\Mklz{i\in I}{r_i>0}$ and $I_-:=\Mklz{i\in I}{r_i<0}$ are facial sets by (i). Hence $\sum_{i\in I_+}\R^+_0 h_i$ and $\sum_{i\in I_-}\R^+_0 h_i$ are faces such that
\begin{equation*}
   \mb{ri}(\sum_{i\in I_+}\R^+_0 h_i) \ni \sum_{i\in I_+}r_i h_i \,=\, \sum_{i\in I_-}(- r_i) h_i \in   \mb{ri}(\sum_{i\in I_-}\R^+_0 h_i).
\end{equation*}
We conclude that $\sum_{i\in I_+}\R^+_0 h_i =\sum_{i\in I_-}\R^+_0 h_i$. From (\ref{intfh}) we find $I_+=I_-$. This is only possible for $I_+=I_-=\emptyset$, so $r=0$.
\qed

For $J\subseteq I$ we set
\begin{eqnarray}
    \overline{F}_J  &:=&   \Mklz{\la\in\h^*}{\la(h_i)=0 \mb{ for } i\in J ,\, \la(h_i)\geq 0 \mb{ for } i\in I\setminus J },\label{openfacet}\\
     F_J &:=&  \Mklz{\la\in\h^*}{\la(h_i)=0 \mb{ for } i\in J ,\, \la(h_i) > 0 \mb{ for } i\in I\setminus J }.\label{closedfacet}
\end{eqnarray}
The lattice of facial sets is isomorphic to the face lattice of $\sum_{i\in I}\R^+_0 h_i$, which is anti-isomorphic to the face lattice of the dual cone $\overline{C}$ by Proposition I.2 (5) of \cite{N}. Thus we have shown:

\begin{prop}\label{fsFa} We obtain an anti-isomorphism of lattices by:
\begin{eqnarray*}
 \Mklz{J\subseteq I}{ J \mb{ facial }  } &\to        &    \Fa{\overline{C}}\\
                          J\qquad\quad          &\mapsto &    \quad \overline{F}_J
\end{eqnarray*}
\end{prop}

The relative interiors and linear hulls of the faces of $\overline{C}$, as well as the isotropy groups of the points of $\overline{C}$ can be described by the facial sets as follows.  
\begin{prop}\label{rihfc} If $J\subseteq I$ is facial then
\begin{equation}\label{ripfC}
    \overline{F}_J = \dot{\bigcup_{K\supseteq J,\;K\text{ facial} }} F_K,
\end{equation}
and $F_J$ is the relative interior of $\overline{F}_J$. In particular, $F_J\neq\emptyset$. The linear hull of $\overline{F}_J$ is given by
\begin{equation}\label{hullfC}
    \overline{F}_J - \overline{F}_J =\Mklz{\la\in\h^*}{\la(h_j)=0\mb{ for all }j\in J\, }.
\end{equation}
For every $\la\in F_J$ we have
\begin{equation}\label{WstabC}
  \mb{stab}_{\mathcal W}(\la) = \We_J.
\end{equation} 
\end{prop}
\Proof From Proposition \ref{fsFa} we find that the faces of $\overline{F}_J$ are given by $\overline{F}_K$, $K$ facial, $K\supseteq J$. It follows that
\begin{equation*}
  \mb{ri} (\overline{F}_J) = \overline{F}_J\setminus\bigcup_{K\supsetneqq J,\;K\text{ facial} }\overline{F}_K.
\end{equation*}
Trivially, we get $F_J\subseteq \mb{ri}(\overline{F}_J)$. If $\la\in\mb{ri}(\overline{F}_J)\subseteq \overline{F}_J$ then $(\sum_{i\in I}\R^+_0 h_i ) \cap (\R\la)^\bot$ is a face of $\sum_{i\in I}\R^+_0 h_i$ with corresponding facial set $ \Mklz{i\in I}{\la(h_i)=0} \supseteq J$. Since $\la$ is not contained in a proper face of $\overline{F}_J$, we conclude that $\Mklz{i\in I}{\la(h_i)=0} = J$. Hence $\la\in F_J$.

The relative interiors of the faces of $ \overline{F}_J $ give a partition of $\overline{F}_J $, which is equation (\ref{ripfC}).
The inclusion ``$\subseteq$" in (\ref{hullfC}) is trivial. To show the reverse inclusion we choose an element $\la_0\in F_J $. Let $\la\in\h^*$ such that $\la(h_j)=0$ for all $j\in J$. Then for some sufficiently large $r\in\R^+$ we have
\begin{equation*}
  (\la + r\la_0)(h_i)=\la(h_i) + r \underbrace{\la_0(h_i)}_{>0} >0
\end{equation*}
for all $i\in I\setminus J$. It follows that
\begin{equation*}
    \la=(\la + r\la_0) -r\la_0\in F_J -\R^+_0\la_0\subseteq   \overline{F}_J -  \overline{F}_J.  
\end{equation*}
Formula (\ref{WstabC}) is obtained from the description of the isotropy group in Theorem \ref{fuprop}.
\qed

Recall that $F_J$ and $\overline{F}_J$ have been defined in (\ref{openfacet}) and (\ref{closedfacet}) for arbitrary subsets $J\subseteq I$, which is sometimes useful. With Proposition \ref{rihfc} it is easy to show:
\begin{prop} \label{fshfree}
Let $J\subseteq I$.
\begin{itemize}
\item[(a)] We have $F_J\neq \emptyset$ if and only if $J$ is facial.
\item[(b)]  $\overline{F}_J = \overline{F}_{J^{fc}}$.
\end{itemize}
\end{prop}

Morphisms of root bases relate the corresponding fundamental chambers and their dual cones in the following way:
\begin{prop}\label{dcpsi} Let $(\h,(h_i)_{i\in I},(\al_i)_{i\in I})$ and  $(\h', (h_i')_{i\in I}, (\al_i')_{i\in I})$ be root bases of linear Coxeter systems. Let  $\phi:\h\to\h'$ be a morphism of root bases. 
\begin{itemize}
\item[(a)] We have $  \sum_{i\in I} \R^+_0 h_i '\subseteq  \mb{im}(\phi)$ and
\begin{equation}\label{dcpi}
   \phi^{-1}(  \sum_{i\in I} \R^+_0 h_i ')  =    \sum_{i\in I} \R^+_0 h_i + \mb{ker}(\phi).
\end{equation}
\item[(b)] We have $\mb{ker}(\phi^*) \subseteq\overline{C}' \cap (-\overline{C}' )$ and
\begin{equation}\label{dcps}
   \phi^*( \overline{C'})  =   \overline{C}\cap \mb{im}(\phi^*).
\end{equation}
\item[(c)] The convex cones (\ref{dcpi}) and (\ref{dcps}) are dual to each other.
\end{itemize}
\end{prop}
\Proof Part (a) follows from $  \phi( \sum_{i\in I} \R^+_0 h_i )= \sum_{i\in I} \R^+_0 \phi(h_i) = \sum_{i\in I} \R^+_0 h_i' $. Part (b) follows from 
\begin{equation*}
    (\phi^*)^{-1}(0) \subseteq (\phi^*)^{-1}(\overline{C}) = 
       \Mklz{\la'\in (\h')^*}{0\leq \phi^*(\la')(h_i)=\la'(\phi^*(h_i))=\la'(h_i')  \mb{ for all } i\in I \,}=\overline{C'}.
\end{equation*}
The closed convex cones $ \sum_{i\in I} \R^+_0 h_i $ and $ \overline{C}$ are dual. The orthogonal linear spaces $\mb{ker}(\phi)$ and $\mb{im}(\phi^*)$ are also dual closed convex cones. Now (c) follows from Proposition I.1.6 in \cite{HHL}. 
\qed

The descriptions in Proposition \ref{dcpsi} (a) and (b) allow to apply Lemmas \ref{conint}, \ref{consum}, and \ref{concap}. The next theorem shows how morphisms of root bases relate the corresponding facial sets. 
\begin{thm}\label{facialmorph} Let $(\h,(h_i)_{i\in I},(\al_i)_{i\in I})$ and  $(\h', (h_i')_{i\in I}, (\al_i')_{i\in I})$ be root bases of linear Coxeter systems. Let  $\phi:\h'\to\h$ be a morphism of root bases. For $J\subseteq I$ we have 
\begin{equation}\label{facmorph}
      \phi^{*}(F_J)  =   F_J'  \cap \mb{im}(\phi^*) \quad \mb{ and } \quad  \phi^{*}(\overline{F}_J)  =   \overline{F'}_J \cap{im}(\phi^*)   ,
\end{equation}
and the following are equivalent:
\begin{itemize}
\item[(i)] $J$ is facial for $(\h ,(h_i)_{i\in I}, (\al_i)_{i\in I})$.
\item[(ii)] $J$ is facial for $(\h',(h'_i)_{i\in I},(\al'_i)_{i\in I})$ and
\begin{equation}\label{fmkBed}
  \mb{ker}(\phi) \cap ( \,\sum_{j\in J } \R h'_j + \sum_{i\in I\setminus J} \R^+_0 h'_i   \,) \subseteq \sum_{j\in J}\R h'_j  .
\end{equation}
\end{itemize}
\end{thm}
\Proof It is easy to check (\ref{facmorph}); we only show the equivalence of (i) and (ii). From the description in Proposition \ref{dcpsi} (a) we get by Lemmas \ref{conint} and \ref{consum} bijections
\begin{equation*}
 \Fa{\sum_{i\in I}\R^+_0 h_i} \to \Fa{ \sum_{i\in I}\R^+_0 h_i'+\mb{ker}(\phi)} \to \{\,F'\in\Fa{\sum_{i\in I}\R^+_0 h_i'}\mid\mb{ker}(\phi)\cap (\sum_{i\in I}\R^+_0 h_i' -F')\subseteq F'-F'\,\},
\end{equation*}
where the maps from the left to the right are given by
\begin{equation*}
   G\mapsto\phi^{-1}(G) \mapsto \phi^{-1}(G)\cap\sum_{i\in I}\R^+_0 h_i'.
\end{equation*}
Furthermore, we have
\begin{equation*}
  \{\,i\in I \mid h_i'\in \phi^{-1}(G)\cap\sum_{j\in I}\R^+_0 h_j'\,\} =\Mklz{i\in I}{\phi(h_i')\in G} =\Mklz{i\in I}{h_i\in G} .
\end{equation*}
This shows that
\begin{equation} \label{bfmor1}
     \Mklz{i\in I}{h_i\in G} \mb{ is facial for } (\h',(h'_i)_{i\in I},(\al'_i)_{i\in I})
\end{equation} 
and 
\begin{equation}\label{bfmor2}
     \phi^{-1}(G)\cap\sum_{i\in I}\R^+_0 h_i '= \sum_{i\in I,\,h_i\in G}\R^+_0 h_i',
\end{equation}
and
\begin{equation}\label{bfmor3}
   \mb{ker}(\phi)\cap \Big(\sum_{i\in I}\R^+_0 h_i'\, - \sum_{i\in I,\,h_i\in G}\R^+_0 h_i' \Big)\subseteq \sum_{i\in I,\,h_i\in G}\R h_i'.
\end{equation}

Suppose that (i) holds. Then $G:=\sum_{j\in J}\R^+_0 h_j$ is a face of $\sum_{i\in I}\R^+_0 h_i$ such that $J=\Mklz{i\in I}{h_i\in G}$. From (\ref{bfmor1}) and (\ref{bfmor3}) we get (ii).

Suppose that (ii) holds. Then there exists a uniquely determined face $G$ of $\sum_{i\in I}\R^+_0 h_i$ such that
\begin{equation*}
       \phi^{-1}(G)\cap\sum_{i\in I}\R^+_0 h_i ' = \sum_{i\in J}\R^+_0 h_i'.
\end{equation*}
From (\ref{bfmor2}) and (\ref{intfh}) we find that $J= \Mklz{i\in I}{h_i\in G}$, which is facial for  $(\h ,(h_i)_{i\in I}, (\al_i)_{i\in I})$.
\qed

We know the facial sets in the particular case described in Proposition \ref{fshfree}. We now use the preceding theorem to obtain a criterion which allows to find the facial sets in the general case. Let $e_1, e_2, \ldots, e_n$ be the canonical base of $\R^n$. For $J\subseteq I$ we set
\begin{equation*}
    \R^J  :=  \bigoplus_{j\in J}\R e_j\subseteq \R^n    \quad \mb{ and }\quad     (\R_0^+)^J :=  \bigoplus_{j\in J}\R^+_0 e_j \subseteq \R^n.
\end{equation*}
\begin{cor} \label{critfa1} Let $(\h,(h_i)_{i\in I},(\al_i)_{i\in I})$ be a root base of a linear Coxeter system with characteristic  $(L_{\bf h},\,L_{{\bf h}^*},\,d)$. For $J\subseteq I$ the following are equivalent:
\begin{itemize}
\item[(i)] $J$ is facial.
\item[(ii)]  $L_{\bf h}\cap \left( \R^J + (\R^+_0)^{I\setminus J} \right)  \subseteq   \R^J $.
\end{itemize}
\end{cor}
\Proof By Theorem \ref{mainer} and \ref{mainmr} (a) there exists a root base $(\h',(h_i')_{i\in I},(\al_i')_{i\in I})$ of a linear Coxeter system with characteristic $ (L_{\bf h'},\,L_{{\bf h'}^*},\,d') =  (\{0\},\,L_{{\bf h}^*},\,d)$, and a surjective morphism $\phi:\h'\to\h$ of root bases with 
\begin{equation}\label{kphifree}
  \mb{ker}(\phi)= [L_{\bf h}]_{\bf h'} =  \left\{r_1 h_1' + \cdots+  r_n h_n'   \left|  r \in L_{\bf h}  \right.\right\}.  
\end{equation}

Since $ L_{\bf h'} =  \{0\}$, we get from Proposition \ref{fshfree} that every subset of $I$ is facial for $(\h',(h_i')_{i\in I},(\al_i')_{i\in I})$. Hence we obtain from Theorem \ref{facialmorph}: $J$ is facial for $(\h,(h_i)_{i\in I},(\al_i)_{i\in I})$ if and only if
\begin{equation*}
  \mb{ker}(\phi) \cap (\, \sum_{j\in J } \R h'_j + \sum_{i\in I\setminus J} \R^+_0 h'_i   \,) \subseteq \sum_{j\in J}\R h'_j  ,
\end{equation*}
where for $\mb{ker}(\phi)$ the expression on the right in (\ref{kphifree}) has to be inserted. Since $h_1'$, \ldots, $h_n'$ are linearly independent, this is equivalent to (ii). 
\qed

It is useful to formulate the criterion of the preceding corollary differently, which can be motivated as follows: The elements $h_1,\,\ldots,\, h_n$ of a root base $(\h,(h_i)_{i\in I},(\al_i)_{i\in I})$ of a linear Coxeter system are only determined up to positive real numbers. Therefore, only the sign vectors of the elements of 
\begin{equation*}
   L_{\bf h}=\Mklz{r\in\R^n}{r_1h_1+\cdots+r_n h_n=0},
\end{equation*}
are relevant for the linear Coxeter system.
For $r\in \R^n$ we set
\begin{equation*}
  I_+(r) := \Mklz{ i\in I }{ r_i > 0 },  \quad  I_0(r) := \Mklz{ i\in I} { r_i = 0 }, \quad  I_-(r) := \Mklz{ i\in I }{ r_i < 0 }.
\end{equation*}
Clearly, we have $I_\pm(-r)=I_\mp(r)$ and $I_0(-r)=I_0(r)$.

\begin{cor}\label{detfacial} Let $(\h,(h_i)_{i\in I},(\al_i)_{i\in I})$ be a root base of a linear Coxeter system with characteristic  $(L_{\bf h},\,L_{{\bf h}^*},\,d)$. For $J\subseteq I$ the following are equivalent:
\begin{itemize}
\item[(i)] $J$ is facial.
\item[(ii)] For every element $r\in L_{\bf h}$ either $ I_+(r),\,I_-(r)\subseteq J $ or $ I_+(r),\, I_-(r) \not\subseteq  J$. 
\item[(iii)] For all  $r\in L_{\bf h}$ the following holds: If $ I_+(r)\subseteq J $ then also $I_-(r)\subseteq J$.
\item[(iv)] For all  $r\in L_{\bf h}$ the following holds: If $ I_-(r)\subseteq J $ then also $I_+(r)\subseteq J$.
\end{itemize}
\end{cor}
\begin{rem}\label{remdetfacial} We have $I_+(0)=I_-(0)=\emptyset$. Hence we can replace $L_{\bf h}$ by $L_{\bf h}\setminus\{0\}$ in (ii), (iii), (iv).
\end{rem}
\Proof By Corollary \ref{critfa1} the set $J$ is facial if and only if $  \Mklz{r\in L_{\bf h}}{I_-(r)\subseteq J}\subseteq \R^J $, which is equivalent to (iv).   $L_{\bf h}$ is a linear space, and for $r\in L_{\bf h}$ we have $I_\pm(-r)=I_\mp(r)$. Hence (iv) is equivalent to (iii). From (iii) and (iv) we get (ii). Trivially, from (ii) we get (iii).
\qed
\begin{rem}\label{hypfa} We can interpret the corollary in a more geometric way: Let $L_{\bf h}$ be k-dimensional and
\begin{equation*}
  b_1=\left(\begin{array}{c}b_{11}  \\
                                         \vdots       \\
                                         b_{1n} 
\end{array}\right)  ,\;
  b_2=\left(\begin{array}{c}b_{21}  \\
                                         \vdots       \\
                                         b_{2n} 
\end{array}\right) ,
\;\cdots,\;
    b_k=\left(\begin{array}{c}b_{k1}  \\
                                         \vdots       \\
                                         b_{kn} 
\end{array}\right) 
\end{equation*}
a base of $L_{\bf h}$. Every $r\in L_{\bf h}$ can be written in the form
\begin{equation*}
   r = s_1 b_1 + s_2 b_2+\cdots s_k b_k 
      = \left(\begin{array}{c} s_1 b_{11}+ s_2 b_{21} + \cdots + s_k b_{k1}  \\
                                         \vdots       \\
                                        s_1 b_{1n} +s_2 b_{2n} +\cdots + s_k b_{kn}
\end{array}\right)  
    \quad \mb{ with } \quad s_1,\,s_2,\,\ldots,\,s_k\in\R.
\end{equation*}
Only the signs $+$, $0$, $-$ of the components of $r$ are relevant for the corollary. We get a hyperplane arrangement in $\R^k$ by the $n$ equations
\begin{equation*}
      \begin{array}{c} s_1 b_{11}+ s_2 b_{21} + \cdots + s_k b_{k1}  = 0,\\
                                         \vdots       \\
                                        s_1 b_{1n} +s_2 b_{2n} +\cdots + s_k b_{kn} = 0,
\end{array}
\end{equation*}
where some equations can give the same hyperplane. The sign-vectors of the elements of $L_{\bf h}$ correspond bijectively to the facets of the hyperplane arrangement and may be indicated as labels on the facets.
\end{rem}

In the situation of Corollary \ref{detfacial} we have $L_{\bf h}\subseteq  \mb{ker}(A^T)$, where $A$ is the generalized Cartan matrix of the root base
 $(\h,(h_i)_{i\in I},(\al_i)_{i\in I})$ of the linear Coxeter system. If $r\in L_{\bf h}\subseteq \mb{ker}(A^T)$ then
\begin{equation}\label{pngl}
\sum_{i\in I_+(r)} \underbrace{ r_i}_{>0} a_{ij}  = \sum_{i\in I_-(r)} \underbrace{( -r_i)}_{>0} a_{ij}   \quad\mb{ for all }\quad j\in I .
\end{equation}
We now use the classification of the generalized Cartan matrices, see Theorem 3 of \cite{V}, or Theorem 4.3 of \cite{K}, to draw some conclusions for $I_+(r)$ and $I_-(r)$ from these equations.

\begin{lem}\label{kerprop1} Let $r\in L_{\bf h} \setminus\{0\}$. Then $I_+(r)$ and $I_-(r)$ are nonempty special sets. For every $j\in I\setminus(I_+(r)\cup I_-(r))$ we have
\begin{equation*}
j\mb{ is adjacent to }I_+(r) \iff j\mb{ is adjacent  to }I_-(r).
\end{equation*}
\end{lem}
\Proof We have $L_{\bf h}\cap (\R^+_0)^I=\{0\}$. Since $L_{\bf h}$ is a linear space, we find $L_{\bf h}\cap (\R^-_0)^I=\{0\}$. It follows that $I_-(r)$ and $I_+(r)$ are nonempty. 
Let $J$ be a connected component of $I_+(r)$. From (\ref{pngl}) we obtain
\begin{equation*}
\sum_{i\in J} \underbrace{(- r_i)}_{<0} a_{ij}  =  \sum_{i\in I_-(r)} \underbrace{r_i}_{<0} \underbrace{a_{ij}}_{\leq 0} \geq 0  \quad\mb{ for all }\quad j\in J.
\end{equation*}
Suppose that $J$ is of finite type. Then from Theorem 4.3 (Fin) of \cite{K} we find $(-r_i)\geq 0 $ for all $i\in J$, which is not possible. Therefore, $I_+(r)$ is special. Similarly, $I_-(r)$ is special. 

For $j\in I\setminus(I_+(r)\cup I_-(r))$ we get 
\begin{equation*}
\sum_{i\in I_+(r)} \underbrace{ r_i}_{>0} \underbrace{a_{ij}}_{\leq0}  = \sum_{i\in I_-(r)} \underbrace{( -r_i)}_{>0} \underbrace{a_{ij}}_{\leq0}
\end{equation*}
from  (\ref{pngl}), from which the remaining part of the lemma follows.
\qed

The following observation is useful to describe how for $r\in L_{\bf h}\setminus \{0\}$ the connected components of $I_+(r)$, $I_-(r)$, and $I_+(r) \cup I_-(r)$ are related: Every connected component $J$ of $I_+(r)\cup I_-(r)$ breaks up into 
\begin{equation*} 
   J_+:=J\cap I_+(r) \quad\mb{ and }\quad   J_-:=J\cap I_-(r) ,
\end{equation*}
and $J_+$,  $J_-$ are unions of connected components of $I_+(r)$, $I_-(r)$, respectively, possibly empty.
Conversely, for every connected component $K$ of $I_+(r)$, $I_-(r)$ there exists a uniquely determined connected component $J$ of $I_+(r)\cup I_-(r)$, such that $K$ is contained in $J_+$, $J_-$, respectively
\begin{lem}\label{kerprop2} Let $r\in L_{\bf h}\setminus\{0\}$. Then $I_+(r)\cup I_-(r)$ is a nonempty special set and we have:
\begin{itemize}
\item[(a)] If $J$ is an affine connected component of $I_+(r)\cup I_-(r)$, then either $J$ is a connected component of $I_+(r)\cup I_-(r)\cup I_-(r)^\bot$ contained in $I_+(r) \cap I_-(r)^\bot $, or $J$ is a connected component of $ I_+(r)\cup I_-(r)\cup I_+(r)^\bot $ contained in $I_-(r) \cap I_+(r)^\bot $.
\item[(b)]  If $J$ is an indefinite connected component of $I_+(r)\cup I_-(r)$ then $ J_+\neq\emptyset$ and $J_-\neq \emptyset$. There are only indefinite connected components of $ J_+$ and $J_-$. Every connected component of $ J_+$ is adjacent to $J_-$, every connected component of $ J_-$ is adjacent to $J_+$. 
\end{itemize}
\end{lem}
\Proof A union of special sets is again a special set. Hence we obtain from Lemma \ref{kerprop1} that $I_+(r)\cup I_-(r)$ is a nonempty special set.

To (a): Suppose that $\emptyset\neq J_+\neq J$. Then $J_+$ is a union of connected components of $I_+(r)$, which are of finite type because $J$ is affine. This is not possible, since $I_+(r)$ is special by Lemma \ref{kerprop1}. Therefore, either $J_+=J$, $J_-=\emptyset$ or $J_+=\emptyset$, $J_-=J$.

Because of $I_\pm(-r)=I_\mp(r)$ it is sufficient to treat the case $J=J_+\subseteq I_+(r)$. Since $J$ is a connected component of $I_+(r)\cup I_-(r)$, we get $J\subseteq I_-(r)^\bot$. Thus $J\subseteq I_+(r)\cap  I_-(r)^\bot$. 

Let $K$ be the connected component of $I_+(r)\cup I_-(r)\cup I_-(r)^\bot$ which contains $J$, and suppose that $J\subsetneqq K$. Then there exists $k\in K\setminus J$ such that $k$ is adjacent to $J$. Since $J$ is a connected component of $I_+(r)\cup I_-(r)$, we find $k\notin I_+(r)\cup I_-(r)$. Since $k$ is adjacent to $J\subseteq I_+(r)$, we find from Lemma \ref{kerprop1} that $k$ is also adjacent to $I_-(r)$, i.e., $k\not\in I_-(r)^\bot$.
So $k\notin I_+(r)\cup I_-(r)\cup I_-(r)^\bot$, which is not possible.

To (b): Suppose that $J_-=\emptyset$. Then $J=J_+\subseteq I_+(r)$. Since $J$ is a connected component of  $I_+(r)\cup I_-(r)$, it follows  from (\ref{pngl}) that
\begin{equation*}
  \sum_{i\in J }\underbrace{ r_i}_{>0} a_{ij}  = 0  \quad\mb{ for all }\quad j\in  J .
\end{equation*}
We find from Corollary 4.3 of \cite{K} that $J$ is of affine type. This contradicts that $J$ is of indefinite type. Thus we have shown that $J_-\neq \emptyset$. Similarly, we get $J_+\neq\emptyset$.

Let $L$ be a connected component of $J_+$. From (\ref{pngl}) we obtain
\begin{equation}\label{affindefdiff0}
\sum_{i\in L} r_i a_{ij}  = \sum_{i\in J_-} \underbrace{( -r_i)}_{>0} \underbrace{a_{ij}}_{\leq 0}\leq 0   \quad\mb{ for all }\quad j\in L .
\end{equation}
Suppose that $L$ and $J_-$ are not adjacent. Then $L$ is also a connected component of $J_+\cup J_-=J$. Since $J$ is connected, we get $J=L\subseteq I_+(r)$, which contradicts $J_-\neq\emptyset$. Thus $L$ and $J_-$ are adjacent, which shows that
\begin{equation}\label{affindefdiff}
    \sum_{i\in L} r_i  a_{ij}  = \sum_{i\in J_-} ( -r_i) a_{ij} < 0   \quad\mb{ for at least one element } j\in L .
\end{equation}
Since $L$ is also a connected component of $I_+(r)$, which is a special set by Lemma \ref{kerprop1}, it is of affine or indefinite type. Assume $L$ to be of affine type. From (\ref{affindefdiff0}) and Theorem 4.3 (Aff) of \cite{K} we obtain
\begin{equation*}
   \sum_{i\in L}  r_i a_{ij}  = 0   \quad\mb{ for all }\quad j\in L ,
\end{equation*}
which contradicts (\ref{affindefdiff}). Therefore, $L$ is of indefinite type. Similarly, every connected component of $J_-$ is of indefinite type, adjacent to $J_+$.
\qed

We have seen in Lemma \ref{kerprop1} that for $r\in L_{\bf h} \setminus\{0\}$ the sets $I_+(r)$, $I_-(r)$ are nonempty, and have only connected components of affine or indefinite type. From Lemma \ref{kerprop2} and the observation preceding this lemma we obtain immediately the following  corollary.

\begin{cor}\label{kerprop4} Let $r\in L_{\bf h} \setminus\{0\}$. We have:
\begin{itemize}
\item[(a)] Every affine connected component of $I_+(r)$ is also a connected component of $I_+(r)\cup I_-(r)$, and even of $I_+(r)\cup I_-(r)\cup I_-(r)^\bot$. No affine connected component of $I_+(r)$ is adjacent to $I_-(r)$.

\item[(b)] No indefinite connected component of $I_+(r)$ is a connected component of $I_+(r)\cup I_-(r)$. Every indefinite connected component of $I_+(r)$ is adjacent to an indefinite connected component of $I_-(r)$.

\end{itemize}
Similar statements hold if $I_+(r)$ and $I_-(r)$ are interchanged.
\end{cor}

For $r\in L_{\bf h} \setminus\{0\}$ every indefinite connected component of $I_+(r)$, $I_-(r)$ is adjacent to some indefinite connected component of $I_-(r)$, $I_+(r)$, respectively. Thus we have also shown:
\begin{cor}\label{kerprop3} Let $r\in L_{\bf h} \setminus\{0\}$. Then the following are equivalent:
\begin{itemize}
\item[(i)] $I_+(r)$ has a connected component of indefinite type.
\item[(ii)] $I_-(r)$ has a connected component of indefinite type.
\end{itemize}
The following are equivalent:
\begin{itemize}
\item[(i)] All connected components of $I_+(r)$ are of affine type.
\item[(ii)] All connected components of $I_-(r)$ are of affine type.
\end{itemize}
\end{cor}

Corollary \ref{detfacial}, Lemmas \ref{kerprop1}, \ref{kerprop2}, and its Corollaries \ref{kerprop4}, \ref{kerprop3} allow to derive quickly some properties of the facial sets. The following theorem, except the implication ``$\Rightarrow$", has been given as Theorem 4 1) of \cite{V}. For symmetric linear Coxeter systems it corresponds to Lemma 2.10 (b) and Lemma 8.5 (a) of \cite{D}.
\begin{thm}\label{fsp} For $J\subseteq I $ we have:
\begin{equation}\label{main fa sp fa}
  J\mb{ is facial } \iff J^\infty \mb{ is facial}.
\end{equation}
In particular, $J$ is facial if $J^0=J$, i.e., if $J=\emptyset$ or all connected components of $J$ are of finite type.
\end{thm}
\Proof For every $r\in L_{\bf h}\setminus\{0\}$ the sets $I_+(r)$ and $I_-(r)$ are special by Lemma \ref{kerprop1}. Therefore,
\begin{equation*}
   I_+(r)\subseteq J \iff  I_+(r)\subseteq J^\infty \quad \mb{ and }\quad I_-(r)\subseteq J \iff  I_-(r)\subseteq J^\infty .
\end{equation*}
Now (\ref{main fa sp fa}) follows from Corollary \ref{detfacial} and Remark \ref{detfacial}. If $J^0=J$ then $J^\infty=\emptyset$, which is a facial set by Proposition \ref{fsFa0} (c). Thus $J$ is facial.
\qed

We say that an indecomposable generalized Cartan matrix $B=(b_{ij})_{i,j\in I}$ is of 0-hyperbolic type, or briefly of hyperbolic type, if $B$ is of indefinite type such that for all $i\in I$ the indecomposable components of the generalized Cartan submatrix $B_{I\setminus\{i\}}$ are of finite or affine type. 
We say that $B$ is of 1-hyperbolic type, if $B$ is of indefinite type such that for all $i_1, i_2\in I$, $i_1\neq i_2$ the indecomposable components of the generalized Cartan submatrix $B_{I\setminus\{i_1,i_2\}}$ are of finite or affine type, and if $B$ is not of $0$-hyperbolic type, 
\begin{prop}\label{prop hyp} Let $B=(b_{ij})_{i,j\in I}$ be an indecomposable generalized Cartan matrix of 0-hyperbolic or 1-hyperbolic type. Let $\emptyset\neq J\subseteq I$ be connected, of nonfinite type. Then we have $|J|\geq |I|-2$. If $|J|= |I|-2$ then $J$ is of affine type.
\end{prop}
\Proof Suppose there exists $\emptyset\neq J\subseteq I$ connected, of nonfinite type with $|J|\leq |I|-3$. Because $I$ is connected there exists 
a sequence 
\begin{equation*}
 J\subsetneqq J\cup \{i_1\} \subsetneqq J\cup \{i_1, i_2\} \subsetneqq \cdots \subsetneqq J\cup \{i_1,i_2,\ldots,i_m\} = I
\end{equation*}
of connected subsets of $I$ with $m\geq 3$. Since $J$ is of nonfinite type, the sets $J\cup\{i_1,\ldots, i_{m-2}\}$ and  $J\cup\{i_1,\ldots, i_{m-1}\}$ are of indefinite type. This contradicts the definition of 1-hyperbolic and 0-hyperbolic.

If $|J|=|I|-2$, then from the definition of 0-hyperbolic and 1-hyperbolic it follows that $J$ is a union of connected components of finite and affine type. We conclude that $J$ is of affine type.
\qed

As a further example how to work with Corollary \ref{detfacial}, Lemmas \ref{kerprop1} and \ref{kerprop2} we prove the following theorem. Part (a) has been given for symmetric linear Coxeter systems as Lemma 8.5 (c) of \cite{D}.
\begin{thm}\label{aff01hyp} Let $J\subseteq I$ be facial. Then:
\begin{itemize}
\item[(a)] Every union of connected components of $J$, which contains all the affine connected components of $J$, is facial.
\item[(b)] Every union of connected components of $J$, which are 0-hyperbolic or 1-hyperbolic, is facial.
\end{itemize}
\end{thm}
\Proof To (a): Let $F$ be a union of components of $J$, which contains all the affine components of $J$. Let $R$ be the union of the remaining components of $J$. 

Let $r\in L_{\bf h}\setminus\{0\}$ such that $I_+(r)\subseteq F$. Since $F\cup R$ is facial, we find $I_-(r) \subseteq F\cup R$ by Corollary \ref{detfacial}. Let $K$ be a component of $I_-(r)$. It is of nonfinite type by Lemma \ref{kerprop1}. Since $F$, $R$ are separated, $K$ has to be contained either in $F$ or in $R$. 

Suppose that $K\subseteq R$.
If $K$ is of affine type, then it is a component of $I_+(r)\cup I_-(r)\cup I_+(r)^\bot$ by Corollary \ref{kerprop4}. Since $K\subseteq R\subseteq F^\bot\subseteq  I_+(r)^\bot\subseteq I_+(r)\cup I_-(r)\cup I_+(r)^\bot$, it is also a component of $R$, which is not possible by the definition of $R$. 
If $K$ is of indefinite type, then by Corollary \ref{kerprop4} it is adjacent to $I_+(r)\subseteq F$, which contradicts that $F$ and $R$ are separated. 

The case $K\subseteq F$ remains. Thus we have shown that $I_-(r)\subseteq F$. From Corollary \ref{detfacial} and Remark \ref{remdetfacial} we find that $F$ is facial.

To (b): Let $U$ be a union of 0-hyperbolic or 1-hyperbolic components of $J$. From Corollary  \ref{detfacial} and Remark \ref{remdetfacial} it follows that $U$ is certainly facial if there are no elements $r\in L_{\bf h}\setminus\{0\}$ such that $I_+(r)\subseteq U$. We prove this indirectly.
Suppose there exists $r\in L_{\bf h}\setminus\{0\}$ such that $I_+(r)\subseteq U$. Since $J$ is facial and $I_+(r)\subseteq U\subseteq J$, we get $I_-(r)\subseteq J$.

Let $K$ be a component of $I_+(r)$. Then $K$ is of nonfinite type by Lemma \ref{kerprop1}, and there exists a 0-hyperbolic or 1-hyperbolic component $H$ of $J$ such that $K\subseteq H\subseteq U$. From Proposition \ref{prop hyp} we obtain $|K|\geq |H|-2$.

If $K=H$ then $K$ is of indefinite type. From Corollary \ref{kerprop4} we find that $K$ is adjacent to $I_-(r)\subseteq J$, which is not possible because $K=H$ is a component of $J$.

Let $K=H\setminus\{i_0\}$ for some $i_0\in H$. Since $H$ is connected, $i_0$ is adjacent to $K$. Since $K$ is a component of $I_+(r)$, also $i_0\notin I_+(r)$. Hence $I_-(r)\cap H \subseteq \{i_0\}$. Note that $I_-(r)\cap H$ is either empty or a union of components of $I_-(r)$, because we have $I_-(r)\subseteq J$ and $H$ is a component of $J$. Since $I_-(r)$ contains no components of finite type by Lemma \ref{kerprop1}, we find $i_0\notin I_-(r)$. From Lemma \ref{kerprop1} it also follows that $i_0$ is adjacent to $I_-(r)\subseteq J$, which is not possible because $H$ is a component of $J$.

Let $K=H\setminus\{i_0,\,i_1\}$ for some $i_0,i_1\in H$ with $i_0\neq i_1$. Then $K$ is of affine type by Proposition \ref{prop hyp}, and $i_0$ or $i_1$ is adjacent to $K$. We may assume that $i_0$ is adjacent to $K$. Since $K$ is a component of $I_+(r)$, also $i_0\notin I_+(r)$. From Corollary \ref{kerprop4} we find $K\subseteq I_+(r)\cap I_-(r)^\bot\subseteq I_-(r)^\bot$. Hence $i_0\notin I_-(r)$, and $I_-(r)\cap H \subseteq \{i_1\}$. Since $I_-(r)$ contains no components of finite type by Lemma \ref{kerprop1}, we get $i_1\notin I_-(r)$. Thus $H\cap I_-(r)=\emptyset$. From Lemma \ref{kerprop1} it also follows that $i_0$ is adjacent to $I_-(r)\subseteq J$, which is not possible because $H$ is a component of $J$.
\qed

Let $J\subseteq I$. From $J\subseteq J^{fc}$ we conclude that $J^\bot \supseteq  (J^{fc})^\bot $.
\begin{prop}\label{jbjfcb} For $J\subseteq I $ we have $J\cup J^\bot\subseteq J^{fc}\cup (J^{fc})^\bot$.
\end{prop}
\Proof Suppose there exists $i_0\in J^\bot\setminus  (J^{fc}\cup (J^{fc})^\bot) $. Because of $i_0\in J^\bot$ and $i_0\notin (J^{fc})^\bot$  we get
$J\subseteq J^{fc}\cap\{i_0\}^\bot\subsetneqq J^{fc}$. We now use Corollary \ref{detfacial} and Remark \ref{remdetfacial} to show that $J^{fc}\cap\{i_0\}^\bot$ is facial, which contradicts the fact that $J^{fc}$ is the smallest facial subset of $I$ which contains $J$. 

Let $r\in L_{\bf h}\setminus\{0\}$ such that $I_+(r)\subseteq J^{fc}\cap\{i_0\}^\bot$, i.e., $I_+(r)\subseteq J^{fc}$ and $i_0\in I\setminus I_+(r)$ not adjacent to $I_+(r)$. Since $J^{fc}$ is facial, we get $I_-(r)\subseteq J^{fc}$. Thus $i_0\notin I_+(r)\cup I_-(r)$. By Lemma \ref{kerprop1} we find  $I_-(r)\subseteq J^{fc}\cap\{i_0\}^\bot$.
\qed

The next theorem can be considered as a generalization of Theorem 4 2) of \cite{V}, which follows from a special case of part (a) by applying Theorem \ref{aff01hyp} (a). 
\begin{thm}\label{fac centralizer} Let $J, J_1 \subseteq I$ be special. 
\begin{itemize}
\item[(a)] Then $J^{aff}\cup J^\bot$ is a facial set which contains $I^0\cup I^{aff}$.
\item[(b)] If $J\subseteq J_1$ then $J^{aff}\cup J^\bot \supseteq J_1^{aff}\cup J_1^\bot$.
\end{itemize}
\end{thm}
\Proof If $J=\emptyset$ then $J^{aff}\cup J^\bot=I$. Thus (a) and (b) hold trivially. Now let $J\ne\emptyset$. To (b): From $J\subseteq J_1$ we get $J^\bot\supseteq J_1^\bot$. Let $K$ be an affine component of $J_1$. We show $J\cup J^\bot\supseteq K$ indirectly. Suppose that there exists $i\in K\setminus(J\cup J^\bot)$.
Since $i\in K$ is adjacent to $J\subseteq J_1$, and $K$ is a component of $J_1$, we conclude that $J\cap K\neq\emptyset$. Moreover, $J\cap K$ is a component of the special set $J\subseteq J_1$, because $K$ is a component of $J_1$. Hence $J\cap K$ is a proper subset of $K$ of affine or indefinite type, which contradicts the fact that $K$ is affine.

To (a): By Theorem 4.3 (Aff), (Ind) of \cite{K} there exists $r\in (\R^+)^J$ such that 
\begin{equation*}
    \sum_{j\in J} a_{ij} r_j \left\{\begin{array}{lcl} 
       = 0  & \mb{if} & i\in J^{aff},\\
        <0  & \mb{if} & i\in J^{ind}.
     \end{array}\right.
\end{equation*}
Clearly, if $i\in J^\bot$ then $\sum_{j\in J} a_{ij} r_j =0$. If $i\in I\setminus (J\cup J^\bot)$ then $a_{ij}<0$ for some $j\in J$, and $a_{ij}\leq 0$ for all $j\in J$. Hence $\sum_{j\in J} a_{ij} r_j <0$. 
Thus $ -\sum_{j\in J} r_j \al_j \in F_{J^{aff}\cup J^\bot}$. Proposition \ref{fshfree} (a) shows that $J^{aff}\cup J^\bot$ is facial. From $J\subseteq I^\infty$ we get $J^{aff}\cup J^\bot \supseteq (I^\infty)^{aff}\cup (I^\infty)^\bot =I^{aff}\cup I^0$ by part (b).
\qed

For $J_1, J_2\subseteq I$ and $\sigma\in\We$ we denote by $J_1\cap\sigma J_2$ the subset of $I$ corresponding to the simple reduced coroots
\begin{equation*}
     \Mklz{\R^+_0 h_i}{i\in J_1}\cap \sigma\Mklz{\R^+_0 h_j}{j\in J_2},
\end{equation*}
or equivalently, see Corollary \ref{redrefl}, to the simple reflections
\begin{equation*}
     \Mklz{\sigma_i}{i\in J_1}\cap \Mklz{\sigma \sigma_j \sigma^{-1}}{j\in J_2}.
\end{equation*}
We will need the following theorem, which corresponds for symmetric linear Coxeter systems to Lemma 2.10 (d) of \cite{D}. We prove it by a direct computation.
\begin{thm}\label{JsJ} Let $J_1, J_2\subseteq I$ be facial and $\sigma\in\mb{}^{J_1}\We^{J_2}$. Then
\begin{equation}\label{kmroots}
(\sum_{i\in J_1}\R _0^+ h_i ) \cap\sigma \sum_{i\in J_2}\R^+_0 h_i =\sum_{i\in J_1\cap\sigma J_2}\R^+_0 h_i  
\end{equation}
is a face of $\sum_{i\in I}\R^+_0 h_i$. In particular, $J_1\cap\sigma J_2$ is facial.
\end{thm}
\begin{rem} For $J_1, J_2\subseteq I$ and $\sigma\in\mb{}^{J_1}\We^{J_2}$ we have the formula 
\begin{equation*}
    \We_{J_1} \cap \sigma\We_{J_2}\sigma^{-1} =\We_{J_1\cap \sigma J_2} 
\end{equation*} 
of R. Kilmoyer, see for example Proposition 2.25 of \cite{AB}. Specialized to $J_1, J_2$ facial, it can be regarded as a supplement to the theorem.
\end{rem}
\Proof (a)  Let $J\subseteq I$ be facial, $K\subseteq I$, and $\tau\in\mb{}^J\We^K$. Let
\begin{equation} \label{kmzwr}
  \tau \sum_{k\in K} s_k h_k = \sum_{k\in K} s_k \tau h_k \in\sum_{j\in J}\R^+_0 h_j
\end{equation}
with coefficients $s_k\in\R^+_0$, $k\in K$. Let $k_0\in K$ such that $s_{k_0}>0$. We show that there exists $j\in J$ such that $\R^+_0 \tau h_{k_0}= \R^+_0 h_j$. We make repeated use of Theorem \ref{ThHN} without mentioning it further.

Because of $l(\tau\sigma_k)>l(\tau)$ we have $  \tau h_k\in\sum_{i\in I}\R^+_0 h_i$ for all $k\in K$. Since $J$ is facial and $s_{k_0}\neq 0$, we find $ \tau h_{k_0}\in\sum_{j\in J}\R^+_0 h_j$ from (\ref{kmzwr}). Therefore, we have
\begin{equation*} 
  \R^+_0 h_{k_0}\subseteq \sum_{j\in J}\R^+_0\tau^{-1} h_j  \quad\mb{with }\quad  \tau^{-1} h_j\in \sum_{i\in I}\R^+_0 h_i
\end{equation*}
because of $l(\tau^{-1}\sigma_j)>l(\tau)$ for all $j\in J$. Since $\R^+_0 h_{k_0}$ is a one-dimensional face of $\sum_{i\in I}\R^+_0 h_i$, there exists $j\in J$ such that $\R^+_0 h_{k_0}=\R^+_0\tau^{-1}h_j$.

(b) We now show formula (\ref{kmroots}). The inclusion ``$\supseteq$" is trivial. To show the reverse inclusion let
\begin{equation*}
  \sum_{i\in J_1}r_i h_i =\sigma \sum_{j\in J_2}s_j h_j 
\end{equation*}
with coefficients $r_i\in\R^+_0$, $i\in J_1$, and $s_j\in\R^+_0$, $j\in J_2$. We conclude from (a) that for every $j\in J_2$ with $s_j\neq 0$, there exists $i\in J_1$, such that $\R^+_0\sigma h_j=\R^+_0h_i$. This proves the inclusion ``$\subseteq$".

(c) It remains to show that the convex cone (\ref{kmroots}) is a face of $\sum_{i\in I}\R^+_0 h_i$. Since $J_1$ is facial, it is sufficient to show that it is a face of $\sum_{i\in J_1}\R^+_0 h_i$. 

Let $h_1=\sum_{i\in J_1} r_i h_i $ and $h_2=\sum_{i\in J_1}s_i h_i$ with coefficients $r_i, s_i \in\R^+_0$, $i\in J_1$, and let $s_{i_0}>0$ for at least one $i_0\in J_1\setminus ( J_1\cap\sigma J_2)$. Suppose that
\begin{equation*}
  h_1+h_2=\sum_{i\in J_1}(r_i+s_i) h_i \in \sum_{i\in J_1\cap\sigma J_2}\R^+_0 h_i  .
\end{equation*}
Then we get
\begin{equation*}
  \sigma^{-1}(h_1+h_2) = \sum_{i\in J_1}(r_i+s_i) \sigma^{-1}h_i \in \sum_{j\in  J_2}\R^+_0 h_j  .
\end{equation*}
Here $r_{i_0}+s_{i_0}>0$. By (a) there exists $j\in J_2$ such that $\R^+_0 \sigma^{-1}h_{i_0} =\R^+_0 h_j$. This contradicts $i_0\notin J_1\cap\sigma J_2$. 
\qed

\subsection{\label{SubsectionEx}An example for the computation of the facial sets}
We consider the generalized Cartan matrix
\begin{equation*}
  A=\left(\begin{array}{rrrrrr}
     2  & -2 &  0 &  0 &  0 & -2\\
    -2  &  2 & -2 &  0 &  0 &  0\\
     0  & -2 &  2 & -2 &  0 &  0\\
     0  &  0 & -2 &  2 & -2 &  0\\
     0  &  0 &  0 & -2 &  2 & -2\\
    -2  &  0 &  0 &  0 & -2 &  2 
\end{array}\right).
\end{equation*}
The Coxeter diagram of the associated Coxeter matrix, as defined in Section 3.2 of \cite{AB}, is given by: 
%
%
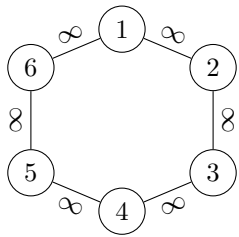
\begin{equation*}
\begin{tikzpicture}[scale=0.35]
\path node (1) at (0, 3.46)[shape=circle, draw]  {$1$}  
         node (2) at (3.46,2) [shape=circle, draw] {$2$}
         node (3) at (3.46,-2) [shape=circle, draw]  {$3$}   
         node (4) at (0,-3.46)[shape=circle, draw]  {$4$} 
         node (5) at (-3.46,-2) [shape=circle, draw]  {$5$} 
         node (6) at (-3.46, 2)[shape=circle, draw]  {$6$}    ;
\draw [-] (1) -- (2) node[midway, sloped, above] {$\infty$}  -- (3) node[midway, sloped, below] {$\infty$} -- (4) node[midway, sloped, below] {$\infty$} --(5) node[midway, sloped, below] {$\infty$}  -- (6) node[midway, sloped, above] {$\infty$} -- (1) node[midway, sloped, above] {$\infty$};
\end{tikzpicture}
\end{equation*}
%
%
%
We find from $A$ that a nonempty connected subset $J\subseteq \{1,\,2,\,\ldots,\,6\}$ is of finite type if $|J|=1$, of affine type if $|J|=2$, and of indefinite type if $|J|\geq 3$.

Since $A$ is of indefinite type, we have  $\mb{ker}(A^T)\cap (\R^+_0)^6 =\{0\}$. We obtain by
\begin{equation*}
\frac{1}{2\sqrt{3}}\left(\begin{array}{c}
  1 \\2\\1\\-1\\-2\\-1
\end{array}\right), \;\;
\frac{1}{2}\left(\begin{array}{c}
  1 \\0\\-1\\-1\\0\\1
\end{array}\right)
\end{equation*}
a linear base of $\mb{ker}(A^T)$, for which the picture of the corresponding labeled hyperplane arrangement, described in Remark \ref{hypfa}, is symmetric:\vspace*{1ex}\\
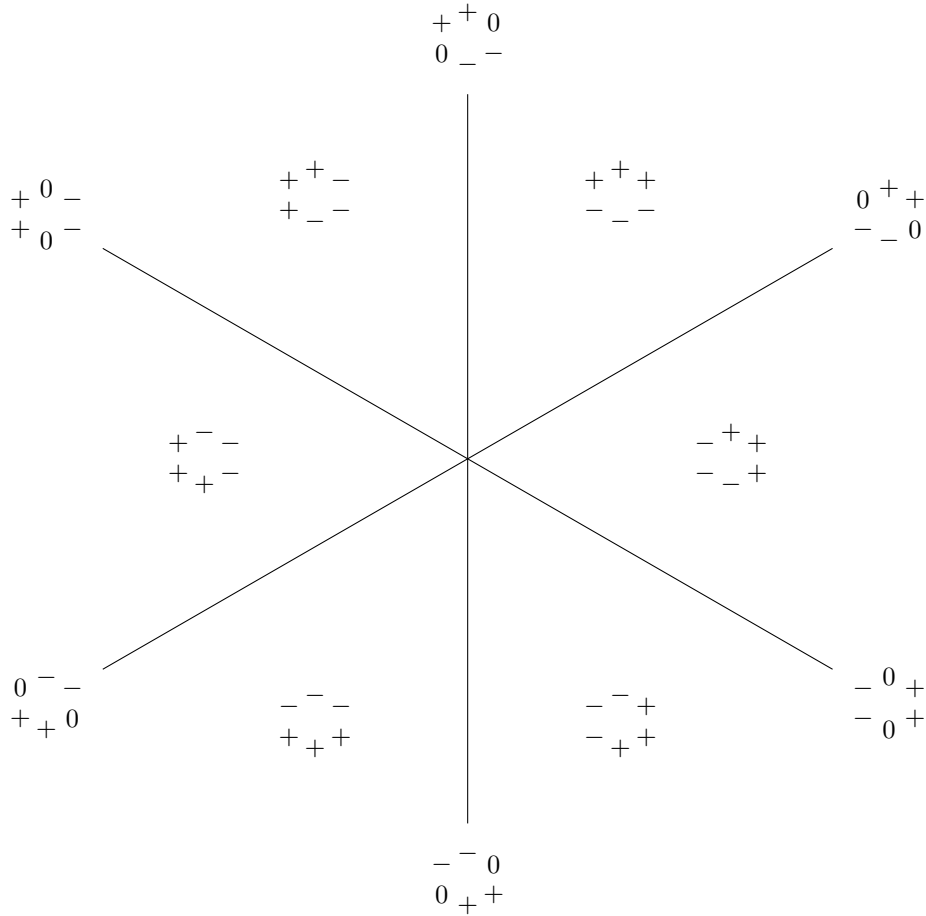
\hspace*{4em}\begin{tikzpicture}[scale=0.7]
\draw (0,-6.93) -- (0, 6.93);
\draw (-6.93, 4) -- (6.93,-4);
\draw (-6.93,-4) -- (6.93,4);
%
%
%
%
\draw (5,0) node{\mbox{\begin{tikzpicture}[scale=0.1]
\path node (1) at (0, 3.464) {$+$}
         node  (2) at (3.46,2) {$+$}
         node (3) at (3.46,-2)  {$+$}   
         node (4) at (0,-3.46) {$ -$} 
         node (5) at (-3.46,-2)  {$-$} 
         node (6) at (-3.46, 2) {$-$}    ;
\end{tikzpicture}}};
\draw (-5,0) node{\mbox{\begin{tikzpicture}[scale=0.1]
\path node (1) at (0, 3.464) {$-$}
         node  (2) at (3.46,2) {$-$}
         node (3) at (3.46,-2)  {$-$}   
         node (4) at (0,-3.46) {$ +$} 
         node (5) at (-3.46,-2)  {$+$} 
         node (6) at (-3.46, 2) {$+$}    ;
\end{tikzpicture}}};
\draw (2.9,5) node{\mbox{\begin{tikzpicture}[scale=0.1]
\path node (1) at (0, 3.464) {$+$}
         node  (2) at (3.46,2) {$+$}
         node (3) at (3.46,-2)  {$-$}   
         node (4) at (0,-3.46) {$ -$} 
         node (5) at (-3.46,-2)  {$-$} 
         node (6) at (-3.46, 2) {$+$}    ;
\end{tikzpicture}}};
\draw (-2.9,-5) node{\mbox{\begin{tikzpicture}[scale=0.1]
\path node (1) at (0, 3.464) {$-$}
         node  (2) at (3.46,2) {$-$}
         node (3) at (3.46,-2)  {$+$}   
         node (4) at (0,-3.46) {$ +$} 
         node (5) at (-3.46,-2)  {$+$} 
         node (6) at (-3.46, 2) {$-$}    ;
\end{tikzpicture}}};
\draw (-2.9, 5) node{\mbox{\begin{tikzpicture}[scale=0.1]
\path node (1) at (0, 3.464) {$+$}
         node  (2) at (3.46,2) {$-$}
         node (3) at (3.46,-2)  {$-$}   
         node (4) at (0,-3.46) {$ -$} 
         node (5) at (-3.46,-2)  {$+$} 
         node (6) at (-3.46, 2) {$+$}    ;
\end{tikzpicture}}};
\draw (2.9,-5) node{\mbox{\begin{tikzpicture}[scale=0.1]
\path node (1) at (0, 3.464) {$-$}
         node  (2) at (3.46,2) {$+$}
         node (3) at (3.46,-2)  {$+$}   
         node (4) at (0,-3.46) {$ +$} 
         node (5) at (-3.46,-2)  {$-$} 
         node (6) at (-3.46, 2) {$-$}    ;
\end{tikzpicture}}};
%
%
%
%
\draw (0,8) node{\mbox{\begin{tikzpicture}[scale=0.1]
\path node (1) at (0, 3.464) {$+$}
         node  (2) at (3.46,2) {$0$}
         node (3) at (3.46,-2)  {$-$}   
         node (4) at (0,-3.46) {$ -$} 
         node (5) at (-3.46,-2)  {$0$} 
         node (6) at (-3.46, 2) {$+$}    ;
\end{tikzpicture}}};
\draw (0,-8) node{\mbox{\begin{tikzpicture}[scale=0.1]
\path node (1) at (0, 3.464) {$-$}
         node  (2) at (3.46,2) {$0$}
         node (3) at (3.46,-2)  {$+$}   
         node (4) at (0,-3.46) {$ +$} 
         node (5) at (-3.46,-2)  {$0$} 
         node (6) at (-3.46, 2) {$-$}    ;
\end{tikzpicture}}};
\draw (8,4.64) node{\mbox{\begin{tikzpicture}[scale=0.1]
\path node (1) at (0, 3.464) {$+$}
         node  (2) at (3.46,2) {$+$}
         node (3) at (3.46,-2)  {$0$}   
         node (4) at (0,-3.46) {$ -$} 
         node (5) at (-3.46,-2)  {$-$} 
         node (6) at (-3.46, 2) {$0$}    ;
\end{tikzpicture}}};
\draw (-8,-4.64) node{\mbox{\begin{tikzpicture}[scale=0.1]
\path node (1) at (0, 3.464) {$-$}
         node  (2) at (3.46,2) {$-$}
         node (3) at (3.46,-2)  {$0$}   
         node (4) at (0,-3.46) {$ +$} 
         node (5) at (-3.46,-2)  {$+$} 
         node (6) at (-3.46, 2) {$0$}    ;
\end{tikzpicture}}};
\draw (8., -4.64) node{\mbox{\begin{tikzpicture}[scale=0.1]
\path node (1) at (0, 3.464) {$0$}
         node  (2) at (3.46,2) {$+$}
         node (3) at (3.46,-2)  {$+$}   
         node (4) at (0,-3.46) {$ 0$} 
         node (5) at (-3.46,-2)  {$-$} 
         node (6) at (-3.46, 2) {$-$}    ;
\end{tikzpicture}}};
\draw (-8,4.64) node{\mbox{\begin{tikzpicture}[scale=0.1]
\path node (1) at (0, 3.464) {$0$}
         node  (2) at (3.46,2) {$-$}
         node (3) at (3.46,-2)  {$-$}   
         node (4) at (0,-3.46) {$ 0$} 
         node (5) at (-3.46,-2)  {$+$} 
         node (6) at (-3.46, 2) {$+$}    ;
\end{tikzpicture}}};
\end{tikzpicture}\vspace*{1ex}\\
In this picture we indicate the sign vectors of the elements of $\mb{ker}(A^T)$ as labels on the vertices of the Coxeter diagram. The sign vector of the zero of  $\mb{ker}(A^T)$ is omitted. 

The cyclic group $C_6$ acts on the labeled hyperplane arrangement. Different subspaces $L_{\bf h}$ of $\mb{ker}(A^T)$ give, up to cyclic symmetry, only four possible systems of facial sets. It is trivial to compute the facial sets corresponding to $L_{\bf h}$ with Corollary \ref{detfacial}. We only list the special facial sets, from which the facial sets can be obtained by Theorem \ref{fsp}:

(a) If $L_{\bf h}=\{0\}$ then every subset of $\{1,2,3,4,5,6\}$ is facial. The special facial sets coincide with the special sets, which are
\begin{eqnarray*}
&\emptyset,&\\
&\{1, 2\},\,\{2, 3\},\,\{3, 4\},\,\{4, 5\},\,\{5, 6\},\,\{1, 6\},& \\
&\{1, 2, 3\},\, \{2, 3, 4\},\,\{3, 4, 5\},\, \{4, 5, 6\},\, \{1, 5, 6\},\,\{1, 2, 6\}, &\\
&\{ 1, 2, 3, 4 \},\,\{2, 3, 4, 5\},\,\{ 3, 4, 5, 6 \},\, \{ 1, 4, 5, 6 \},\, \{1, 2, 5, 6\},\, \{ 1, 2, 3, 6 \}, &\\
&\{1, 2\}\cup\{4, 5\},\,   \{2,\,3\}\cup\{5, 6\},\,   \{1,\, 6\}\cup\{ 3, 4\}, &\\
&\{1, 2, 3, 4, 5\},\,  \{2, 3, 4, 5, 6\},\,  \{1, 3, 4, 5, 6\},\,   \{1, 2, 4, 5, 6\},\,  \{1, 2, 3, 5, 6\},\,  \{1, 2, 3, 4, 6\}, &\\
&\{1, 2, 3, 4, 5, 6\}.&
\end{eqnarray*}
There are 29 special facial sets and 64 facial sets.

(b) $L_{\bf h}$ is spanned by an element of $\mb{ker}(A^T)$, which corresponds to a point in an open chamber of the hyper\-plane arrangement of $\mb{ker}(A^T)$. Take for example $L_{\bf h}=\R (1,2,1,-1,-2,-1)^T$. Here the special facial sets are
\begin{eqnarray*}
 & \emptyset, &\\
& \{1, 2\},\,\{2, 3\},\,\{3, 4\},\,\{4, 5\},\,\{5, 6\},\,\{1, 6\},& \\
&\{2, 3, 4\},\,\{3, 4, 5\},\,\{1, 5, 6\},\,\{1, 2, 6\}, &\\
&\{2, 3, 4, 5\},\,   \{1, 2, 5, 6\},&\\
&\{1, 2\}\cup\{4, 5\},\,  \{2, 3\}\cup\{5, 6\},\,\{1, 6\}\cup\{ 3, 4 \},&\\
&\{1,\,2,\,3,\,4,\,5,\,6\}.&
\end{eqnarray*}
There are 17 special facial sets and 50 facial sets.

(c) $L_{\bf h}$ is spanned by an element of $\mb{ker}(A^T)$, which corresponds to a point in an open panel of the hyper\-plane arrangement of $\mb{ker}(A^T)$. Take for example $L_{\bf h}=\R (1,0,-1,-1,0,1)^T$. Here the special facial sets are
\begin{eqnarray*}
&\emptyset,&\\
&\{1, 2\},\,\{2, 3\},\, \{4, 5\},\,\{5, 6\},& \\
&\{1, 2, 3\},\, \{4, 5, 6\}, &\\
&\{1, 2\}\cup\{4, 5\},\,   \{2,\,3\}\cup\{5, 6\},\,   \{1,\, 6\}\cup\{ 3, 4\},&\\
& \{1, 3, 4, 5, 6\},\, \{1, 2, 3, 4, 6\},&\\
&\{1, 2, 3, 4, 5, 6\} .&
\end{eqnarray*}
There are 13 special facial sets and 40 facial sets.

(d) If $L_{\bf h}=\mb{ker}(A^T)$ then the special facial sets are
\begin{eqnarray*}
   &\emptyset ,&\\
   &\{1, 2\}\cup\{4, 5\},\,   \{2,\,3\}\cup\{5, 6\},\,   \{1,\, 6\}\cup\{ 3, 4\},&\\
    &\{1,2,3,4,5,6\}.&
\end{eqnarray*}
There are 5 special facial sets and 22 facial sets.
%
%
%
%
%
\subsection{\label{SubsectionFT}The faces of the Tits cone}
%
%
%
%
%
%
%
%
In this subsection we determine the faces of the Tits cone of a linear Coxeter system and describe the lattice operations of its face lattice. These results generalize results which have been obtained for free integral linear Coxeter systems by E. Looijenga, P. Slodowy, and the author, and for symmetric linear Coxeter systems by M. Dyer. Our approach is motivated by the theory of Coxeter group invariant convex cones in Section \ref{SecInvsub}.

Let $(\h,(h_i)_{i\in I},(\al_i)_{i\in I})$ be a root base of the linear Coxeter system $(\h,\,(H_i, \,L_i)_{i\in I})$. For $J\subseteq I$ and $\sigma\in \We$ we call $\sigma F_J$ an open facet, and $\sigma \overline{F}_J$ a closed facet of the Tits cone of type $J$. Theorem \ref{fuprop} and Proposition \ref{rihfc} imply that the partition of the fundamental chamber
\begin{equation*} 
   \overline{C} = \dot{ \bigcup_{J\text{ facial}} } F_J
\end{equation*} 
induces a $\We$-invariant partition of the Tits cone $X$: 
\begin{prop} (a) Let $J$, $L$ be facial sets and $\sigma, \tau\in\We$. Then we have $\sigma F_J \cap \tau F_L\neq\emptyset$ if and only if  $\sigma F_J =\tau F_L$, if and only if $J=L$ and $\sigma^{-1}\tau\in\We_{L}$.

(b) $X$ is the disjoint union of the open facets, whose types are facial sets.
\end{prop}

Recall that for $\Th\subseteq I $ we have $\mb{ri}( \sum_{i\in\Th}\R^+_0 h_i )=\sum_{i\in\Th}\R^+ h_i $. The classification of the generalized Cartan matrices in Theorem 3 of \cite{V}, or Theorem 4.3 of \cite{K} implies:
\begin{prop}\label{Ksp} Let $\Th\subseteq I $. Then $\Th$ is special if and only if $(\sum_{i\in\Th}\R^+ h_i ) \cap(-\overline{C}^\vee)\neq\emptyset$.
\end{prop}

The following Lemma \ref{LooR}, and in essential also the Corollaries \ref{specialfaceX} and \ref{stabRTh} are due to E. Looijenga for free integral linear Coxeter systems, see Lemma 2.2, and Corollary 2.3 of \cite{Loo}. Actually, E. Looijenga did not investigate the faces of the Tits cone but was interested in the linear spaces
\begin{equation*}
   \Mklz{\la\in \h^*}{\la(h_i)=0\,\mb{ for all }\,i\in\Th} ,\quad\Th\mb{ special},
\end{equation*}
and their stabilizers for the construction of partial compactifications of certain orbit spaces. These results have been used by P. Slodowy to investigate the faces of the Tits cone in Chapter 6.2 of \cite{Sl}. Corollaries \ref{specialfaceX} and \ref{stabRTh} can be obtained for symmetric linear Coxeter systems as a combination of parts of Theorem 10.3 and Lemma 10.4 of \cite{D}.
\begin{lem}\label{LooR} Let $\Th\subseteq I$ be special facial and $h\in\left(\sum_{i\in\Th}\R^+\,h_i\right)\cap (-\Cq^\vee)$. Let $\la\in X$ such that $\la(h)\leq 0$. Then we have $\la(h)=0$ and $\la\in\We_{\Th^\bot}\FTq $. Moreover, $\la(h_i)=0$ for all $i\in\Th$.
\end{lem}
\Proof The proof is similar to the inductive proof of Lemma 2.2 in \cite{Loo}, if the characterization of the Tits cone $X$ given as formula (\ref{chTcrf1}) in Theorem \ref{fuprop} is used. \qed

\begin{cor}\label{specialfaceX} Let $\Th\subseteq I$ be special facial. Then 
\begin{equation*}
  R(\Th):=\We_{\Th^\bot}\FTq
\end{equation*}
is an exposed face of $X$, i.e., for every element $h\in\left(\sum_{i\in\Th}\R^+\,h_i\right)\cap (-\Cq^\vee)$ we have
\begin{equation*}
   X     \subseteq  \Mklz{\la\in \h^*}{\la(h)\geq 0} \quad { and } \quad 
  R(\Th)  =        X\cap \Mklz{\la\in \h^*}{\la(h)=0}.
\end{equation*}
Its linear hull is 
\begin{equation*} 
        R(\Th) -  R(\Th) = \Mklz{\la\in \h^*}{\la(h_i)=0\,\mb{ for all }\,i\in\Th} .
\end{equation*}
\end{cor}
\begin{rem} In particular, $R(\emptyset)= X$. We will see in Corollary \ref{descsf} that $R(I^\infty)$ coincides with the smallest face of $X$.
\end{rem}

\begin{cor}\label{stabRTh}  Let $\Th\subseteq I$ be special facial. Then the pointwise stabilizer of $R(\Th)$ and the stabilizer of $R(\Th)$ as a whole are given by 
\begin{equation*}
    Z_{\mathcal W}(R(\Th))  = \We_{\Th}  \quad \mb{ and }\quad  N_{\mathcal W}(R(\Th)) =\We_{\Th\cup \Th^\bot} .
\end{equation*}
\end{cor}

The following proposition has been given for minimal free integral linear Coxeter systems as Proposition 2.5 of \cite{M2}. It can be proved similarly as in \cite{M2}. For symmetric linear Coxeter sytems it can be obtained as a combination of parts of Theorem 10.3 and Lemma 10.4 of \cite{D}.
\begin{prop}\label{Rspecialinc} Let $\sigma,\sigma'\in\We$, and $\Th,\Th'\subseteq I$ be special facial. Then:
\begin{equation*}
   \sigma'R(\Th') \subseteq \sigma R(\Th)\quad \iff \quad \Th' \supseteq \Th \mb{ and } \sigma^{-1}\sigma' \in \We_{\Th^\bot}\We_{\Th'}.  
\end{equation*}
\end{prop}

D. Krammer introduced in Section 2.2 of \cite{Kr} certain projectors, which may be adopted to our situation as follows: For $J_f\subseteq I$ such that $(J_f)^0 = J_f$ define a linear map $\mb{mid}_{J_f}: \h ^*\to \h ^*$ by
\begin{equation*}
   \mb{mid}_{J_f}(\la):=   \frac{1}{|{\cal W}_{J_f}|}\sum_{\sigma\in{\cal W}_{J_f}}\sigma\la,  \quad \la \in \h^* . 
\end{equation*}
In particular, $ \mb{mid}_{\,\emptyset}$ is the identity on $\h^*$. For $i\in I$ the map $ \mb{mid}_{\{i\}}$ is the projection corresponding to the reflection $\sigma_i$, i.e., it is the projection onto the hyper\-plane $H_i^{=0}$ along the line $L_i\cup -L_i$.
\begin{prop} \label{midelm}The linear map $\mb{mid}_{J_f}$ is a projector, i.e., $(\mb{mid}_{J_f})^2 = \mb{mid}_{J_f}$. Its kernel and its image are 
\begin{eqnarray}
  \mb{ker}\left(\mb{mid}_{J_f}\right) &=& \mb{span}\Mklz{\al_i}{i\in J_f}  , \label{kprmid}\\
   \mb{im}\left(\mb{mid}_{J_f}\right) &=&  \Mklz{\la\in \h^*}{\la(h_i)=0 \mb{ for all }i\in J_f } . \label{imprmid}
\end{eqnarray}
\end{prop}
\Proof The proposition holds trivially for $J_f=\emptyset$. Let $J_f\neq\emptyset$. The projector property $(\mb{mid}_{J_f})^2 = \mb{mid}_{J_f}$ is easy to check.
%
Since $J_f$ is a union of components of finite type, the submatrix $(a_{ij})_{i,\,j\in J_f}$ of $A$ is invertible. This implies that $(h_i)_{i\in J_f}$ as well as $(\al_i)_{i\in J_f}$ are linearly independent. Hence the dimensions of the linear spaces in (\ref{kprmid}) and (\ref{imprmid}) on the right add up to the dimension of $\h^*$. Therefore, it is sufficient to show only the inclusions ``$\supseteq$" in  (\ref{kprmid}) and (\ref{imprmid}).
Let $i\in J_f$. Then
\begin{equation*}
   \sum_{\sigma\in{\cal W}_{J_f}}\sigma\al_i=  \sum_{\sigma\in{\cal W}_{J_f}}(\sigma\sigma_i)\al_i =
 - \sum_{\sigma\in{\cal W}_{J_f}}\sigma\al_i
\end{equation*}
shows that $\mb{mid}_{J_f}(\al_i)=0$. Thus, the inclusion ``$\supseteq$" holds in (\ref{kprmid}). Since $\We_{J_f}$ acts trivially on the linear space in  (\ref{imprmid}) on the right, the inclusion ``$\supseteq$" holds in (\ref{imprmid}).
\qed

Let $p(t)$ be a statement depending on elements $t\in\R^+$. We say that $p(t)$ holds for all sufficiently small $t\in\R^+$ if there exists $t_0\in\R^+$ such that $p(t)$ holds for all $t\in\R^+$, $t\leq t_0$.

For $\emptyset\neq J_f\subseteq I$ such that $(J_f)^0 = J_f$ it follows from the classification of the generalized Cartan matrices, see Theorem 4.3 (Fin) of \cite{K}, that there exists $\gamma\in \sum_{i\in J_f}\R^+ \al_i$ such that $\gamma(h_i)>0$ for all $i\in J_f$. For $J_f=\emptyset$ this condition is satisfied by definition only for $\gamma=0$.

Formula (\ref{midsurj0}) of the following lemma has been given for $J=\emptyset$ and symmetric linear Coxeter systems as Theorem 2.2.3 in \cite{Kr}. Our proof is quite different from that in \cite{Kr}.
\begin{lem}\label{midpoint} Let $J\subseteq I$ be facial. Let $J_f\subseteq J^\bot$ such that $(J_f)^0 = J_f$. Then we have 
\begin{equation}\label{midsurj0}
      \mb{mid}_{J_f}\left( F_J\right)= F_{J\cup J_f}.
\end{equation}
Furthermore, if $\la\in F_{J\cup J_f}$ and $\gamma\in \sum_{i\in J_f}\R^+ \al_i$ such that $\gamma(h_i)>0$ for all $i\in J_f$, then we have
\begin{equation}\label{midsurj}
     \la+t\gamma\in F_J \quad \mb{ and }\quad \mb{mid}_{J_f}\left(\la+t \gamma\right) =   \la
\end{equation}
for all sufficiently small $t\in\R^+$.
\end{lem}
\Proof (a) Let $\la\in F_J$. We show that $\mb{mid}_{J_f}(\la)\in F_{J\cup J_f}$:
For all $i\in J_f$ we have
\begin{equation*}
    -\mb{mid}_{J_f}(\la)(h_i)=\mb{mid}_{J_f}(\la)(\sigma_i h_i) =  \frac{1}{|{\cal W}_{J_f}|}\sum_{\sigma\in{\cal W}_{J_f}}(\sigma_i^{-1}\sigma\la)(h_i)=\mb{mid}_{J_f}(\la)(h_i),
\end{equation*}
which shows that $\mb{mid}_{J_f}(\la)(h_i)=0$.

It is $J_f\subseteq J^\bot$. For all  $i\in J$ we find
\begin{equation*}
   \mb{mid}_{J_f}(\la)(h_i)=   \frac{1}{|{\cal W}_{J_f}|}\sum_{\sigma\in{\cal W}_{J_f}}\la(\underbrace{\sigma^{-1}h_i}_{=h_i}) =\la(h_i)=0 .
\end{equation*}

Let $i\in I\setminus (J\cup J_f)$. If $\sigma\in\We_{J_f}$ then from Proposition 2.15 of \cite{AB} we obtain $l(\sigma^{-1}\sigma_i)>l(\sigma^{-1})$, and from Theorem \ref{ThHN} we get $\sigma^{-1}h_i\in\sum_{j\in I} \R^+_0 h_j$. We conclude that
\begin{equation*}
   \mb{mid}_{J_f}(\la)(h_i)=   \frac{1}{|{\cal W}_{J_f}|}\Big(\underbrace{\la(h_i)}_{>0} + \sum_{\sigma\in{\cal W}_{J_f},\,\sigma\neq 1}\underbrace{\la(\sigma^{-1}h_i)}_{\geq 0} \Big)  > 0 . 
\end{equation*}

(b) The equation in (\ref{midsurj}) holds because of $\gamma\in \mb{ker}(\mb{mid}_{J_f})$ and $\la\in  \mb{im}(\mb{mid}_{J_f})$. 
If $i\in J$ then $(\la +  t\gamma)(h_i)= \la(h_i)=0$.
If $i\in J_f$ then $(\la +  t\gamma)(h_i)= t\gamma(h_i)>0$ for all $t\in\R^+$.
If $i\in I\setminus (J\cup J_f)$ then $(\la +  t\gamma)(h_i)= \la(h_i) + t\gamma(h_i)$ and $\la(h_i)>0$. Hence for all sufficiently small $t\in\R^+$ we have
\begin{equation*}
     (\la +  t\gamma)(h_i)= \la(h_i) + t\gamma(h_i) >0
\end{equation*}
for all $i\in I\setminus (J\cup J_f)$. 
\qed

The following corollary has been given for free integral linear Coxeter systems by P. Slodowy in \cite{Sl}, Chapter 6.2 as Corollary 2. The author could not understand the proof and gave an alternative proof for minimal free integral  linear Coxeter systems, compare Theorem 4.3 of \cite{M}, which can be easily adapted to the case of general root bases. We give here a short proof with Lemma \ref{midpoint}, similar to the proof of the corresponding result for $X=R(\emptyset)$ in  \cite{Kr}, Corollary 2.2.5.

\begin{cor} \label{riRTht} Let $\Th\subseteq I$ be special facial. Let $\Th_f\subseteq \Th^\bot$ such that $(\Th_f)^0=\Th_f$. Then we have
\begin{equation*}
 F_{\Th\cup\Th_f}\subseteq \mb{\rm ri}(R(\Th) ) .
\end{equation*}
\end{cor}
\Proof We have $\overline{F}_\Th \subseteq R(\Th)$ and $\overline{F}_\Th-\overline{F}_\Th=R(\Th)-R(\Th)$ by  Corollary \ref{specialfaceX} and Proposition \ref{rihfc}. Hence $  F_\Th=\mb{ri}(\overline{F}_\Th)\subseteq \mb{ri}(R(\Th)) $. Since $\We_{\Th_f}\subseteq N_{\mathcal W}(R(\Th))$ by Corollary \ref{stabRTh}, we get
\begin{equation*} 
   \sigma F_\Th\subseteq \sigma  \, \mb{ri}(R(\Th)) = \mb{ri}(\sigma R(\Th)) = \mb{ri}(R(\Th)) 
\end{equation*}
for all $\sigma\in \We_{\Th_f}$. With Lemma \ref{midpoint} we find $ F_{\Th\cup\Th_f} = \mb{mid}_{\Th_f}(F_\Th) \subseteq \mb{ri}(R(\Th))$. 
\qed

For free integral  linear Coxeter systems  the next theorem coincides in essential with Corollary 3 of Chapter 6.2 of \cite{Sl}, where only exposed faces have been considered. For minimal free integral linear Coxeter systems it can be found in \cite{M}, Corollaries 4.4 - 4.6. For symmetric linear Coxeter systems it can be obtained as a combination of a part of Theorem 10.3 and Lemma 10.9 of \cite{D}.
\begin{thm}\label{FXf}  (a) Let $\Th\subseteq I$ be special facial. The relative interior of $R(\Th)$ is given by 
\begin{equation*}
 \mb{\rm ri}(R(\Th))  = \,\We_{\Th^\bot} \bigcup_{\Th_f\subseteq \Th^{\bot}\atop (\Th_f)^0=\Th_f} F_{\Th\cup\Th_f} .
\end{equation*}

(b) Every $\We$-orbit of $\Fa{X}$ contains exactly one of the faces $R(\Th)$, $\Th$ special facial. A face $R\in \Fa{X}$ is of the form $R(\Th)$, $\Th$ special facial, if and only if $\mb{\rm ri}(R)\cap \Cq\neq \emptyset$. 
\end{thm}
\Proof To (b): Let $R\in\Fa{X}$ and $\la\in\mb{ri}(R)$. Then there exist $\sigma\in\We$ and $J\subseteq I$ facial such that
\begin{equation}\label{fullFaX}
  \la\in\sigma F_J=\sigma F_{J^\infty\cup J^0}\subseteq \sigma\, \mb{ri}(R(J^\infty) ) = \mb{ri}(\sigma R(J^\infty)) .
\end{equation}
Here $J^\infty$ is facial by Theorem \ref{fsp}, and the inclusion holds by Corollary \ref{riRTht}. Since the relative interiors of the faces of $X$ are a partition of $X$, we conclude that $R=\sigma R(J^\infty)$. From Proposition \ref{Rspecialinc} we find that $R(J^\infty)$ is uniquely determined.

If there exists $\la\in\mb{ri}(R)\cap \overline{C}$ we may set $\sigma=1$ in (\ref{fullFaX}) and get $R=R(J^\infty)$. From Corollary \ref{riRTht} we obtain $\mb{ri}(R(\Th))\cap\overline{C}\supseteq F_\Th\neq\emptyset$ for all $\Th$ special facial.

To (a): The inclusion ``$\supseteq$" follows from Corollary \ref{riRTht} and $N_{\mathcal W}(R(\Th))=\We_{\Th\cup\Th^\bot}$. To show the reverse inclusion let $\la\in\mb{ri}(R(\Th))$. As in the proof of (b) we find that there exist $\sigma\in\We$ and $J\subseteq I$ facial such that $R(\Th)=\sigma R(J^\infty)$. From Proposition \ref{Rspecialinc} we get $J^\infty=\Th$ and $\sigma\in\We_{\Th\cup\Th^\bot}$. Hence 
\begin{equation*}
\la\in \sigma F_J\in\We_{\Th\cup\Th^\bot} F_{\Th\cup J^0}= \We_{\Th^\bot} F_{\Th\cup J^0}
\end{equation*}
with $J^0\subseteq\Th^\bot$, $(J^0)^0=J^0$.
\qed

The smallest face $X\cap(-X)$ of $X$ is a $\We$-invariant linear subspace of $\h^*$. The following description is obtained immediately from Proposition \ref{Rspecialinc}, Theorem \ref{FXf} (b), and Corollary \ref{specialfaceX}.
\begin{cor}\label{descsf} We have
\begin{equation*}
   X\cap(-X) =  R(I^\infty) = \{ \la \in \h^* \mid \la(h_i)=0 \mb{ for all } i\in I^\infty \} .
\end{equation*}
In particular, $X=\h^*$ if and only if $I^\infty=\emptyset$.
\end{cor}

We define a function $\mb{red}:\We\to I$ as follows: We set $\red{1}:=\emptyset$. If $\sigma\in\We\setminus\{1\}$ and $\sigma=\sigma_{i_1}\cdots \sigma_{i_k}$ is a reduced expression with simple reflections $\sigma_{i_1},\,\ldots,\,\sigma_{i_k}$ we set
\begin{equation*}
\red{\sigma}:= \{i_1,\,\ldots,\,i_k\}\subseteq I,
\end{equation*} 
which is independent of the chosen reduced expression by a result of J. Tits, compare for example with Proposition 2.16 or Remark 2.34 of \cite{AB}.

For minimal free integral linear Coxeter systems the following theorem has been given as Theorem 4 (b) of \cite{M6}. For symmetric linear Coxeter systems it can be obtained as a combination of Proposition 10.2, Remark 10.2, Proposition 2.17, and Theorem 10.3 (c) of \cite{D}. 
\begin{thm}\label{FXlmj} Let $\tau\in\mb{}^{\Th_1\cup\Th_1^\bot}\We^{\Th_2\cup\Th_2^\bot}$, and $\Th_1,\Th_2\subseteq I$ be special facial. Then $\Th_1\cup\Th_2\cup \text{\rm red}(\tau)$ is a special set and 
\begin{eqnarray}
   R(\Th_1)\cap \tau R(\Th_2)   &=&  R (   ( \Th_1\cup\Th_2\cup \text{\rm red}(\tau) )^{fc}  ) , \label{int}\\
   R(\Th_1)\vee\tau R(\Th_2)   &=&  R((\Th_1\cap \tau\Th_2)^\infty). \label{join}
\end{eqnarray}
\end{thm}
\begin{rem} The lattice meet and lattice join of two arbitrary faces of $X$ can be easily reduced to the formulas (\ref{int}), (\ref{join}) by Theorem \ref{FXf} (b) and Corollary \ref{stabRTh}.
\end{rem}

\Proof The proof that $\Th_1\cup\Th_2\cup \red{\tau}$ is a special set is identical to the corresponding part of the proof of Theorem 4 (b) in \cite{M6}. The proofs of (\ref{int}) and (\ref{join}) are similar to the corresponding parts of the proof of Theorem 4 (b) in \cite{M6}, the faces $R(\Th)$ now parametrized by special facial sets $\Th$ instead of special sets $\Th$. For (\ref{join}), Theorem \ref{JsJ} is used in addition. 
\qed

%
%
%
\subsection{\label{SubsectionFI}The faces of the imaginary cone}
%
%
%
The imaginary cone has been introduced for symmetric linear Coxeter systems by M. Dyer in Section 3 of \cite{D}, generalizing the convex cone generated by the imaginary roots of a Kac-Moody algebra as described for example in Chapter 5 of \cite{K}. Its definition can be adopted for linear Coxeter systems without change. In this section we determine the faces of the imaginary cone of a linear Coxeter system and describe the lattice operations of its face lattice, generalizing the results obtained in Section 10 of \cite{D}. Our approach is motivated by the theory of invariant convex cones in Section \ref{SecInvsub}.

Let $(\h,(h_i)_{i\in I},(\al_i)_{i\in I})$ be a root base of the linear Coxeter system $(\h,\,(H_i, \,L_i)_{i\in I})$. Set
\begin{equation}\label{defKim}
  K := (\sum_{i\in I}\R^+_0 h_i ) \cap(-\overline{C}^\vee) ,
\end{equation}
which is a finitely generated closed convex cone. The proof of the following proposition, in which Proposition \ref{fcsdiff} is used, is similar to the proof of Proposition 3.2 in \cite{D}.

\begin{prop}\label{defimcone}  We get a convex cone, which is called the imaginary cone, by
\begin{equation}\label{defimcone1}
 Z:=\bigcup_{\sigma\in {\cal W}} \sigma K = (-X^\vee)\cap \bigcap_{\sigma\in{\cal W}} \sigma (\,\sum_{i\in I}\R^+_0 h_i \,) .
\end{equation}
\end{prop}
\begin{rem} Recall from Propositions \ref{dlinCox1} and \ref{dlinCox2} that $X^\vee$ may be interpreted as a Tits cone with $\We_{I_1}$ as Coxeter group, where $I_1$ is defined in (\ref{indexexc}). We may replace $\We$ by $\We_{I_1}$ in (\ref{defimcone1}). It is not possible to replace $I$ by $I_1$ in (\ref{defKim}), (\ref{defimcone1}). For the example given in Remark \ref{exaffkm} we have $Z=K= \R^+_0(h_1+h_2)$, but $\sum_{i\in I_1}\R^+_0 h_i=\{0\}$.
\end{rem}

Recall that the faces of $\sum_{i\in I}\R^+_0 h_i$ are given by $\sum_{i\in \Th}\R^+_0 h_i$, $\Th \subseteq I$ facial. Its relative interiors give a partition of $\sum_{i\in I}\R^+_0 h_i$, i.e.,
\begin{equation*}
     \sum_{i\in I}\R^+_0 h_i =  \dot{ \bigcup_{\Th \text{ facial}} }  \mb{ri} (\,\sum_{i\in \Th}\R^+_0 h_i \,) =  \dot{ \bigcup_{\Th \text{ facial} } } \; \sum_{i\in \Th}\R^+ h_i ,
\end{equation*}
which in turn induces a partition of $K$. Thus
\begin{equation}\label{part K in rel int}
  K=\dot{\bigcup_{\Th \text{ special facial} }} K(\Th)
\end{equation}
where by Proposition \ref{Ksp} the convex cone 
\begin{equation}
  K(\Th) :=  (\sum_{i\in \Th}\R^+ h_i ) \cap(-\overline{C}^\vee) \neq\emptyset  
\end{equation}
for all $\Th$ special facial. Only $K(\emptyset)=\{0\}$ contains the zero. Moreover, for $\Th$ special facial, 
\begin{equation}\label{clKThdef}
  \overline{ K(\Th) } :=  (\sum_{i\in \Th}\R_0^+ h_i ) \cap(-\overline{C}^\vee)    = \dot{\bigcup_{\Xi \text{ special facial},\; \Xi\subseteq \Th }} K(\Xi)
\end{equation}
is a face of $K$. In general, there also exist faces of $K$ not of this form, but these are not of interest for the following considerations.  
Note that $\overline{K(\Th)}$ is the closure of $K(\Th)$, which can be seen as follows: $\overline{K(\Th)}$ is partitioned by the relative interiors of its faces. From (\ref{clKThdef}) we find that $\overline{K(\Th)}$ is closed and $\mb{ri}(\overline{K(\Th)})\subseteq K(\Th)\subseteq \overline{K(\Th)}$.
Now we use Theorem 6.3 of \cite{Ro}.

The following theorem introduces some faces of the imaginary cone. For symmetric linear Coxeter systems it can be obtained as a combination of Lemma 10.4 (b) and Lemma 10.9 of \cite{D}.
\begin{thm}\label{specialfaceZ} Let $\Th\subseteq I$ be special facial. Then 
\begin{equation*}
  F(\Th):=\We_\Th  \overline{K(\Th)}
\end{equation*}
is an exposed face of $Z$, i.e., for every element $\la\in F_{\Th\cup\Th_f}$ where $\Th_f\subseteq\Th^\bot$ such that $(\Th_f)^0=\Th_f$ we have
\begin{eqnarray}
   Z     & \subseteq & \{ h\in \h \mid \la(h)\geq 0\} ,\label{Zexp1}\\
  F(\Th) &  =        & Z\cap \{ h\in \h \mid \la(h)=0\}.\label{Zexp2}
\end{eqnarray}
\end{thm}
\begin{rem} In particular, $F(I^\infty)=Z$, which follows from the theorem for $\la=0\in F_{I^\infty\cup I^0}$, and $F(\emptyset)=\{0\}$.
\end{rem}
\Proof The elements of $Z$ are of the form $\sigma h$ with $\sigma\in\We$, $h\in K(\Xi)$, $\Xi$ special facial. From Corollary \ref{specialfaceX} we obtain
\begin{eqnarray*}
  \la(\sigma h)= (\sigma^{-1}\la) (h) \, \left\{\begin{array}{lcl}
   =0 & \mb{ if } & \sigma^{-1}\la\in R(\Xi), \\
    >0 &\mb{ if } & \sigma^{-1}\la\notin R(\Xi) . 
\end{array}\right.
\end{eqnarray*}
In particular, this shows (\ref{Zexp1}).
From Corollary \ref{riRTht} we find $\la\in \mb{ri}(R(\Th))$. Hence, $\sigma^{-1}\la\in R(\Xi)$ if and only if $\sigma^{-1}R(\Th)\subseteq R(\Xi)$. 
By Proposition \ref{Rspecialinc} this is equivalent to $\Th\supseteq\Xi$ and $\sigma\in\We_{\Th}\We_{\Xi^\bot}$, from which we get
\begin{equation*}
\sigma h\in \We_{\Th}\We_{\Xi^\bot} K(\Xi)=\We_{\Th} K(\Xi)\subseteq \We_{\Th} \overline{K(\Th)}.
\end{equation*} 
This shows the inclusion ``$\supseteq$" in (\ref{Zexp2}). 
To show the reverse inclusion let $\sigma h\in F(\Th)$ where $\sigma\in\We_\Th$, $h\in\overline{K(\Th)}$. Since $\sigma\in Z_{\mathcal W}(F_{\Th\cup\Th_f})$ and $h\in\sum_{i\in\Th}\R^+_0 h_i$, we find 
\begin{equation*}
  \la(\sigma h)= (\sigma^{-1}\la)(h)=\la(h)=0.
\end{equation*}
\qed

In the preceding theorem the face $F(\Th)$ of the imaginary cone $Z$ has been described as an exposed face, dual to the face $R(\Th)$ of the Tits cone $X$. The imaginary cone $Z$ is contained in $\sum_{i\in I}\R^+_0 h_i$, and $F(\Th)$ can also be obtained as an intersection of $Z$ with a face of $\sum_{i\in I}\R^+_0 h_i$. For symmetric linear Coxeter systems this follows from Lemma 3.4 of \cite{D}.
\begin{cor}\label{FZdiffint} Let $\Th\subseteq I$ be special facial. Then we have 
\begin{equation*}
       F(\Th) =Z\cap \sum_{i\in\Th}\R^+_0 h_i.
\end{equation*}
\end{cor}
\Proof Since $Z\subseteq \sum_{i\in I}\R^+_0 h_i$ by Proposition \ref{defimcone}, we find from Theorem \ref{specialfaceZ} for $\la\in F_{\Th}$ that
\begin{equation*}
   F(\Th) = Z \cap \{\,h\in \h \mid \la(h)=0 \,\} = Z \cap \{\,h\in  \sum_{i\in I}\R^+_0 h_i  \mid \la(h)=0 \,\} =  Z\cap \sum_{i\in\Th}\R^+_0 h_i. 
\end{equation*}
\qed

For symmetric linear Coxeter systems the following Proposition has been given as Lemma 10.4 (d) and (e) of \cite{D}.
\begin{prop} \label{Fspecialinc} Let $\sigma,\sigma'\in\We$, and $\Th,\Th'\subseteq I $ be special facial. Then:
\begin{equation}\label{sub}
   \sigma' F(\Th') \supseteq \sigma F(\Th)\quad \iff \quad \Th' \supseteq \Th \mb{ and } \sigma^{-1}\sigma' \in \We_{\Th^\bot}\We_{\Th'}.  
\end{equation}
In particular, the stabilizer of $F(\Th)$ as a whole is given by
\begin{equation} \label{sub2}
    N_{\mathcal W}(F(\Th))  = \We_{\Th\cup\Th^\bot}  .
\end{equation}
\end{prop}
\Proof (\ref{sub2}) is a special case of (\ref{sub}). We first show the implication ``$\Leftarrow$" in (\ref{sub}). Decompose $\sigma^{-1}\sigma'=ab$ with $a\in\We_{\Th^\bot}$ and $b\in\We_{\Th'}$. Since $\Th\subseteq\Th'$, we find
\begin{equation*}
  a^{-1}F(\Th)=a^{-1}\We_{\Th}\overline{K(\Th)}= \We_{\Th}\,a^{-1}\overline{K(\Th)} = \We_{\Th} \overline{K(\Th)}   \subseteq \We_{\Th'} \overline{K(\Th')} =b\We_{\Th'} \overline{K(\Th')} = b F(\Th').
\end{equation*}
We conclude that $\sigma F(\Th)\subseteq  \sigma' F(\Th')$. 

Now we show the implication ``$\Rightarrow$" in (\ref{sub}). Let $h\in K(\Th) = \sum_{i\in \Th}\R^+ h_i  \cap(-\overline{C}^\vee) $. Then $\sigma h\in\sigma F(\Th)\subseteq \sigma' F(\Th')$, and from Theorem \ref{specialfaceZ} we obtain for $\la\in F_{\Th'}$ that
\begin{equation*}
    0=\la((\sigma')^{-1}\sigma h)= (\sigma^{-1}\sigma'\la)(h).
\end{equation*}
From Theorem \ref{specialfaceX} we get $\sigma^{-1}\sigma'\la\in R(\Th)$. Since $\la\in\mb{ri}(R(\Th'))$ by Corollary \ref{riRTht}, this is equivalent to $\sigma^{-1}\sigma' R(\Th')\subseteq R(\Th)$. By Proposition \ref{Rspecialinc} we get $ \Th' \supseteq \Th$ and $ \sigma^{-1}\sigma' \in \We_{\Th^\bot}\We_{\Th'}$.
\qed

Up to this point we used the results for the faces of the Tits cone $X$ to obtain dually results for the faces of the imaginary cone $Z$. To show
\begin{equation*}
    \Fa{Z}   = \{\sigma F(\Th)  \mid  \sigma\in\We,\,\Th \mb{ special facial}  \}
\end{equation*}
is more complex. Since $-Z$ is an invariant convex subcone of the Tits cone $X^\vee$, we can use an approach motivated by the theory of invariant convex cones in Section \ref{SecInvsub}. Instead of working with the system of facial sets for $\overline{C}^\vee$, which is in general different from that for $\overline{C}$, taking into consideration also Propositions \ref{dlinCox1}, \ref{dlinCox2}, and \ref{dlinCox3}, we transfer some results for free linear Coxeter systems to general linear Coxeter systems by morphisms.

We set $S^\vee_\emptyset:=\{0\}\subseteq\h$. If $\Th\subseteq I$ is of affine type, then by Theorem 4.3 (Aff) of \cite{K} the kernel of $A_\Th^T$ is one-dimensional, and there exists $k\in (\R^+)^\Th$ such that $A_\Th^T k=0$. We define
\begin{equation*}
  D_\Th^\vee:= \R^+\sum_{i=1}^n k_i h_i \quad \mb{ and } \quad  S_\Th^\vee:= \R \sum_{i=1}^n k_i h_i ,
\end{equation*}
which are independent of the chosen root base of the linear Coxeter system. If $\emptyset\neq \Th\subseteq I$ has only affine connected components, say $\Th_1$,\ldots, $\Th_m$, we set
\begin{equation*}
   S_\Th^\vee := \sum_{i=1}^m S_{\Th_i}^\vee.
\end{equation*}

For $J_f\subseteq I$ such that $(J_f)^0 = J_f$ we define a linear projector $\mb{mid}_{J_f}^\vee: \h \to \h $ by
\begin{equation*}
   \mb{mid}_{J_f}^\vee (h):=   \frac{1}{|{\cal W}_{J_f}|}\sum_{\sigma\in{\cal W}_{J_f}}\sigma h,  \quad h \in \h . 
\end{equation*}
It has similar properties as the projector $ \mb{mid}_{J_f}:\h^*\to\h^*$, to which it is dual.

Recall from Proposition \ref{facial-h-free} that for a free linear Coxeter system every subset of $I$ is facial. Our first aim is to describe for $\Th\subseteq I$ special the position of the nonempty set $K(\Th)= (\sum_{i\in\Th}\R^+ h_i)\cap (-\overline{C}^\vee)$ in $-\overline{C}^\vee$ in detail, which will be reached in Theorem \ref{mainKTh}.
For symmetric linear Coxeter systems part (a) of the following theorem has been given in Section 8.6 of \cite{D}, and the property $K(\Th) \cap (-F_{\Th^\bot}^\vee)\neq\emptyset$ of part (b) follows from  Lemma 10.4 (b) of \cite{D}.

\begin{thm}\label{mainThK0} Let the linear Coxeter system be free. Let $\emptyset\neq \Th\subseteq I$ be special and connected. 
\begin{itemize}
\item[(a)] If $\Th$ is of affine type then $K(\Th)= D_\Th^\vee\subseteq  -F_{\Th\cup\Th^\bot}^{\vee}$.
\item[(b)] If $\Th$ is of indefinite type then
\begin{equation}\label{kkfacetseq1}
  K(\Th) = \dot{\bigcup_{\Th_f\subseteq \Th\atop (\Th_f)^0=\Th_f} }  \underbrace{ K(\Th) \cap (-F_{\Th_f\cup \Th^\bot}^\vee)  }_{\neq\emptyset} ,
\end{equation}
and for $\Th_f\subseteq \Th$ such that $(\Th_f)^0=\Th_f$ we have
\begin{equation}\label{kkfacetseq2}
   \mb{mid}_{\Th_f}^{\,\vee} \left(  K(\Th) \cap  \left(-F_{\Th^\bot}^\vee\right)  \right) = K(\Th) \cap  (-F_{\Th_f\cup \Th^\bot}^\vee) .
\end{equation}
The relative interiors of $\overline{K(\Th)}$ and $K(\Th)$ are given by
\begin{equation}\label{kkfacetseq3}
    \mb{\rm ri}(\overline{K(\Th)}) = \mb{\rm ri}(K(\Th)) = K(\Th) \cap  (-F_{\Th^\bot}^\vee) .
\end{equation}
The linear hulls of  $\overline{K(\Th)}$ and $K(\Th)$ are given by
\begin{equation}\label{kkfacetseq4}
   \overline{K(\Th)} -\overline{K(\Th)} = K(\Th) - K(\Th)=\mb{span}\Mklz{ h_i }{ i\in\Th } .
\end{equation}
\end{itemize}
\end{thm}
\Proof It is easy to see that for $J\subseteq I$, and $h \in \sum_{i\in J}\R^+ h_i  $ with $\al_j(h)\leq 0$ for all $j\in J$ we have
\begin{equation}\label{el k in facet}
   h  \in - F^\vee_{J_f\cup J^\bot} \quad \mb{ with }\quad J_f: =  \Mklz{j\in J}{\al_j (h)=0} \subseteq J ,
\end{equation}
compare also the proof of Theorem \ref{fac centralizer} (a). 

Now part (a) follows from the definition of $K(\Th)$, Theorem 4.3 (Aff) of \cite{K}, and (\ref{el k in facet}). To (b): 
We first show the inclusion ``$\subseteq$" in (\ref{kkfacetseq1}). We obtain from (\ref{el k in facet}) that for $h\in K(\Th)=(\sum_{i\in \Th} r_i h_i )\cap (-\overline{C}^\vee)$ we have $h\in - F_{\Th_f\cup\Th^\bot}^\vee$ with $\Th_f\subseteq \Th$. From Theorem 4.3 (Ind) of \cite{K} we get $\Th_f\neq\Th$, and as in \cite{K}, Exercise 5.9, we find $\Th_f^0=\Th_f$.

From the definition of $K(\Th)$ and from the inclusion ``$\subseteq$" in  (\ref{kkfacetseq1}) we get
\begin{equation*}
   K(\Th) = (\sum_{i\in\Th}\R^+ h_i )\cap (-\overline{C}^\vee) = (\sum_{i\in\Th}\R^+ h_i )\cap( -\overline{F}_{\Th^\bot}^\vee).
\end{equation*}
By Theorem 4.3 (Ind) of \cite{K} there exists $h\in\sum_{i\in \Th}\R^+ h_i $ such that $\al_j (h)=\sum_{i\in\Th}r_i a_{ij}<0$ for all $j\in \Th$. From (\ref{el k in facet}) we obtain $h\in -F_{\Th^\bot}^\vee$. Therefore,
\begin{equation*}
   \mb{ri}  (\,\sum_{i\in\Th}\R^+ h_i \, )  \cap \mb{ri}  ( -\overline{F}_{\Th^\bot}^\vee )  
     = (\sum_{i\in\Th}\R^+ h_i )\cap (-F_{\Th^\bot}^\vee ) = K(\Th)\cap (-F_{\Th^\bot}^\vee )\neq \emptyset .
\end{equation*}
By Theorem 6.3 and 6.5 of \cite{Ro} we find $\mb{ri}(\overline{K(\Th)}) = \mb{ri}(K(\Th)) = K(\Th)\cap (-F_{\Th^\bot}^\vee )$.

Now we show for $\Th_f\subseteq \Th$ such that $(\Th_f)^0=\Th_f$ the inclusion
\begin{equation}\label{zwmidvee}
   \mb{mid}_{\Th_f}^\vee \left(    K(\Th) \cap ( -F_{\Th^\bot}^\vee) \right)  \subseteq (\sum_{i\in \Th}\R^+ h_i ) \cap ( -F^\vee_{\Th_f\cup\Th^\bot} ) = K(\Th) \cap  (-F_{\Th_f\cup \Th^\bot}^\vee) ,
\end{equation}
which in particular implies $K(\Th) \cap (-F^\vee_{\Th_f\cup\Th^\bot} )\neq\emptyset$.
Let $h\in  K(\Th) \cap ( -F_{\Th^\bot}^\vee )\subseteq \sum_{i\in \Th}\R^+ h_i $. We find from Proposition \ref{fcsdiff} that
\begin{equation*}
    \sigma h\in h+(  \sum_{i\in I}\R^+_0 h_i \cap \sum_{i\in \Th_f}\R h_i ) = h+\sum_{i\in \Th_f}\R^+_0 h_i \subseteq \sum_{i\in\Th}\R^+h_i
\end{equation*}
for all $\sigma\in \We_{\Th_f}$, from which we get $\mb{mid}_{\Th_f}^\vee (h)\in  \sum_{i\in\Th}\R^+h_i$. By the analogue of Lemma \ref{midpoint} for $\mb{mid}_{\Th_f}^\vee$ we obtain $\mb{mid}_{\Th_f}^\vee (h)\in  -F^\vee_{\Th_f\cup\Th^\bot}$. 

We next show that in (\ref{zwmidvee}) we have equality. This holds trivially for $\Th_f=\emptyset$. Let $\Th_f\neq\emptyset$ and $h\in  K(\Th) \cap  (-F_{\Th_f\cup \Th^\bot}^\vee)$. It follows from Theorem 4.3 (Fin) of \cite{K} that there exists $k\in\sum_{i\in \Th_f}\R^+ h_i$ such that $\al_i(k)>0$ for all $i\in \Th_f$. By the analogue of Lemma \ref{midpoint} for $\mb{mid}_{\Th_f}^{\vee}$ we obtain
\begin{equation*}
    h  - t k\in - F_{\Th^\bot}^\vee \quad\mb{ and }\quad \mb{mid}_{\Th_f}^\vee (h - t k ) =  h
\end{equation*}
for all sufficiently small $t\in\R^+$. Since $ h \in \sum_{i\in\Th} \R^+ h_i$ and $\Th_f\subseteq \Th$, we find $h - t k \in \sum_{i \in\Th}\R^+ h_i$ for all sufficiently small $t\in\R^+$. Thus $h - tk\in  \sum_{i \in\Th}\R^+ h_i \cap (- F_{\Th^\bot}^\vee )= K(\Th) \cap (- F_{\Th^\bot}^\vee )$ for all sufficiently small $t\in\R^+$.

From the definitions of $K(\Th)$, $\overline{K(\Th)}$ we get
\begin{equation*}
   K(\Th)-K(\Th)\subseteq \overline{K(\Th)}-\overline{K(\Th)}\subseteq \mb{span}\Mklz{h_i}{i\in\Th}.
\end{equation*}
To show the reverse inclusions let $h\in  \mb{span}\Mklz{h_i}{i\in\Th}$. Choose $k\in K(\Th)\cap (-F_{\Th^\bot}^\vee)\subseteq \sum_{i\in \Th}\R^+ h_i $. Then
\begin{equation*}
 h= h+t k- t k \quad  \mb{ with } \quad h+tk\in\sum_{i\in \Th}\R^+ h_i \quad\mb{ and } \quad -tk\in - K(\Th)
\end{equation*}
for all sufficiently big $t\in\R^+$. Clearly, if $j\in I\setminus\Th$ then $\al_j(h+tk)\leq 0$ for all sufficiently big $t\in\R^+$. If $j\in\Th$ then
\begin{equation*}
 \al_j(h+tk)=\al_j(h) +t \underbrace{\al_j(k)}_{<0} <0
\end{equation*}
for all sufficiently big $t\in\R^+$. Hence $h+tk\in (\sum_{i\in \Th}\R^+ h_i ) \cap(-\overline{C}^\vee)= K(\Th)$ for all sufficiently big $t\in\R^+$.
\qed

It is trivial to show the following:
\begin{prop}\label{KThZer}  Let the linear Coxeter system be free. Let $\emptyset\neq\Th\subseteq I$ be special, and let $\Th=\Th_1\cup\cdots \cup \Th_m$ be its decomposition in connected components. Then we have
\begin{equation*}
    K(\Th) =  K(\Th_1) + \cdots + K(\Th_m)  . 
\end{equation*}
\end{prop}

From the previous proposition and Theorem \ref{mainThK0} we obtain the position of the nonempty set $K(\Th)$ in $-\overline{C}^\vee$ also for special sets $\Th$, which are not connected. For symmetric linear Coxeter systems the property $K(\Th) \cap (-F_{\Th^{aff}\cup\Th^\bot}^\vee)\neq\emptyset$ of the following theorem follows from  Lemma 10.4 (b) of \cite{D}.
\begin{thm} \label{mainKTh}  Let the linear Coxeter system be free. Let $\Th\subseteq I$ be special. Then 
\begin{equation}\label{kfacetseq1}
  K(\Th) = \dot{\bigcup_{\Th_f\subseteq \Th^{ind}\atop (\Th_f)^0=\Th_f} }  \underbrace{ K(\Th) \cap   (-F_{\Th^{aff}\cup \Th_f\cup \Th^\bot}^\vee )   }_{\neq\emptyset} 
\end{equation}
and for $\Th_f\subseteq \Th^{ind}$ such that $(\Th_f)^0=\Th_f$ we have
\begin{equation}\label{kfacetseq2}
    \mb{mid}_{\Th_f}^{\,\vee}\left(  K(\Th) \cap  \left(-F_{\Th^{aff}\cup \Th^\bot}^\vee\right) \right)  =  K(\Th) \cap  (-F_{\Th^{aff}\cup\Th_f\cup \Th^\bot}^\vee ) .
\end{equation}
The relative interiors of $\overline{K(\Th)}$ and $K(\Th)$ are given by
\begin{equation}\label{kfacetseq3}
  \mb{\rm ri} (\overline{K(\Th)} ) =  \mb{\rm ri}(K(\Th)) = K(\Th) \cap  \left(-F_{\Th^{aff}\cup \Th^\bot}^\vee\right) .
\end{equation}
The linear hulls of $\overline{K(\Th)}$ and $K(\Th)$ are given by
\begin{equation}\label{kfacetseq4}
   \overline{K(\Th)} -\overline{K(\Th)} = K(\Th)-K(\Th)=  S^\vee_{\Th^{aff}} + \mb{span}\Mklz{ h_i }{ i\in\Th^{ind} } .
\end{equation}
\end{thm}
\Proof Clearly, the theorem holds for $K(\emptyset)=\{0\}$. Let $\Th\neq\emptyset$, and let $\Th=\Th_1\cup\cdots\cup\Th_m$ be its decomposition in connected components. We first show that for $J_1\subseteq \Th_1$, \ldots, $J_m\subseteq \Th_m$ we have
\begin{equation}\label{KThFZer}
\left(K(\Th_1)\cap(-F^\vee_{J_1\cup\Th_1^\bot})\right)+\cdots +\left(K(\Th_m)\cap(-F^\vee_{J_1\cup\Th_m^\bot})\right) =K(\Th) \cap (-F^\vee_{J_1\cup\cdots\cup J_m\cup\Th^\bot} ),
\end{equation}
where a sum with the empty set as a summand is defined to be empty. 
We have
\begin{equation*}
    (-F^\vee_{J_1\cup\Th_1^\bot}) +\cdots + (-F^\vee_{J_1\cup\Th_m^\bot}) =  -F^\vee_L \quad \mb{ with }\quad L:=\bigcap_{i=1}^m (J_i\cup\Th_i^\bot) . 
\end{equation*}
Since $J_i\cap J_j=\emptyset$ and $J_i\subseteq\Th_j^\bot$ for $i\neq j$, we obtain
\begin{eqnarray*}
     &&L = \bigcap_{i=1}^m (J_i\cup\Th_i^\bot) = \bigcup_{T\subseteq\{1, \ldots, m\} }(\,\bigcap_{i\in T}\,J_i\cap\bigcap_{j\in \{1,\ldots, m\}\setminus T}\Th_j^\bot  \,) = (\,\bigcap_{j=1}^m\Th_j^\bot \,) \cup J_1\cup\cdots\cup J_m\\
&&\quad = J_1\cup\cdots\cup J_m\cup (\Th_1\cup\cdots\cup\Th_m)^\bot .
\end{eqnarray*}
Now the inclusion ``$\subseteq$" in (\ref{KThFZer}) follows with Proposition \ref{KThZer}.
To show the reverse inclusion let
\begin{equation*}
k\in K(\Th) \cap (-F^\vee_{J_1\cup\cdots\cup J_m\cup\Th^\bot} ) .
\end{equation*} 
According to Proposition  \ref{KThZer} we have $k=\sum_{i=1}^m k_i $ with $k_i\in K(\Th_i)$. By Theorem \ref{mainThK0} we find $k_i\in -F^\vee_{L_i\cup \Th_i^\bot}$ for some $L_i\subseteq\Th_i$. By the inclusion ``$\subseteq$" of (\ref{KThFZer}) with $J_i$ replaced by $L_i$ for all $i$ we obtain
\begin{equation*}
    k\in K(\Th) \cap (-F^\vee_{L_1\cup\cdots\cup L_m\cup\Th^\bot} ).
\end{equation*} 
Hence $L_1\cup\cdots\cup L_m\cup\Th^\bot=J_1\cup\cdots\cup J_m\cup\Th^\bot$, from which we get $L_1=J_1$, \ldots,\, $L_m=J_m$.

Now (\ref{kfacetseq1}) follows from (\ref{KThFZer}), and from Theorem \ref{mainThK0} (a) and (\ref{kkfacetseq1}) in (b).  Let $i\in\{1,\,\ldots,\,m\}$ and $k\in K(\Th_i)$. Since $\We_{\Th_f} =\prod_{j=1}^m \We_{\Th_f\cap\Th_j}$, we find 

\begin{equation*}
   \mb{mid}_{\Th_f}^\vee (k)
     =   \frac{1}{\prod_{j=1}^m |\We_{\Th_f\cap\Th_j}|} \sum_{\tau_j\in {\mathcal W}_{\Th_f\cap\Th_j}\atop j=1,\ldots, m}\underbrace{ \tau_1\cdots\tau_m k}_{=\tau_i k} 
     =  \frac{1}{|\We_{\Th_f\cap\Th_i}|} \sum_{\tau_i\in {\mathcal W}_{\Th_f\cap\Th_i} }\tau_i k =\mb{mid}_{\Th_f\cap\Th_i}^\vee (k).
\end{equation*}
Thus (\ref{kfacetseq2}) follows from (\ref{KThFZer}) and from (\ref{kkfacetseq2}) in Theorem \ref{mainThK0}  (b).

Proposition \ref{KThZer} also implies
\begin{eqnarray*}
   \mb{ri}(K(\Th)) &=&  \mb{ri}(K(\Th_1)) + \cdots + \mb{ri}(K(\Th_m) ) ,\\
    K(\Th)- K(\Th)  &=& \left( K(\Th_1) - K(\Th_1) \right)+ \cdots + \left(K(\Th_m) - K(\Th_m) \right).
\end{eqnarray*}
From Theorem  \ref{mainThK0} (a) and  (\ref{kkfacetseq3}) in (b), from (\ref{KThFZer}), and from Theorem 6.3 of \cite{Ro} we get  (\ref{kfacetseq3}). Since $K(\Th)-K(\Th)$ is a closed linear space which contains $K(\Th)$, we find $\overline{K(\Th)}-\overline{K(\Th)}=K(\Th)-K(\Th)$. From Theorem  \ref{mainThK0} (a) and  (\ref{kkfacetseq4}) in (b) we obtain (\ref{kfacetseq4}).
\qed

\begin{cor} \label{zenfree}  Let the linear Coxeter system be free. Let $\Th\subseteq I$ be special. Then
\begin{equation*} 
   Z_{\mathcal W}(F(\Th))  =  Z_{\cal W}(\overline{K(\Th)}) =\We_{\Th^{aff}\cup\Th^\bot}  .
\end{equation*}
Moreover, we have $ F(\Th) =\We_{\Th^{ind}}\overline{K(\Th)}$.
\end{cor}
\Proof From Theorem 6.3 of \cite{Ro} and Theorem \ref{mainKTh} we obtain 
\begin{equation*} 
   Z_{\cal W}(\overline{K(\Th)}) =  Z_{\cal W}(\mb{ri}(K(\Th))) = \We_{\Th^{aff}\cup\Th^\bot}.
\end{equation*} 
Hence $ F(\Th)=\We_{\Th}\overline{K(\Th)} =\We_{\Th^{ind}}\We_{\Th^{aff}}\overline{K(\Th)}=\We_{\Th^{ind}}\overline{K(\Th)}$, from which we find
\begin{equation*}
   Z_{\cal W}(F(\Th) )= \bigcap_{\sigma\in{\mathcal W}_{\Th^{ind}} }\sigma  Z_{\cal W}(\overline{K(\Th)}) \sigma^{-1} = \bigcap_{\sigma\in{\mathcal W}_{\Th^{ind}} }\sigma  \We_{\Th^{aff}\cup\Th^\bot}  \sigma^{-1} =  \We_{\Th^{aff}\cup\Th^\bot}  .
\end{equation*}
\qed

\begin{cor}\label{hullfree}  Let the linear Coxeter system be free. Let $\Th\subseteq I$ be special. Then
\begin{equation*} 
        F(\Th) -  F(\Th) = \overline{K(\Th)} - \overline{K(\Th)} = S_{\Th^{aff}}^\vee + \mb{span}\Mklz{ h_i }{i\in\Th^{ind}}  .
\end{equation*}
\end{cor}
\Proof By Theorem \ref{mainKTh} we find
\begin{equation*}
\pm \overline{K(\Th)}\subseteq \overline{K(\Th)} - \overline{K(\Th)} =  S_{\Th^{aff}}^\vee + \mb{span}\Mklz{ h_i }{ i\in\Th^{ind} } \subseteq F(\Th)- F(\Th).
\end{equation*}
Since $F(\Th)= \We_{\Th^{ind}}\overline{K(\Th)}$ by Corollary \ref{zenfree}, it remains to show that $ \overline{K(\Th)} - \overline{K(\Th)}$ is invariant under $\We_{\Th^{ind}}$. 
Clearly, $\mb{span}\Mklz{ h_i }{ i\in\Th^{ind} }$ is invariant, and $S_{\Th^{aff}}^\vee$ is invariant because of $\Th^{aff}\subseteq (\Th^{ind})^\bot$. 
\qed

The proof of part (a) of the following corollary is similar to the proof of Corollary \ref{riRTht}. Part (b) follows from part (a) and (\ref{kfacetseq1}).
\begin{cor}\label{prikZf} Let the linear Coxeter system be free. Let $\Th\subseteq I $ be special.
\begin{itemize} 
\item[(a)] For $\Th_f\subseteq \Th^{ind}$ such that $(\Th_f)^0=\Th_f$ we have
\begin{equation*}
        K(\Th) \cap    (-F_{\Th^{aff}\cup\Th_f\cup \Th^\bot}^\vee )  \subseteq \mb{\rm ri}(F(\Th)) .
\end{equation*}
\item[(b)] We have $K(\Th)\subseteq \mb{\rm ri}(F(\Th))$.
\end{itemize}
\end{cor}

We now use morphisms of root bases to transfer Corollaries \ref{zenfree}, \ref{hullfree}, and \ref{prikZf} (b) to general linear Coxeter systems. 
\begin{prop}\label{morphZ1} Let $(\h,(h_i)_{i\in I},(\al_i)_{i\in I})$ and  $(\h', (h_i')_{i\in I}, (\al_i')_{i\in I})$ be root bases of linear Coxeter systems. Let  $\phi:\h'\to\h$ be a morphism of root bases. If $\Th\subseteq I$ is special facial for $(\h,(h_i)_{i\in I},(\al_i)_{i\in I})$, then it is also special facial for $(\h', (h_i')_{i\in I}, (\al_i')_{i\in I})$, and we have
\begin{equation*}
    \phi(K'(\Th) )=K(\Th)\quad \mb{ and }\quad  \phi(\overline{ K'(\Th) } )=\overline{ K(\Th) }\quad \mb{ and }\quad \phi(F'(\Th) )=F(\Th) .
\end{equation*}
\end{prop}
\Proof If $\Th$ is special facial for $(\h,(h_i)_{i\in I},(\al_i)_{i\in I})$, then it is also special facial for $(\h',(h'_i)_{i\in I},(\al'_i)_{i\in I})$ by Theorem \ref{facialmorph}. By definition, 
\begin{eqnarray*}
    K'(\Th) &=&  \bigg\{ \sum_{i\in \Th}r_i h_i'  \,\big| \, { r_i\in\R^+ \mb{ for all }i\in I  \;,\; \sum_{i\in \Th}r_i a_{ij}\leq 0 \mb{ for all }j\in I } \bigg\} ,\\
    \overline{K'(\Th)} &=&  \bigg\{ \sum_{i\in \Th}r_i h_i'  \,\bigg|\, { r_i\in\R^+_0 \mb{ for all }i\in I \;,\;\sum_{i\in \Th}r_i a_{ij}\leq 0 \mb{ for all }j\in I } \bigg\},
\end{eqnarray*}
and $K(\Th)$, $\overline{K(\Th)}$ are described similarly. Trivially, $\phi(K'(\Th))=K(\Th)$, $\phi( \overline{K'(\Th)})= \overline{K(\Th)}$. If we identify the Coxeter groups of the root bases by Proposition \ref{realCM2}, then $\phi$ is $\We$-equivariant. Hence we get $\phi (F'(\Th))= F(\Th)$.
\qed

As before, $(\h,(h_i)_{i\in I},(\al_i)_{i\in I})$ is a root base of the linear Coxeter system $(\h,\,(H_i , \,L_i)_{i\in I})$. 
\begin{thm}\label{incl K in ri F} Let $\Th\subseteq I$ be special facial. Then we have: 
\begin{itemize}
\item[(a)] $F(\Th) - F(\Th)=  \overline{K(\Th)} - \overline{ K(\Th) } = S_{\Th^{aff}}^\vee + \mb{span}\Mklz{ h_i }{i\in\Th^{ind}}$.
\item[(b)] $Z_{\mathcal W}(F(\Th))  = Z_{\mathcal W}(\overline{K(\Th)})  = \We_{\Th^{aff}\cup\Th^\bot} $.
\item[(c)] $K(\Th)\subseteq \mb{\rm ri}(F(\Th)) $.
\end{itemize}
\end{thm}
\Proof We describe $(\h ,(h_i)_{i\in I},(\al_i)_{i\in I})$ as a subquotient of a free root base $(\h'', (h_i'')_{i\in I},(\al_i'')_{i\in I})$ as in Corollary \ref{subquot}. With the notation of this corollary we have the following morphisms of root bases: 
\begin{equation*}
\xymatrix{ \h''\ar@{->>}[rd]_<<<<<<\phi
&
&\h\ar@{^{(}->}[ld] ^<<<<<<<\psi\\
&\h'}
\end{equation*}
Since the map $\psi$ is injective, the facial sets for $(\h ,(h_i)_{i\in I},(\al_i)_{i\in I})$ and $(\h' ,(h'_i)_{i\in I},(\al'_i)_{i\in I})$ coincide by Theorem \ref{facialmorph}. We conclude from Proposition \ref{morphZ1} that
\begin{equation}\label{mor K Kq F}
   \phi(K''(\Th))= \psi( K(\Th) ) \quad \mb{ and }\quad
   \phi(\overline{ K''(\Th)}) = \psi( \overline{ K(\Th) }) \quad\mb{ and }\quad
   \phi(F''(\Th)) = \psi( F(\Th) ) .
\end{equation}

Now (a) follows from the description of $ F''(\Th) - F''(\Th)=  \overline{K''(\Th)} - \overline{ K''(\Th) } $ in Corollary \ref{hullfree}, (\ref{mor K Kq F}), and the injectivity of $\psi$. 
Furthermore, we conclude from (a) that $Z_{\mathcal W}(F(\Th))=Z_{\mathcal W}(\overline{ K(\Th) } )$. We identify the Coxeter groups of the root bases by Proposition \ref{realCM2}. If we show
\begin{equation*}
Z_{\mathcal W}(\overline{ K(\Th) }) = Z_{\mathcal W}(\overline{ K''(\Th) }) ,
\end{equation*} 
then (b) follows from the description of $Z_{\mathcal W}(\overline{ K''(\Th) })$ in Corollary \ref{zenfree}.
Since $ \psi( \overline{ K(\Th) }) = \phi(\overline{ K''(\Th) }) $ by (\ref{mor K Kq F}), it is sufficient to show that for $z\in \overline{ K(\Th) }$, $z''\in \overline{ K''(\Th) }$ such that $\psi(z)=\phi(z'')$ we have 
\begin{equation*}
\mb{stab}_{\mathcal W}(z)=\mb{stab}_{\mathcal W}(z'').
\end{equation*}
Let $\sigma\in\We$. The morphisms $\phi$ and $\psi$ are $\We$-equivariant. Since $\psi$ is injective, $\sigma z=z$ is equivalent to $\sigma\phi(z'')=\phi(z'')$. By the description of $ \mb{ker}(\phi) $ in Corollary \ref{subquot} (a) and by Proposition \ref{fcsdiff} this is equivalent to 
\begin{equation*}
\sigma z'' - z''  \in   \{\,\sum_{i=1}^n r_i h_i'' \mid r\in L_{\bf h} \,\}\cap \sum_{i=1}^n \R^+_0 h_i'' =\{0\},
\end{equation*}
where we used that $h_1''$, \ldots, $h_n''$ are linearly independent, and $ L_{\bf h}\cap (\R^+_0)^n=\{0\}$.

It remains to show (c). From Corollary \ref{prikZf} (b) we obtain
\begin{equation*}
   K''(\Th)+\mb{ker}(\phi)\subseteq \mb{ri}( F''(\Th) )+\mb{ker}(\phi) = \mb{ri}(F''(\Th)+\mb{ker}(\phi)).
\end{equation*}
From (\ref{mor K Kq F}), Lemma \ref{conint}, and the injectivity of $\psi$ we find
\begin{equation*}
  \underbrace{ \psi^{-1}(\phi(K''(\Th)+\mb{ker}(\phi))) }_{= K(\Th)} \subseteq \psi^{-1}(\phi(  \mb{ri}(F''(\Th)+\mb{ker}(\phi))))
    =\mb{ri}(\underbrace{ \psi^{-1}(\phi(F''(\Th)+\mb{ker}(\phi)))}_{=F(\Th)}) .
\end{equation*}
\qed

The next theorem can be proved similary as Theorem \ref{FXf} if (\ref{part K in rel int}) and Theorem \ref{incl K in ri F} (a) is used. For symmetric linear Coxeter systems it is obtained as a combination of Theorem 10.3 and 10.7 of \cite{D}.
\begin{thm} \label{FZf}  (a) Let $\Th\subseteq I$ be special facial. Then
\begin{equation*}
       \mb{\rm ri}(F(\Th)) =\We_\Th\, K(\Th) .
\end{equation*}

(b) Every $\We$-orbit of $\Fa{Z}$ contains exactly one of the faces $F(\Th)$, $\Th$ special facial. A face $F\in\Fa{Z}$ is of the form $F(\Th)$, $\Th$ special facial, if and only if $\mb{\rm ri}(F)\,\cap \,K=\mb{\rm ri}(F)\cap (-\overline{C}^\vee)\neq \emptyset$. 
\end{thm}

Comparing Proposition \ref{Fspecialinc} and Proposition \ref{Rspecialinc} we obtain the following corollary. For symmetric linear Coxeter systems it is given as a part of Theorem 4.3 of \cite{D}.
\begin{cor} We obtain a $\We$-equivariant anti-isomorphism of partially ordered sets from $\Fa{X}$ to $\Fa{Z}$ by mapping $\sigma R(\Th)$ to $\sigma F(\Th)$, $\sigma\in\We$, $\Th$ special facial. In particular, the face lattices $\Fa{X}$ and $\Fa{Z}$ are anti-isomorphic. 
\end{cor}

From Theorem \ref{FXlmj} we get the following corollary. For symmetric linear Coxeter systems it can be obtained by combining Proposition 10.2, Remark 10.2, Proposition 2.17, and Theorem 10.3 of \cite{D}.
\begin{cor} Let $\tau\in\mb{}^{\Th_1\cup\Th_1^\bot}\We^{\Th_2\cup\Th_2^\bot}$, and $\Th_1,\Th_2\subseteq I$ be special facial. Then $\Th_1\cup\Th_2\cup \text{\rm red}(\tau)$ is a special set and
\begin{eqnarray}
   F(\Th_1)\cap\tau F(\Th_2)    &=&  F((\Th_1\cap \tau\Th_2)^\infty), \label{intim} \\
   F(\Th_1)\vee \tau F(\Th_2)   &=&    F (  ( \Th_1\cup\Th_2\cup \text{\rm red}(\tau) )^{fc}  ) .\label{joinim}
\end{eqnarray}
\end{cor}
\begin{rem} The lattice meet and lattice join of two arbitrary faces of $Z$ can be reduced to the formulas (\ref{intim}), (\ref{joinim}) by Theorem \ref{FZf} (b) and the stabilizer formula (\ref{sub2}) of Proposition \ref{Fspecialinc}. 
\end{rem}

We denote by the upper script ``$d$" the dual cone. From the definitions of $X$, $Z$, from Corollary \ref{specialfaceX}, Corollary \ref{stabRTh}, and Theorem \ref{FZf} (a) we get $X\subseteq Z^d$, and if $\Th\subseteq I$ is special facial then
\begin{equation*}
  R(\Th) = X\cap F(\Th)^\bot = X\cap (\R h)^\bot \quad\mb{ for all }\quad h\in\mb{ri}(F(\Th)).
\end{equation*}
Similarly, from the definitions of $X$, $Z$, from Theorem \ref{FXf} (a), Theorem \ref{specialfaceZ}, and Proposition \ref{Fspecialinc} we get $Z\subseteq X^d$, and if $\Th\subseteq I$ is special facial then
\begin{equation*}
  F(\Th) = Z\cap R(\Th)^\bot = Z\cap (\R \la)^\bot \quad\mb{ for all }\quad \la\in\mb{ri}(R(\Th)).
\end{equation*}
In particular, $X$ and $Z$ are a pair of semidual cones as defined in A.6 of \cite{D}.

%


%
\section{\label{SecInvsub} Invariant convex subcones of the Tits cone}
%
%
%
%
%
%
In this section we investigate the faces and the face lattices of Coxeter group invariant convex subcones of the Tits cone of a linear Coxeter system. As explained in the introduction, certain Weyl group invariant polyhedral cones underlie the theory of normal reductive linear algebraic monoids. The cross section lattice and the type maps introduced by M. S. Putcha and L. E. Renner in the theory of reductive linear algebraic monoids
can be used for the description of the faces and the face lattice of such a cone; compare with the books \cite{PB}, \cite{R3}, in particular, with Sections 7.2 and 7.5.1 of \cite{R3}. We use similar data for the description of the faces and the face lattices of Coxeter group invariant convex subcones of the Tits cone.

In this section $(\h,\,(H_i, \,L_i)_{i\in I})$ is a linear Coxeter system, and $(\h,(h_i)_{i\in I},(\al_i)_{i\in I})$ is a root base of the linear Coxeter system with generalized Cartan matrix $A$. For the whole section, except for Lemmas \ref{sfix} and \ref{sfix2}, $Y$ is a $\We$-invariant convex subcone of the Tits cone $X\subseteq \h^*$ which contains the zero. 

For $h\in\h$ we set $H_h^{=0}:=\Mklz{\la\in\h^*}{\la(h)= 0}$ and $H_h^{\geq 0}$, $H_h^{>0}$, $H_h^{\leq 0}$, $H_h^{<0}$ are defined analogously.
%
%
%
\subsection{A single face\label{a single face}}
%
%
%
For every face of a convex $\We$-invariant cone with zero in $\h^*$ the coroot system can be divided naturally into four parts according to the following lemmas:
\begin{lem} \label{sfix} Let $Y$ be a convex $\We$-invariant cone in $\h^*$ which contains the zero. Let $R\in \Fa{Y}$. For $h\in\Phi^\vee$ we have:
\begin{itemize}
\item[(a)] The reflection $\sigma_h$ fixes $R$ pointwise if and only if $R-R\subseteq H_h^{=0}$. 
\item[(b)] The reflection $\sigma_h$ fixes $R$ as a whole, but not pointwise if and only if $\al_h\in R-R$.
\end{itemize}
\end{lem}
\Proof (a) is trivial. To (b): Let $\al_h\in R-R$. Then $\sigma_h R\subseteq Y$ and $\sigma_h R\subseteq R-R$, from which we get $\sigma_h R\subseteq Y\cap (R-R)=R$. Furthermore, $R\subseteq   (\sigma_h)^{-1}R=\sigma_h R$. Since $\al_h(h)=2$, we conclude from (a) that $\sigma_h$ does not fix $R$ pointwise. 

Now suppose that $\sigma_h$ fixes $R$ as a whole, but not pointwise. Then it follows from (a) that there exists an element $\la_0\in R$ such that $\la_0(h)\neq 0$. We obtain
\begin{equation*}
  \underbrace{\sigma_h \la_0}_{\in R} =\underbrace{\la_0}_{\in R}-\underbrace{\la_0(h)}_{\neq 0} \al_h,
\end{equation*}
which shows that $\al_h\in R-R$.
\qed

\begin{lem}\label{sfix2}  Let $Y$ be a convex $\We$-invariant cone in $\h^*$ which contains the zero. Let $R\in \Fa{Y}$. Let $h\in\Phi^\vee$ such that the reflection $\sigma_h$ does not fix $R$ as a whole. Then exactly one of the following possibilities holds:
\begin{itemize}
\item[($p$)] $R\subseteq H_h^{\geq 0}$ with $\mb{ri}(R)\subseteq H_h^{>0}$.   
\item[($n$)] $R\subseteq H_h^{\leq 0}$ with $\mb{ri}(R)\subseteq H_h^{<0}$. 
\end{itemize}
\end{lem}
 \Proof Suppose that there exists $\la_0\in \mb{ri}(R)$ such that $\la_0(h)=0$. Then $\sigma_h \la_0=\la_0$. We conclude that $\mb{ri}(\sigma_h R)\cap \mb{ri(R)}\neq\emptyset$, which implies $\sigma_h R=R$. However, this has been excluded. 

Suppose that there exist $\la_1, \la_2\in\mb{ri}(R)$ such that $\la_1(h)>0$ and $\la_2(h)<0$. The closed line segment $\overline{\la_1\la_2}$ is contained in $\mb{ri}(R)$. Since the evaluation in $h$ is continuous, there exists $\la_0\in \overline{\la_1\la_2}\setminus\{\la_1,\,\la_2\}$ such that $\la_0(h)=0$, which is not possible.

We may assume $\la(h)> 0$ for all $\la\in \mb{ri}(R)$. Let $\mu\in R$. Choose an element $\la\in\mb{ri}(R)$. Then  $\overline{\la\mu}\setminus\{\mu\}$ is contained in $\mb{ri}(R)$. Since the evaluation in $h$ is continuous, we find $\mu(h)\geq 0$.
\qed
\mb{}\\

In the remaining section $Y$ is a $\We$-invariant convex subcone of the Tits cone $X$ which contains the zero. We call
\begin{equation*}
  \Upsilon:=\Mklz{R\in\Fa{Y}}{\mb{ri}(R)\cap\overline{C}\neq\emptyset }
\end{equation*}
the cross section lattice of $\Fa{Y}$ obtained by the chamber $\overline{C}$. This name will be justified by Corollary \ref{csp} and Theorem \ref{csublattice} below. 
For $R\in\Upsilon$ we set
\begin{equation*}
  \upsilon(R):= \Mklz{i\in I}{\sigma_i R=R }.
\end{equation*}
By Lemma \ref{sfix} we may define
\begin{eqnarray*}
  &&\upsilon_*(R):=  \Mklz{i\in I}{  \sigma_i \mb{ fixes }R\mb{ pointwise }  } = 
\Mklz{i\in I}{   R-R\subseteq H^{=0}_i  },\\
  &&\upsilon^*(R):= \upsilon(R)\setminus\upsilon_*(R) = \Mklz{i\in I}{ L_i\subseteq R-R  }.
\end{eqnarray*}
We have the disjoint union
\begin{equation*}
   \upsilon(R) = \upsilon_*(R) \,\dot{\cup} \,\upsilon^*(R).
\end{equation*}
Clearly, $\upsilon_*(R)$ and $\upsilon^*(R)$ are separated. In particular,
\begin{equation*}
  \We_{\upsilon(R)}=\We_{\upsilon_*(R)}\We_{\upsilon^*(R)}=\We_{\upsilon^*(R)}\We_{\upsilon_*(R)}.
\end{equation*}
We denote by ${\mathcal P(I)}$ the power set of $I$. The map $\upsilon:\Upsilon\to {\mathcal P}(I)$ is called the type map. The maps $\upsilon_*,\,\upsilon^*:\Upsilon\to {\mathcal P}(I)$ are called the lower and upper type maps.

\begin{thm}\label{maincsx} Let $R\in\Upsilon$.
\begin{itemize}
\item[(a)] The pointwise stabilizer of $R$ and the stabilizer of $R$ as a whole are given by
\begin{equation*}
  Z_{\mathcal W}(R)= \We_{\upsilon_*(R)} \quad \mb{ and }\quad N_{\mathcal W}(R) =\We_{\upsilon(R)}  .
\end{equation*}
\item[(b)] We have $R\cap\overline{C} = R\cap\overline{F}_{\upsilon_*(R)}$, and
\begin{equation*}
  R=N_{\mathcal W}(R) (R\cap\overline{C}) = \We_{\upsilon^*(R)}(R\cap\overline{F}_{\upsilon_*(R)})
\end{equation*}
where $\mb{ri}(R)\cap F_{\upsilon_*(R)}\neq\emptyset $. In particular, the set $\upsilon_*(R)$ is facial.
\end{itemize}
\end{thm}
\Proof The proof of the theorem is divided into the parts (i) - (vi). Since $R\in \Upsilon$ there exists an element $\la_0\in \mb{ri}(R)\cap\overline{C}$, which we fix for the whole proof.

(i) Let $\tau\in\We$ and $\la_1\in R\cap\tau\overline{C}$. We first show that
\begin{equation*}
\la_1\in\We_{\upsilon(R)}(R\cap \overline{C})
\end{equation*}
by induction over the length $l(\tau)$ of $\tau\in\We$.
This is trivial for $l(\tau)=0$, i.e., $\tau=1$. Let $l(\tau)=m>1$. Write $\tau$ in the form $\tau=\sigma_i \sigma$ with $i\in I$ and $\sigma\in\We$ such that $l(\sigma)=l(\tau)-1$. 

If $i\in\upsilon(R)$ then $\sigma_i$ fixes $R$ as a whole, and we get
\begin{equation*}
   \sigma_i \la_1  \in  \sigma_i (R\cap\tau\overline{C})=\sigma_i R\cap\sigma_i\tau\overline{C} = R\cap\sigma\overline{C} .
\end{equation*}
By induction, $\sigma_i\la_1\in \We_{\upsilon(R)}(R\cap \overline{C})$. Hence $\la_1\in \We_{\upsilon(R)}(R\cap \overline{C})$.

If $i\notin\upsilon(R)$ then $\sigma_i$ does not fix $R$ as a whole. Since $\la_0\in\mb{ri}(R)\cap \overline{C}$, we have $\la_0(h_i)\geq 0$. Hence we obtain from Lemma \ref{sfix2} that $\la_0(h_i)>0$ and $\la_1(h_i)\geq 0$. It is $l(\tau^{-1}\sigma_i)<l(\tau^{-1})$. By Theorem \ref{ThHN}, this implies $\tau^{-1}h_i\in -\sum_{j\in I}\R^+_0 h_j$, from which we find
\begin{equation*}
  \la_1(h_i)= \underbrace{(\tau^{-1}\la_1)}_{\in\overline{C}}(\tau^{-1}h_i)\leq 0.
\end{equation*}
We conclude that $\la_1(h_i)=0$ and
\begin{equation*}
    R\ni \la_1= \sigma_i\la_1 \in\sigma_i\tau  \overline{C} =\sigma \overline{C}.
\end{equation*}
By induction we get $\la_1\in \We_{\upsilon(R)}(R\cap \overline{C})$.

(ii)  From the definition of $\upsilon(R)$ we obtain $\We_{\upsilon(R)}\subseteq N_{\mathcal W}(R)$. By (i) we find
\begin{equation*}
  R\subseteq \We_{\upsilon(R)}(R\cap \overline{C}) \subseteq  \We_{\upsilon(R)} R\subseteq R.
\end{equation*}
Hence $ R = \We_{\upsilon(R)}(R\cap \overline{C}) $.

(iii) Now we use (ii) to show that $N_{\mathcal W}(R) \subseteq  \We_{\upsilon(R)}$.
Let $\sigma\in N_{\mathcal W}(R) $. Then we have
\begin{equation*}
   \sigma\la_0\in R=\We_{\upsilon(R)}(R\cap \overline{C}) .
\end{equation*}
Since every $\We$-orbit in $X$ intersects $\overline{C}$ in exactly one point, there exists $\tau\in\We_{\upsilon(R)}$ such that $\sigma\la_0=\tau \la_0$. It follows that
\begin{equation*}
     \sigma \in \tau \We_J \quad\mb{ with }\quad J:=\Mklz{j \in I}{\sigma_j \la_0=\la_0}.
\end{equation*}
If $j\in J$ then $\mb{ri}(\sigma_j R)\cap \mb{ri}(R)\neq\emptyset $, which implies $\sigma_j R=R$, so $j\in\upsilon(R)$. Thus we obtain $\sigma\in \We_{\upsilon(R)}$.

(iv) The set 
\begin{equation*}
    J_{min}:=\bigcap_{J\text{ facial},\;R\cap F_J\neq\emptyset} J   
\end{equation*}
is facial by Proposition \ref{fsFa0} (a), and 
\begin{equation*}
 R\cap\overline{C}=\bigcup_{J \text{ facial},\;R\cap F_J\neq\emptyset} R\cap F_J \subseteq R\cap\overline{F}_{J_{min}} \subseteq  R\cap\overline{C}
\end{equation*}
shows that $R\cap \overline{C} =  R\cap\overline{F}_{J_{min}}$. For every facial set $J$ with $R\cap F_J \neq\emptyset$ we choose an element $\la_J\in R\cap F_J$. Then
\begin{equation*}
   \mu:=\sum_{J\text{ facial},\;R\cap F_J\neq\emptyset} \la_J\;\in \; R\cap F_{J_{min}} . 
\end{equation*}
Furthermore, $\mu+\la_0\in F_{J_{min}}+\overline{F}_{J_{min}}\subseteq F_{J_{min}}$ and $\mu+\la_0\in\mb{ri}(R)$.

(v) We now show that $Z_{\mathcal W}(R)=Z_{\mathcal W}(R\cap\overline{C})$. The inclusion ``$\subseteq$" holds trivially; we only have to show the inclusion ``$\supseteq $".
We find from (ii) that
\begin{equation*}
  Z_{\mathcal W}(R) =\bigcap_{\sigma\in {\mathcal W}_{\upsilon(R)}}  \sigma Z_{\mathcal W}(R\cap\overline{C})\sigma^{-1}.
\end{equation*}
Every $\sigma\in\We_{\upsilon(R)}\setminus\{1\}$ can be written in the form $\sigma=\sigma_{i_1}\cdots\sigma_{i_{m-1}}\sigma_{i_m}$ with $i_1,\,\cdots,\,i_{m-1},\,i_m\in\upsilon(R)$ and we have
\begin{equation*}
   \sigma Z_{\mathcal W}(R\cap\overline{C})\sigma^{-1}=  \sigma_1\left(\cdots\left( \sigma_{m-1}\left(\sigma_m  Z_{\mathcal W}(R\cap\overline{C})\sigma_m^{-1}\right)\sigma_{m-1}^{-1}\right)\cdots\right)\sigma_1^{-1}.
\end{equation*}
Therefore, it is sufficient to show that
\begin{equation*}
      \sigma_i Z_{\mathcal W}(R\cap\overline{C}) \sigma_i^{-1}= Z_{\mathcal W}(R\cap\sigma_i\overline{C}) \supseteq  Z_{\mathcal W}(R\cap\overline{C})
\end{equation*}
for all $i\in\upsilon(R)$.

Let $\la\in R\cap\sigma_i \overline{C}$ where $i\in\upsilon(R)$. The closed line segment $\overline{\la\, \sigma_i\la}$, which may also consist of a single point, is contained in $R$ with endpoints $\la\in R\cap\sigma_i \overline{C}$, $\sigma_i\la\in R\cap\overline{C}$, and midpoint
\begin{equation*}
    \mu:= \frac{1}{2}(\la+\sigma_i\la)\in R\cap (\overline{C}\cap\sigma_i\overline{C}).
\end{equation*}
Because $\sigma_i\la$ and $\mu$ are fixed by the elements of $Z_{\mathcal W}(R\cap\overline{C})$ also $\la= 2\mu-\sigma_i\la$ is fixed.

(vi) If we combine (iv) and (v) we find 
\begin{equation*}
     Z_{\mathcal W}(R) = Z_{\mathcal W}(R\cap\overline{C})  = \We_{J_{min}},
\end{equation*}
from which we get $  \upsilon_*(R)=\Mklz{i\in I }{\sigma_i \in  Z_{\mathcal W}(R) } =J_{min}$.
\qed

\begin{cor}\label{mainlh} The linear hull of $R\in\Upsilon$ is given by
\begin{equation*}
    R-R = \left( R\cap\overline{C} \right) -  \left( R\cap\overline{C} \right) =  \left(R\cap\overline{F}_{\upsilon_*(R)} \right) - \left( R\cap\overline{F}_{\upsilon_*(R)} \right) .
\end{equation*}
\end{cor}
\Proof The inclusions ``$\supseteq$" hold trivially. It is sufficient to show that
\begin{equation*}
   R   \subseteq \left(R\cap\overline{F}_{\upsilon_*(R)} \right) - \left( R\cap\overline{F}_{\upsilon_*(R)} \right) .
\end{equation*}
By Theorem \ref{maincsx} (b) there exists $\la_0\in R \cap  F_{\upsilon_*(R)}$. For $\la\in R$ we get
\begin{equation*}
   \la=(\la+ t\la_0)- t\la_0 \in \left(R\cap F_{\upsilon_*(R)} \right) - \left( R\cap F_{\upsilon_*(R)} \right) \subseteq \left(R\cap\overline{F}_{\upsilon_*(R)} \right) - \left( R\cap\overline{F}_{\upsilon_*(R)} \right) 
\end{equation*}
for all sufficiently big $t\in\R^+$, because of $R\subseteq \{\la\in \h^*\mid \la(h_i)=0 \mb{ for all }i\in \upsilon_*(R) \}$. 
\qed

Since $Y$ is partitioned by the relative interiors of the faces of $Y$, the definition of $\Upsilon$ implies 
\begin{equation}\label{rifd1}
     Y \cap \overline{C} = \dot{\bigcup_{R\in\Upsilon}} \;\mb{ri}(R)\cap\overline{C}  .
\end{equation}
We find from Theorem \ref{maincsx} that the relative interior of $R\in\Upsilon$ is given by
\begin{equation}\label{rifd2}
   \mb{ri}(R)=\We_{\upsilon^*(R)}\, ( \mb{ri}(R)\cap\overline{C} ) .
\end{equation}
The intersections $\mb{ri}(R)\cap\overline{C}$, $R\in\Upsilon$, which appear in (\ref{rifd1}), (\ref{rifd2}), can be described as follows:
\begin{thm}\label{maincsx2} Let $R\in\Upsilon$. Then
\begin{equation}\label{riallg1}
  \mb{\rm ri}(R) \cap\overline{C} = \dot{\bigcup_{J_f\subseteq \upsilon^*(R),\, (J_f)^0=J_f}}\underbrace{\mb{\rm ri}(R)\cap F_{\upsilon_*(R)\cup J_f}}_{\neq\emptyset}
\end{equation}
and for $ J_f\subseteq \upsilon^*(R),\, (J_f)^0=J_f $ we have
\begin{equation}\label{riallg2}
    \mb{mid}_{J_f}(\mb{\rm ri}(R\cap\overline{C}))=\mb{mid}_{J_f}(\mb{\rm ri}(R\cap\overline{F}_{\upsilon_*(R)})) = \mb{\rm ri}(R)\cap F_{\upsilon_*(R)\cup J_f} .
\end{equation}
\end{thm}
\Proof  We know from Theorem \ref{maincsx} (b) that
\begin{equation*}
    \mb{ri}(R)\cap\overline{C} = \bigcup_{\upsilon_*(R)\subseteq L} \mb{ri}(R)\cap F_L
\end{equation*}
and $\mb{ri}(R)\cap F_{\upsilon_*(R)}\neq\emptyset$.
 
(a) Let $\upsilon_*(R)\subseteq L $ and $\mb{ri}(R)\cap F_L\neq \emptyset$. We first show that $L\subseteq \upsilon(R)$. Let $\la\in \mb{ri}(R)\cap F_L$. For every $l\in L$ we have $\sigma_l\la=\la$. We conclude that $ \mb{ri}(\sigma_l R)\cap \mb{ri}(R)\ne\emptyset$, which implies $\sigma_l R= R$. Thus $l\in\upsilon(R)$.

(b) Let $J\subseteq \upsilon^*(R)$ and $\mb{ri}(R)\cap F_{\upsilon_*(R)\cup J}\neq \emptyset$. We now show that $J^0 =J$. The set $\upsilon_*(R)\cup J$ is facial. Since $\upsilon_*(R)$ and $J\subseteq \upsilon^*(R)$ are separated, we find $(\upsilon_*(R)\cup J)^\infty= \upsilon_*(R)^\infty\cup J^\infty$, which is also a facial set by Theorem \ref{fsp}. Recall from Corollary \ref{specialfaceX} that we get a face of the Tits cone $X$ by 
\begin{equation*}
    R(\upsilon_*(R)^\infty\cup J^\infty) =\We_{(\upsilon_*(R)^\infty\cup J^\infty)^\bot} \overline{F}_{\upsilon_*(R)^\infty\cup J^\infty}.
\end{equation*}
Hence $R\cap R(\upsilon_*(R)^\infty\cup J^\infty)$ is a face of $R$, and 
\begin{equation*}
  \mb{ri}(R)\cap ( R\cap R(\upsilon_*(R)^\infty\cup J^\infty)) \supseteq \mb{ri}(R)\cap F_{\upsilon_*(R)\cup J}\neq\emptyset.
\end{equation*}
It follows that 
\begin{equation*}
   R\subseteq  R\cap R(\upsilon_*(R)^\infty\cup J^\infty)\subseteq  R(\upsilon_*(R)^\infty\cup J^\infty),
\end{equation*}
from which we obtain
\begin{equation*}
    R\cap F_{\upsilon_*(R)}\subseteq R\cap\overline{C}\subseteq R(\upsilon_*(R)^\infty\cup J^\infty)\cap\overline{C}= \overline{F}_{\upsilon_*(R)^\infty\cup J^\infty}.
\end{equation*}
Because of $R\cap F_{\upsilon_*(R)}\neq\emptyset$ we get $\upsilon_*(R)\supseteq \upsilon_*(R)^\infty\cup J^\infty\supseteq J^\infty$. Moreover, we have $J^\infty\subseteq J\subseteq \upsilon^*(R)$. Since $\upsilon_*(R)$ and $ \upsilon^*(R)$ are disjoint, we find $J^\infty=\emptyset$.

(c) We get $R\cap\overline{C}= R\cap\overline{F}_{\upsilon_*(R)}$ and $\mb{ri}(R)\cap \mb{ri}(\overline{F}_{\upsilon_*(R)}) = \mb{ri}(R)\cap F_{\upsilon_*(R)}\neq\emptyset$ from Theorem \ref{maincsx} (b) and Proposition \ref{rihfc}. We conclude from Theorem 6.5 of \cite{Ro} that
\begin{equation*}
  \mb{ri}(R\cap\overline{C} )= \mb{ri}( R\cap\overline{F}_{\upsilon_*(R)}) =  \mb{ri}(R)\cap F_{\upsilon_*(R)} .
\end{equation*} 

To show the inclusion ``$\subseteq$" in (\ref{riallg2}) let $\la\in \mb{ri}(R)\cap F_{\upsilon_*(R)}$. Since $\We_{J_f}\subseteq \We_{\upsilon(R)}= N_{\mathcal W}(R)$, we find
\begin{equation*}
   \sigma\la\in \sigma\,\mb{ri}(R)=\mb{ri}(\sigma R)=\mb{ri}(R) \quad \mb{ for all }\quad\sigma\in\We_{J_f},
\end{equation*}
which implies $\mb{mid}_{J_f}(\la)\in \mb{ri}(R)$. Since $\upsilon_*(R)$ is facial, and $J_f\subseteq \upsilon^*(R)\subseteq \upsilon_*(R)^\bot$, we obtain
\begin{equation*}
    \mb{mid}_{J_f}(\la) \in \mb{ri}(R)\cap F_{\upsilon_*(R)\cup J_f}
\end{equation*}
from the first part of Lemma \ref{midpoint}. In particular, $ \mb{ri}(R)\cap F_{\upsilon_*(R)\cup J_f}\neq\emptyset$.

We use the second part of Lemma \ref{midpoint} to show the inclusion ``$\supseteq$" in (\ref{riallg2}). Let $\la\in\mb{ri}(R)\cap F_{\upsilon_*(R)\cup J_f}$. Let $\gamma\in\sum_{i\in J_f}\R^+\al_i$ such that $\gamma(h_i)>0$ for all $i\in J_f$. Then
\begin{equation*}
  \la+t\gamma\in F_{\upsilon_*(R)} \quad\mb{ and }\quad  \mb{mid}_{J_f}(\la+t\gamma) = \la 
\end{equation*}
for all sufficiently small $t\in\R^+$. From $J_f\subseteq \upsilon^*(R)=\Mklz{i\in I}{\al_i\in R-R}$ we obtain $\gamma\in R-R$. Since $\la\in\mb{ri}(R)$, and $\mb{ri}(R)$ is open in $R-R$, we find $\la+t\gamma\in \mb{ri}(R)$ for all sufficiently small $t\in\R^+$. 
\qed
%
%
%
\subsection{Inclusion of faces\label{inclusion of faces}}
%
%
%
The next aim is to describe the inclusion of faces, which are $\We$-translates of faces from $\Upsilon$. To do this we first describe the inclusion of faces from $\Upsilon$.
\begin{prop}\label{inclspeziell} Let $R_1, R_2\in\Upsilon$. The following are equivalent:
\begin{itemize}
\item[(i)] $R_1\subseteq R_2$.
\item[(ii)] $R_1\cap\overline{C} \subseteq R_2\cap\overline{C}$.
\item[(iii)] $R_1\cap\overline{F}_{\upsilon_*(R_1)}\subseteq R_2\cap\overline{F}_{\upsilon_*(R_2)}$.
\end{itemize}
In this case we have $\upsilon_*(R_1)\supseteq\upsilon_*(R_2)$ and $\upsilon^*(R_1)\subseteq\upsilon^*(R_2)$. 
\end{prop}
\begin{rem} In general the pair of functions $\upsilon_*$, $\upsilon^*$ does not separate the elements of $\Upsilon$. For some extreme examples see (1) and (2) of Lemma \ref{vextrval}.
\end{rem}
\Proof (i) implies (ii), which in turn is equivalent to (iii) by Theorem \ref{maincsx} (b). From (iii) and Corollary \ref{mainlh} we find $R_1-R_1\subseteq R_2-R_2$. If we intersect with $Y$ on both sides we obtain (i). The properties of $\upsilon_*$ and $\upsilon^*$ follow immediately from 
the definitions of $\upsilon_*$ and $\upsilon^*$.
\qed
 
Now we can show:
\begin{thm}\label{inclallg} Let $R_1, R_2\in\Upsilon$ and $\sigma_1, \sigma_2\in\We$. The following are equivalent:
\begin{itemize}
\item[(i)]  $\sigma_1 R_1\subseteq \sigma_2 R_2$.
\item[(ii)] $R_1\subseteq R_2$ and $\sigma_1 N_{\mathcal W}(R_1) \cap \sigma_2 N_{\mathcal W}(R_2)\neq \emptyset$.
\item[(iii)] $R_1\subseteq R_2$ and $\sigma_1^{-1}\sigma_2\in N_{\mathcal W}(R_1) N_{\mathcal W}(R_2) = \We_{\upsilon_*(R_1)}\We_{\upsilon^*(R_2)}$.
\end{itemize}
\end{thm}
\Proof For $R_1\subseteq R_2$ we obtain from Theorem \ref{maincsx} (a) and Proposition \ref{inclspeziell} that
\begin{eqnarray*}
  && N_{\mathcal W}(R_1) N_{\mathcal W}(R_2) 
       = \We_{\upsilon_*(R_1)}\We_{\upsilon^*(R_1)}\We_{\upsilon^*(R_2)}\We_{\upsilon_*(R_2)} 
       =  \We_{\upsilon_*(R_1)}\We_{\upsilon^*(R_2)}\We_{\upsilon_*(R_2)}\\
  && = \We_{\upsilon_*(R_1)}\We_{\upsilon_*(R_2)}\We_{\upsilon^*(R_2)} = \We_{\upsilon_*(R_1)}\We_{\upsilon^*(R_2)}.
\end{eqnarray*}

The equivalence of (ii) and (iii) is trivial. We next show the implication ``(i) $\Rightarrow$ (iii)": By Theorem \ref{maincsx} (b) we have
\begin{equation*}
   \sigma_1\We_{\upsilon^*(R_1)}(R_1\cap\overline{F}_{\upsilon_*(R_1)})\subseteq  \sigma_2\We_{\upsilon^*(R_2)}(R_2\cap\overline{F}_{\upsilon_*(R_2)})
\end{equation*}
and $R_1\cap F_{\upsilon_*(R_1)}\neq\emptyset$. Let $\la\in R_1\cap F_{\upsilon_*(R_1)}$. Since every $\We$-orbit in $X$ intersects $\overline{C}$ in exactly one point, there exists $\tau\in\We_{\upsilon^*(R_2)}$ such that $\sigma_1\la=\sigma_2\tau\la$. We conclude that $\sigma_1^{-1}\sigma_2\tau\in \We_{\upsilon_*(R_1)}$. Hence
\begin{equation*}
  \sigma_1^{-1}\sigma_2\in \We_{\upsilon_*(R_1)}\tau^{-1}\subseteq  \We_{\upsilon_*(R_1)}\We_{\upsilon^*(R_2)}.
\end{equation*}

Decompose $ \sigma_1^{-1}\sigma_2=ab$ with $a\in\We_{\upsilon_*(R_1)}$ and $b\in\We_{\upsilon^*(R_2)}$. By Theorem \ref{maincsx} (a) the inclusion $\sigma_1R_1\subseteq \sigma_2 R_2$ is equivalent to
\begin{equation*}
  R_1\subseteq \sigma^{-1}\sigma_2R_2 =ab R_2=a R_2,
\end{equation*}
which is equivalent to $R_1=a^{-1}R_1\subseteq R_2$.
This chain of equivalences also shows the implication ``(iii) $\Rightarrow$ (i)".
\qed

From the preceding theorem we obtain immediately:
\begin{cor}\label{fincleq} Let $R\in\Fa{Y}$ and $\sigma\in \We$. Then we have $\sigma R=R$, if and only if $\sigma R\subseteq R$, if and only if $\sigma R\supseteq R$.
\end{cor}

%
\subsection{The cross section lattice\label{the cross section lattice}}
%
%
%
The following corollary of Theorem \ref{inclallg} is the reason why $\Upsilon$ has been called a cross section.
\begin{cor}\label{csp} Let $R'\in\Fa{Y}$. Then there exist $\sigma\in\We$ and a uniquely determined $R\in\Upsilon$ such that $R'=\sigma R$.
\end{cor}
\Proof Let $\la'\in\mb{ri}(R')$. There exists $\sigma\in\We$ such that $\sigma^{-1}\la'\in\overline{C}$. Since $\sigma^{-1}\la'\in\mb{ri}(\sigma^{-1}R')$, we have
\begin{equation*}
    R'=\sigma (\sigma^{-1}R')  \quad\mb{ with }\quad \sigma^{-1}R'\in\Upsilon .
\end{equation*}
Let $R'=\sigma_1 R_1=\sigma_2 R_2$ with $R_1,R_2\in\Upsilon$ and $\sigma_1,\sigma_2\in\We$. Then we get $R_1=R_2$ from Theorem \ref{inclallg} .
\qed

Corollary \ref{csp} can be generalized as follows:
\begin{thm}\label{chainconj} Let $R_1\subseteq R_2\subseteq\cdots\subseteq R_m$ be a chain of faces of $Y$. Then there exists uniquely a chain $S_1\subseteq S_2\subseteq\cdots\subseteq S_m$ of faces from $\Upsilon$, there exists an element $\sigma\in\We$ such that $R_i=\sigma S_i$, $i=1,2,\ldots, m$.
\end{thm}
\Proof We prove the theorem by induction over $m$. For $m=1$ the theorem coincides with Corollary \ref{csp}. Let  $m>1$ and $R_1\subseteq R_2\subseteq\cdots\subseteq R_m$ be a chain of faces of $Y$. By Corollary \ref{csp} there exist $S_1\in\Upsilon$ and $\tau\in\We$ such that $R_1=\tau S_1$. By induction there exists a chain $S_2\subseteq\cdots\subseteq S_m$ of faces from $\Upsilon$ and $\sigma\in\We$ such that $R_i=\sigma S_i$, $i=2,\ldots, m$. Thus, the chain of faces of $Y$ is of the form
\begin{equation*}
   \tau S_1\subseteq \sigma S_2\subseteq\cdots\subseteq \sigma S_m \quad\mb{ with } \quad S_2\subseteq\cdots\subseteq S_m .
\end{equation*}
From Theorem \ref{inclallg} and Proposition \ref{inclspeziell} we obtain
\begin{equation*}\label{chain ind step}
   S_1\subseteq S_2 \quad\mb{ and } \quad \tau^{-1}\sigma \in\We_{\upsilon_*(S_1)}\We_{\upsilon^*(S_2)} \quad\mb{ where } \quad \We_{\upsilon^*(S_2)}\subseteq \cdots\subseteq \We_{\upsilon^*(S_m)} .
\end{equation*}
In particular, we get $S_1\subseteq S_2\subseteq\cdots\subseteq S_m$. We write $\tau^{-1}\sigma$ in the form $\tau^{-1}\sigma=ab$ with $a\in \We_{\upsilon_*(S_1)}$ and $b\in\We_{\upsilon^*(S_2)}\subseteq\cdots\subseteq\We_{\upsilon^*(S_m)} $. By Theorem \ref{maincsx} (a) we find
\begin{equation*}
 R_1=\tau S_1 =\sigma b^{-1}a^{-1}S_1=\sigma b^{-1}S_1 \quad\mb{ and }\quad
 R_k=\sigma S_k= \sigma b^{-1}S_k\quad \mb{ for }\quad k=2, \ldots, m.
\end{equation*}
The uniqueness of $S_1, S_2, \ldots, S_m\in\Upsilon$ follows from Corollary \ref{csp}.
\qed

The following theorem is the reason why $\Upsilon$ has been called a lattice.
\begin{thm}\label{csublattice} $\Upsilon$ is a complete sublattice of $\Fa{Y}$.
\end{thm}
\Proof We first show that  $\Upsilon$ contains the biggest face of $Y$, and the smallest face of $Y$. By Corollary \ref{csp} there exists $\sigma\in\We$ such that $\sigma^{-1}Y\in\Upsilon$. It is $\sigma^{-1}Y=Y$. By Corollary \ref{csp} there exists $\tau\in\We$ such that $\tau^{-1}(Y\cap (-Y))\in\Upsilon$. It is $\tau^{-1}(Y\cap (-Y))=(\tau^{-1}Y)\cap (-\tau^{-1}Y)=Y\cap (-Y)$.

Now let $(R_l)_{l\in L}$ be a nonempty family of elements in $\Upsilon$. By Corollary \ref{csp} there exist $R\in \Upsilon$ and $\sigma\in\We$ such that 
$  \bigcap_{l\in L} R_l =\sigma R$. For every $l\in L$ we have $\sigma R\subseteq R_l$, from which we obtain $R\subseteq R_l$ by Theorem \ref{inclallg}. We conclude that
$ R\subseteq \bigcap_{l\in L} R_l =\sigma R$. From Corollary \ref{fincleq} we find $ \bigcap_{l\in L}R_l= R\in\Upsilon$.
Similarly, it follows that $\vee_{l\in L} R_l\in\Upsilon$.
\qed

\begin{prop} Let $R, S\in\Upsilon$. Then $ \upsilon(R)\cap\upsilon(S) \subseteq  \upsilon(R\vee S)\cap\upsilon(R\cap S)$.
\end{prop}
\Proof Let $i\in \upsilon(R)\cap\upsilon(S)$. Then $\sigma_i R=R$ and $\sigma_i S=S$ by definition. Hence $\sigma_i (R\cap S)=\sigma_i R\cap\sigma_i S=R\cap S$, which shows that $i\in\upsilon(R\cap S)$. Similarly, we find $i\in \upsilon (R\vee S)$.
\qed

\begin{prop} If the convex cone $Y\cap\overline{C}$ has only finitely many faces, then the cross section lattice $\Upsilon$ is finite. 
\end{prop}
\begin{rem}In particular, if $Y\cap\overline{C}$ is a finitely generated convex cone, then the cross section lattice $\Upsilon$ is finite.
\end{rem}
\Proof By Theorem \ref{maincsx} (b), $R\in\Upsilon$ is of the form $ R= \We_{\upsilon^*(R)}(R\cap\overline{C}) $. Moreover, $R\cap \overline{C}$ is a face of $Y\cap\overline{C}$, because $R$ is a face of $Y$. There exist only finitely many faces of $Y\cap\overline{C}$, and only finitely many standard parabolic subgroups of $\We$. Therefore, $\Upsilon$ is finite.
\qed

%
\subsection{The relation to some structures from the literature\label{relation to literature}}
%
%
Let $({\mathcal U},\leq)$ be a lower semilattice which has a biggest element $1\in {\mathcal U}$. A map $\la:{\mathcal U}\to {\mathcal P}(I)$ with
$\la(1)=I$ is called a regular type map if we obtain on the set
\begin{equation*}
 {\mathcal F}({\mathcal U},\la):=\Mklz{(u,\sigma \We_{\la(u)}) }{u\in{\mathcal U},\,\sigma\in\We }
\end{equation*}
the structure of a lower semilattice with biggest element $(1,\We)$ by
\begin{equation*}
    (u_1,\sigma_1\We_{\la(u_1)})\leq  (u_2,\sigma_2\We_{\la(u_2)}) \quad :\iff \quad u_1\leq u_2 \;\mb{ and }\;\sigma_1\We_{\la(u_1)}\cap \sigma_2\We_{\la(u_2)}\neq\emptyset .
\end{equation*}
Regular type maps are important to describe the structure of regular monoids on groups with BN-pairs in \cite{P2}. A characterization of regular type maps has been obtained in \cite{AP}. 
\begin{cor} The map $\upsilon:\Upsilon\to {\mathcal P}(I)$ is a regular type map, such that $ {\mathcal F}(\Upsilon,\upsilon)$ is isomorphic to the face lattice $\Fa{Y}$.
\end{cor}
\Proof $\Upsilon$ is a complete lattice by Theorem \ref{csublattice}, and $Y$ is its biggest element. Trivially, we have $\upsilon(Y)=I$. From Theorem \ref{maincsx} (a), Theorem \ref{inclallg}, and Corollary \ref{csp} it follows that 
\begin{eqnarray*}
   {\mathcal F}(\Upsilon,\upsilon)\quad &\to & \Fa{Y} \\
     (R,\sigma \We_{\upsilon(R)}) &\mapsto &\;\; \sigma R
\end{eqnarray*}
 is a well defined isomorphism of partially ordered sets. Moreover, $\Fa{Y}$ is a lattice.
\qed

Let $\mathcal R$ be a monoid. We equip its set of idempotents $E({\mathcal R})$ with the partial order
\begin{equation*}
e\leq f :\iff ef=fe=e,
 \end{equation*}
where $e,f\in E({\mathcal R})$. Its unit group ${\mathcal R}^\times$ acts on $(E({\mathcal R}),\leq)$ by conjugation.

A monoid $\mathcal R$ is called factorisable, if $E({\mathcal R})$ is a commutative submonoid of $\mathcal R$ and
\begin{equation*}
   { \mathcal R } =  {\mathcal R}^\times E({\mathcal R}) =  E({\mathcal R})  {\mathcal R}^\times  .
\end{equation*}
This implies that $(E({\mathcal R}),\leq)$ is a lower semilattice which has the identity $1$ as biggest element; see Proposition 1.3.2 of \cite{Ho}.

As a common concept for various sorts of monoids, which have all been called Renner monoids in the literature, E. Godelle introduced in Definition 1.4 of \cite{Go} generalized Renner-Coxeter systems:
A triple ($\mathcal R$, $\Lambda$, $\mathcal S$) is called a generalized Renner-Coxeter system if the following holds:
\begin{itemize}
\item[(a)] $\mathcal R$ is a factorisable monoid.
\item[(b)] (${\mathcal R}^\times$, $\mathcal S$) is a Coxeter system.
\item[(c)] $\Lambda$ is a sub-semilattice of $E({\mathcal R}) $, and a cross section for action of ${\mathcal R}^\times$ on $E({\mathcal R}) $ by conjugation.
\item[(d)] For every pair $e_1\leq e_2$ in $E({\mathcal R}) $ there exists $x\in {\mathcal R}^\times$ and a pair $f_1\leq f_2$ in $\Lambda$, such that
$e_1=x f_1 x^{-1}$ and $e_2=x f_2 x^{-1}$.
\item[(e)] For every $e\in\Lambda$ the subgroups
    \begin{equation*}
  C_{{\mathcal R}^\times}(e) :=  \{ x \in {\mathcal R}^\times \mid x e = e x\} \quad \text{ and }\quad  S_{{\mathcal R}^\times}(e) :=  \{ x \in {\mathcal R}^\times \mid x e = e x=e\}
\end{equation*}
of ${\mathcal R}^\times$ are standard parabolic.
\item[(f)] The map $\la^*:(\Lambda,\leq) \to ({\mathcal P}({\mathcal S}),\subseteq)$, defined on $e\in\Lambda$ by
\begin{equation*}
    \la^*(e) := \{ s\in {\mathcal S} \mid s e = e s \neq e \},
\end{equation*}
is a morphism of partially ordered sets.
\end{itemize}
The monoid $\mathcal R$ is called a generalized Renner monoid. To describe the centralizers and stabilizers in (e) define maps $\la, \la_*: \Lambda \to {\mathcal P}({\mathcal S}) $ on $e\in\Lambda$ by
\begin{equation*}
    \la(e) := \{ s\in {\mathcal S} \mid s e = e s  \} \quad\text{ and }\quad \la_*(e) := \{ s\in {\mathcal S} \mid s e = e s = e \}.
\end{equation*}
Then $ C_{{\mathcal R}^\times}(e)$, $S_{{\mathcal R}^\times}(e) $ are the standard parabolic subgroups associated to $\la(e)$, $\la_*(e)$, respectively. Moreover, $\la(e)=\la_*(e)\,\dot{\cup}\,\la^*(e)$. We call $\la, \la_*, \la^*:\Lambda \to {\mathcal P}({\mathcal S}) $ the type map, the lower type map, the upper type map, respectively. 

We now associate to the $\We$-invariant convex subcone $Y$ of the Tits cone $X$ a generalized Renner monoid by the following construction. The semidirect product $\We\ltimes \Fa{Y}$ consists of the set $\We\times\Fa{Y}$ equipped with the structure of a monoid with unit $(1,Y)$ by 
\begin{equation*}
  (\sigma, R )\cdot(\tau, T) \;:=\; (\sigma\tau, \tau^{-1}R\cap T) .
\end{equation*}
It is easy to see that we get a congruence relation on $\We\ltimes \Fa{Y}$ by
\begin{equation*}
   (\sigma, R) \sim (\sigma', R')   \quad :\iff \quad  R=R'\quad\mb{ and }\quad \sigma^{-1}\sigma' \in Z_{\mathcal W}(R) .
\end{equation*}
We denote the quotient monoid by
\begin{equation*}
   {\mathcal R}(Y):=( \We\ltimes \Fa{Y} )/\sim.
\end{equation*}
We denote the congruence class of $(\sigma, R)$ by $\sigma \ve{R}$. 

To state the next Corollary in an easy way, we restrict to the case where $\We$ acts faithfully on $Y$.
\begin{cor} If $\We$ acts faithfully on $Y$, then (${\mathcal R}(Y)$, $\Mklz{1\ve{R}}{R\in\Upsilon}$, $\Mklz{\sigma_i \ve{Y}}{i\in I} $) is a generalized Renner-Coxeter system. 

The unit group ${\mathcal R}(Y)^\times$ is isomorphic to the Coxeter group $\We$, and the semilattice of idempotents $E({\mathcal R}(Y))$ is isomorphic to the face lattice $\Fa{Y}$ by the maps
\begin{eqnarray*}
 \begin{array}{ccc}
         \We &\to     & {\mathcal R}(Y)^\times   \\                       
      \sigma &\mapsto & \sigma \ve{Y}\;
  \end{array}
 & \quad \mbox{ and } \quad &
  \begin{array}{ccc}      
    \Fa{Y} &\to     & E({\mathcal R}(Y))\\
      R  &\mapsto & 1\,\ve{R}
  \end{array} 
 \;\;. 
 \end{eqnarray*}

If we identify $\Mklz{1\ve{R}}{R\in\Upsilon}$ with $\Upsilon$, and $\Mklz{\sigma_i\ve{Y}}{i\in I}$ with $I$, then the type, lower type, upper type maps of the generalized Renner-Coxeter system coincide with $\upsilon,\upsilon_*,\upsilon^*:\Upsilon\to {\mathcal P}(I)$, respectively.
\end{cor}
\Proof The Coxeter group $\We$ acts faithfully on $Y$ if and only if $\upsilon_*(Y)=\emptyset$. Now the proposition is easy to check. In particular, part (c) of the definition of a generalized Renner-Coxeter system follows from Corollary \ref{csp} and Theorem \ref{csublattice}; part (d) follows from the theorem on chains, Theorem \ref{chainconj}, for the length $m=2$; part (f) follows from Proposition \ref{inclspeziell}.
\qed

%
%
\subsection{The lattice operations of the face lattice\label{lattice operations of the face lattice}}
%
%
%
One can show that the type map $\la:\Lambda\to P({\mathcal S})$ of a generalized Renner-Coxeter system ($\mathcal R$, $\Lambda$, $\mathcal S$) as defined in \cite{Go} is a regular type map as defined in \cite{P2}, \cite{AP}, the associated semilattice ${\mathcal F}(\Lambda,\la)$ isomorphic to the semilattice of idempotents $E(\mathcal R)$. In the proof of Theorem 1 of \cite{AP}, as well as in Corollary 1.13 (iii) of \cite{Go}, there has been given a  description of certain lattice meets, which is, transfered to our situation, the following: If $R_1,\,R_2\in\Upsilon$ and $\tau\in\mb{}^{\upsilon(R_1)}\We^{\upsilon(R_2)}$ then $ R_1\cap \tau  R_2$ is the biggest element of 
\begin{equation}\label{lop0}
    \Mklz{R\in\Upsilon}{ R\subseteq R_1,\;R\subseteq R_2,\;\red{\tau}\subseteq \upsilon(R)} .   
\end{equation}
We improve this formula slightly, replacing $v(R)$ by a smaller set. Our proof is obtained by extracting the parts of the proof of Theorem \ref{FXlmj} for the face lattice of the Tits cone $X$, which work also for the face lattice of the invariant convex subcone $Y$ of $X$. 
\begin{thm} Let $R_1,\,R_2\in\Upsilon$ and $\tau\in\mb{}^{\upsilon(R_1)}\We^{\upsilon(R_2)}$. Then we have:
\begin{itemize}
\item[(a)] $ R_1\cap \tau  R_2$ is the biggest element of 
\begin{equation}\label{lop1}
    \Mklz{R\in\Upsilon}{ R\subseteq R_1,\;R\subseteq R_2,\;\text{\rm red}(\tau)\subseteq \upsilon_*(R)} .   
\end{equation}
In particular, $ R_1\cap \tau  R_2\in\Upsilon$ with
\begin{equation}
     \upsilon_*(R_1)\cup\upsilon_*(R_2)\cup\text{\rm red}(\tau) \subseteq \upsilon_*(R_1\cap\tau R_2).
\end{equation}
\item[(b)] $ R_1\vee \tau  R_2$ is the smallest element of 
\begin{equation}\label{lop2}
    \Mklz{R\in\Upsilon}{ R_1\subseteq R,\;R_2\subseteq R,\;\text{red}(\tau)\subseteq \upsilon^*(R)}.
\end{equation}
In particular, $ R_1\vee \tau  R_2\in\Upsilon$ with
\begin{equation}
     \upsilon^*(R_1)\cup\upsilon^*(R_2)\cup\text{\rm red}(\tau) \subseteq \upsilon^*(R_1\vee\tau R_2).
\end{equation}
\end{itemize}
\end{thm}
\begin{rem} The following proof shows that (a) even holds for $\tau\in\mb{}^{\upsilon^*(R_1)}\We^{\upsilon^*(R_2)}$, and (b) even holds for $\tau\in\mb{}^{\upsilon_*(R_1)}\We^{\upsilon_*(R_2)}$.

The set (\ref{lop1}) is contained in the set (\ref{lop0}) and, trivially, all elements of the set (\ref{lop0}) are smaller than or equal to $R_1\cap\tau R_2$. Therefore, the description of $R_1\cap\tau R_2$ as the biggest element of (\ref{lop1}) implies the description as the biggest element of (\ref{lop0}).

The lattice meet and lattice join of two arbitrary faces of $Y$ can be easily reduced to (a) and (b) by Corollary \ref{csp} and Theorem \ref{maincsx} (a).
\end{rem}
\Proof (a) and (b) can be proved similarly; we only show (a). Let $R\in\Upsilon$ be contained in the set (\ref{lop1}). From Theorem \ref{maincsx} (a) we get $\tau^{-1}R=R\subseteq R_2$., which shows that $R\subseteq \tau R_2$. We conclude that 
\begin{equation}\label{lopmax}
      R\subseteq R_1\cap\tau R_2 .
\end{equation} 
Furthermore, we obtain $\upsilon_*(R)\supseteq \upsilon_*(R_1)$ and  $\upsilon_*(R)\supseteq \upsilon_*(R_2)$ from Proposition \ref{inclspeziell}. Hence
\begin{equation*}
     \upsilon_*(R_1)\cup\upsilon_*(R_2)\cup\red{\tau} \subseteq \upsilon_*(R).
\end{equation*}

It remains to prove that $R_1\cap\tau R_2$ is contained in the set (\ref{lop1}). By Corollary \ref{csp} there exist $\sigma\in\We$ and $R\in\Upsilon$ such that $R_1\cap\tau R_2 =\sigma R$. In particular, $\sigma R\subseteq R_1$ and $\sigma R\subseteq\tau R_2$. From Theorem \ref{inclallg} we find
\begin{eqnarray*}
  && R\subseteq R_1 \quad \mb{ and }\quad \sigma^{-1}\in\We_{\upsilon_*(R)}\We_{\upsilon^*(R_1)},  \label{lopincls1}\\
  && R\subseteq R_2 \quad \mb{ and } \quad\sigma^{-1}\tau\in\We_{\upsilon_*(R)}\We_{\upsilon^*(R_2)} . \label{lopincls2}
\end{eqnarray*}
Therefore, there exist $a_1, a_2 \in \We_{\upsilon_*(R)}$, $b_1\in\We_{\upsilon^*(R_1)}$, $b_2\in\We_{\upsilon^*(R_2)}$ such that
$ \sigma^{-1} = a_1 b_1$ and $\sigma ^{-1}\tau =a_2 b_2$. By eliminating $\sigma$ we obtain
\begin{equation*}
   \We_{\upsilon_*(R)}\ni a_1^{-1}a_2 =b_1 \tau b_2^{-1} \in \We_{\upsilon^*(R_1)}\tau \We_{\upsilon^*(R_2)} \subseteq \We_{\upsilon(R_1)}\tau \We_{\upsilon(R_2)}.
\end{equation*}
Since $\tau$ is a minimal double coset representative, there exist $c_1\in \We_{\upsilon(R_1)}$, $c_2\in \We_{\upsilon(R_2)}$ such that
\begin{equation*}
     \We_{\upsilon_*(R)}\ni a_1^{-1}a_2 =c_1 \tau c_2 \quad \mb{ and }\quad l (c_1\tau c_2) = l(c_1)+  l(\tau) + l(c_2),
\end{equation*}
see for example Proposition 2.23 of \cite{AB}. We conclude that $\red{\tau}\subseteq \upsilon_*(R)$. Thus $R$ is contained in the set (\ref{lop1}), and from (\ref{lopmax}) we find
\begin{equation*}
R\subseteq R_1\cap\tau R_2 =\sigma R    .
\end{equation*}
From Corollary \ref{fincleq} we get $R=R_1\cap \tau R_2$.
\qed
%
%
\subsection{The values of $\upsilon_*$ and $\upsilon^*$ on the smallest and the biggest face}
%
%
$Y\cap(-Y)$ is the smallest face, and $Y$ is the biggest face of $Y$. Since $Y$ is $\We$-invariant, we have
\begin{equation} \label{type on min max}
\upsilon(Y\cap(-Y)) = I \quad \mb{ and } \quad \upsilon(Y)=I .
\end{equation}
The next two lemmas describe the values of $\upsilon_*$ and $\upsilon^*$ on these faces. 
\begin{lem}\label{vextrval} Let $I$ be connected. Then exactly one of the following three cases holds:
\begin{itemize}
\item[(1)] $I$ is of finite type and $\mb{span}\Mklz{\al_i}{i\in I}\subseteq Y$. In this case we have
\begin{equation*}
   \upsilon_*(R) = \emptyset    \quad \mb{ and }\quad  \upsilon^*(R) = I \quad \mb{ for all }\quad R\in\Upsilon  .
\end{equation*}
\item[(2)] $Y \subseteq \left(\mb{span}\Mklz{h_i}{i\in I}\right)^\bot$. In this case we have
\begin{equation*}
   \upsilon_*(R)  =  I\quad \mb{ and }\quad  \upsilon^*(R) = \emptyset \quad \mb{ for all }\quad R\in\Upsilon .
\end{equation*}
\item[(3)]  $\upsilon_*(Y\cap(-Y)) = I$, $\upsilon^*(Y\cap(-Y)) =\emptyset$ and $\upsilon_*(Y) = \emptyset$,   $\upsilon^*(Y) = I$.
\end{itemize}
\end{lem}
\begin{rem} In (1) the Coxeter group $\We$ acts trivially on the linear space $\h^*/\mb{span}\Mklz{\al_i}{i\in I}$, and $Y/\mb{span}\Mklz{\al_i}{i\in I}$ may be an arbitrary convex cone in $\h^*/\mb{span}\Mklz{\al_i}{i\in I}$.

In (2) the Coxeter group $\We$ acts trivially on the linear space 
\begin{equation*}
    \left(\mb{span}\Mklz{h_i}{i\in I} \right)^\bot = \Mklz{\la\in\h^*}{\la(h_i)=0 \mb{ for all }i\in I},
\end{equation*} 
and $Y$ may be an arbitrary convex cone in $(\mb{span}\Mklz{h_i}{i\in I})^\bot$.

\end{rem}
\Proof From (\ref{type on min max}) we get the following decompositions of $I$ into separated subsets:
\begin{equation}\label{pvsmbig}
  I =   \upsilon_* (Y\cap(-Y)) \cup \upsilon^*(Y\cap(-Y))   \quad \mb{ and }\quad I = \upsilon_*(Y)\cup \upsilon^*(Y) .
\end{equation}
Since $I$ is connected, $\upsilon_*(Y\cap(-Y))$, $\upsilon^*(Y\cap(-Y))$, $\upsilon_*(Y)$, $\upsilon^*(Y)$ can only take the values $\emptyset$, $I$.

(a) The minimal face $Y\cap (-Y)$ of $Y$ is a linear space. Let $\upsilon^*(Y\cap(-Y))=I$, which is equivalent to $\mb{span}\Mklz{\al_i}{i\in I}\subseteq Y\cap(-Y)$, which in turn is equivalent to $\mb{span}\Mklz{\al_i}{i\in I}\subseteq Y$. 
Suppose that $I$ is of affine or indefinite type. Then we get
\begin{equation*}
   \mb{span}\Mklz{\al_i}{i\in I}\subseteq Y \cap(-Y)\subseteq X\cap(-X) =\Mklz{\la\in \h^*}{\la(h_i)=0 \mb{ for all }i\in I}
\end{equation*}
from Corollary \ref{descsf}, which contradicts that $\al_i(h_i)=2$ for all $i\in I$.

For every $R\in\Upsilon$ we have $Y \cap(-Y)\subseteq R $, from which we find $I= \upsilon^*(Y\cap(-Y))\subseteq \upsilon^*(R)$ by Proposition \ref{inclspeziell}. Thus $\upsilon^*(R)=I$ and $ \upsilon_*(R) = \upsilon(R) \setminus\upsilon^*(R)= \emptyset$. 

(b) Let $\upsilon_*(Y) =I$, which is equivalent to $Y\subseteq \Mklz{\la\in\h^*}{\la(h_i)=0 \mb{ for all }i\in I}$. For every $R\in\Upsilon $ we have
$R\subseteq Y$, from which we obtain $\upsilon_*(R)\supseteq \upsilon_*(Y)=I$ by Proposition \ref{inclspeziell}. Thus $\upsilon_*(R)  = I $ and $\upsilon^*(R) = \upsilon(R) \setminus\upsilon_*(R)=\emptyset$.

(c) If $\upsilon^*(Y\cap(-Y))\neq I$ and $\upsilon_*(Y) \neq I$ then we have $\upsilon^*(Y\cap(-Y))=\emptyset $ and $\upsilon_*(Y) =\emptyset$.
From (\ref{pvsmbig}) we find $\upsilon_*(Y\cap(-Y))=I $ and $\upsilon^*(Y) =I $.

(d) We have seen in (a), (b), (c) that at least one of the cases (1), (2), (3) holds. Since the values of the function $v_*$ are different in (1), (2), (3), at most one of these cases holds.
\qed

The following lemma, where we do not assume that $I$ is connected, can be proved similarly as Lemma \ref{vextrval}. 
\begin{lem} The sets $\upsilon^*(Y\cap(-Y)) $ and $\upsilon_*(Y\cap(-Y))=I\setminus \upsilon^*(Y\cap(-Y))$ are unions of connected components of $I$. Furthermore, we have $\upsilon^*(Y\cap(-Y))^0 = \upsilon^*(Y\cap(-Y))$.

The sets $\upsilon^*(Y)$ and $\upsilon_*(Y)=I\setminus \upsilon^*(Y) $ are unions of connected components of $I$. 
\end{lem}
\begin{rem} The Coxeter group $\We$ acts faithfully on $Y$ if and only if  $\upsilon_*(Y)=\emptyset$, if and only if $\upsilon^*(Y)= I$. In particular, this is satisfied if $Y^\circ\neq\emptyset$, or equivalently $Y-Y=\h^*$.
\end{rem}

%
\subsection{The dimension function}
%
%
Let $R\in Fa(Y)$. The dimension $\mb{dim}(R)$ of $R$ is defined as the dimension of its linear hull $R-R$.
\begin{prop}\label{dimfct} The dimension function $\mb{dim}:\Fa{Y}\to \Nn$ is a strongly monotone $\We$-invariant function.
\end{prop}
\Proof Trivially, this function is monotone and $\We$-invariant. Let $R_1,\,R_2\in\Fa{Y}$ such that $R_1 \subseteq R_2$ and $\mb{dim}(R_1)=\mb{dim}(R_2)$. Since $R_1-R_1\subseteq R_2-R_2$, and $\mb{dim}(R_1-R_1)=\mb{dim}(R_2-R_2)$, we find $R_1-R_1=R_2-R_2$. If we intersect with $Y$ on both sides we get $R_1=R_2$.
\qed
%
%
%
%
\subsection{Intervals\label{Intervals}}
%
%
%
Let $R_1,\,R_2\in\Fa{Y}$ such that $R_1\subseteq R_2$. The interval with endpoints $R_1$, $R_2$ is defined by
\begin{equation*}
  [R_1,R_2]:=\Mklz{R\in\Fa{Y}}{R_1\subseteq R \subseteq R_2}.
\end{equation*}
Recall that for the action of $\We$ on $\Fa{Y}$ the pointwise stabilizer of $[R_1,R_2]$, the stabilizer of $[R_1,R_2]$ as a whole are 
\begin{eqnarray*}
  Z_{\mathcal W}([R_1,R_2])  &:=& \Mklz{\sigma\in\We}{\sigma R=R \mb{ for all }R\in[R_1,R_2]},\\
  N_{\mathcal W}([R_1,R_2] ) &:=& \Mklz{\sigma\in\We}{\sigma[R_1,R_2]=[R_1,R_2]},
\end{eqnarray*}
and that $ Z_{\mathcal W}([R_1,R_2])$ is normal in $N_{\mathcal W}([R_1,R_2])$. The quotient group $N_{\mathcal W}([R_1,R_2] )/ Z_{\mathcal W}([R_1,R_2]) $ acts faithfully on $[R_1,R_2]$.

We may restrict our investigation to intervals with endpoints in $\Upsilon$: By Theorem \ref{chainconj} there exist $\tau\in\We$ and uniquely determined $R_1',\,R_2'\in\Upsilon$ such that $R_1=\tau R_1'$ and $R_2=\tau R_2'$. Clearly, we have
\begin{equation*}
   [R_1,R_2] =\tau [R_1',R_2'] .
\end{equation*}
\begin{prop}\label{stabchain1} Let  $R_1,\,R_2\in\Upsilon$ such that $R_1\subseteq R_2$. Then we have
\begin{equation}\label{stab int norm}
   N_{\mathcal W}([R_1,R_2])= \We_{\upsilon(R_1)\cap \upsilon(R_2)}.
\end{equation}
Moreover,
\begin{equation}\label{stab int index dec}
   \upsilon(R_1)\cap\upsilon(R_2) = \upsilon^*(R_1) \, \dot{\cup}\,  \upsilon_*(R_2) \,\dot{\cup}\, (\upsilon_*(R_1)\cap \upsilon^*(R_2))
\end{equation}
is a decomposition into three separated parts, and we have
\begin{equation}\label{stab int cent}
    Z_{\mathcal W}([R_1,R_2])=\We_{\upsilon^*(R_1)\cup  \upsilon_*(R_2)\cup L},
\end{equation}
where $L$ is either empty or a nonempty union of connected components of $\upsilon_*(R_1)\cap \upsilon^*(R_2)$.
\end{prop}
\Proof From Theorem \ref{maincsx} (a) we get $ N_{\mathcal W}(R_1)=\We_{\upsilon(R_1)}$ and $N_{\mathcal W}(R_2) = \We_{\upsilon(R_2)}$.
To prove (\ref{stab int norm}) we have to show
\begin{equation*}
  N_{\mathcal W}([R_1,R_2]) = N_{\mathcal W}(R_1)\cap N_{\mathcal W}(R_2).
\end{equation*}
Let $\sigma\in N_{\mathcal W}([R_1,R_2]) $. Then $\sigma R_1, \sigma R_2 \in [R_1,R_2]$. In particular, $R_1\subseteq \sigma R_1$ and $\sigma R_2\subseteq R_2$, from which we find $\sigma R_1= R_1$ and $\sigma R_2=R_2$ by Corollary \ref{fincleq}.
Now let $\sigma\in  N_{\mathcal W}(R_1)\cap N_{\mathcal W}(R_2) $ and $R\in[R_1,R_2]$. From $R_1\subseteq R\subseteq R_2$ we get
\begin{equation*}
  R_1=\sigma R_1\subseteq \sigma R\subseteq\sigma R_2= R_2,
\end{equation*}
which shows that $\sigma R\in [R_1,R_2]$.

It is easy to see that (\ref{stab int index dec}) is a decomposition whose parts are separated. We have
\begin{equation*}
    Z_{\mathcal W}([R_1,R_2]) = \bigcap_{R\in [R_1,R_2]} N_{\mathcal W}(R)  \subseteq  \We_{\upsilon(R_1)\cap \upsilon(R_2)} ,
\end{equation*}
where $N_{\mathcal W}(R)$ are parabolic subgroups of $\We$ for all $R\in [R_1,R_2]$. Hence $Z_{\mathcal W}([R_1,R_2])$ is a parabolic subgroup of 
$\We_{\upsilon(R_1)\cap \upsilon(R_2)}$, see for example Section 3 of \cite{Nu}.
Since $Z_{\mathcal W}([R_1,R_2])$ is also normal in $\We_{\upsilon(R_1)\cap \upsilon(R_2)}$, it follows easily that it is a standard parabolic subgroup given by the empty set or by a nonempty union of connected components of $\upsilon(R_1)\cap \upsilon(R_2)$.

Let $R\in [R_1,R_2]$. If $i\in\upsilon^*(R_1)$ then $\al_i\in R_1-R_1\subseteq R-R$. Hence $\sigma_i$ fixes $(R-R)\cap Y=R$ as a whole. If $i\in\upsilon_*(R_2)$ then $\sigma_i$ fixes $R_2$, and also $R\subseteq R_2$ pointwise. Since the decomposition (\ref{stab int index dec}) is separated, we have shown (\ref{stab int cent}).
\qed

If $\emptyset\neq J\subseteq I$ then $(\h,\,(H_i, \,L_i)_{i\in J})$ is a linear Coxeter system, and ($\h$, $(h_j)_{j\in J}$, $(\al_j)_{j\in J}$) is a root base of the linear Coxeter system with generalized Cartan matrix  $A_J =(a_{ij})_{i,j\in J}$. Its linear Coxeter group is $\We_J$. We denote by $\overline{C}(J)$ its fundamental chamber and by $X(J)$ its Tits cone.
To avoid case differentiations in the formulation of some results we define
\begin{equation}\label{Cq X emptyset}
     \overline{C}(\emptyset):=\overline{C}-\overline{C}=\h^* \quad\mb{ and } \quad X(\emptyset) := X-\overline{C}=\h^*.
\end{equation}

Some parts of the following Proposition have been given for integral free linear Coxeter systems as Lemma 1.13 of \cite{Loo}. A part of (a) has been given for symmetric linear Coxeter systems as Lemma 4.2 (f) of \cite{D}.
\begin{prop}\label{titssub} Let $J\subseteq I$ be facial. 
\begin{itemize}
\item[(a)]  We have $\overline{C}(J)=\overline{C}-\overline{F}_J$. Moreover, $K\subseteq J$ is facial for $(\h,\,(H_i, \,L_i)_{i\in J})$ if and only if it is facial for $(\h,\,(H_i, \,L_i)_{i\in I})$. In this case  
\begin{equation}\label{titssub eq 1}
       F_K (J) =  F_K + ( \overline{F}_J - \overline{F}_J ) \quad \mb{ and }\quad \overline{F}_K(J)=\overline{F}_K-\overline{F}_J .
\end{equation}
\item[(b)]  We have $X(J)= X-\overline{F}_J$.
\end{itemize}
\end{prop}
\Proof The proposition holds for $J=\emptyset$ by the definitions (\ref{Cq X emptyset}). Let $\emptyset\neq J\subseteq I$ be facial and $K\subseteq J $. We first show that
\begin{equation}\label{titssub eq 3}
       F_K (J) =  F_K + ( \overline{F}_J - \overline{F}_J )  .
\end{equation}
Since $J$ is facial we find from Proposition \ref{rihfc} that there exists $\mu_0\in F_J=\mb{ri}(\overline{F_J})$, and we have $ \overline{F}_J - \overline{F}_J   =\{\la\in\h^*\mid \la(h_j) = 0 \mb{ for all } j\in J \}$.
Hence the inclusion ``$\supseteq$" in (\ref{titssub eq 3}) is trivial. Clearly, the reverse inclusion ``$\subseteq$" in (\ref{titssub eq 3}) holds for $F_K (J)=\emptyset$. If $\la\in  F_K (J)$ then we have  
\begin{equation*}
   \la= (\la+ t \mu_0 ) - t\mu_0 \in F_K  - \R^+\mu_0 \subseteq F_K + ( \overline{F}_J - \overline{F}_J )
\end{equation*}
all sufficiently big $t\in\R^+$. 

From (\ref{titssub eq 3}) and Proposition \ref{fshfree} (a) we conclude that $K\subseteq J$ is facial for $(\h,\,(H_i, \,L_i)_{i\in J})$ if and only if it is facial for $(\h,\,(H_i, \,L_i)_{i\in I})$. Moreover, the first formula of (\ref{titssub eq 1}) holds. The second formula of  (\ref{titssub eq 1}) can be proved similarly.

Part (b) can be proved similarly to the corresponding part of Lemma 1.13 of \cite{Loo} if the description (\ref{chTcrf1}) of the Tits cone given in Theorem \ref{fuprop} is used. 
\qed

The following theorem is the main theorem on intervals. It identifies an interval with the face lattice of some invariant convex subcone of a Tits cone, for which the whole theory developed so far can be applied.
\begin{thm}\label{mainint} Let $R_1,\,R_2\in\Upsilon$ such that $R_1\subseteq R_2$ and set $J:=\upsilon_*(R_1)\cap\upsilon^*(R_2)$. 
\begin{itemize}
\item[(a)] The convex cone $R_2-R_1$ is a $\We_J$-invariant subcone of the Tits cone $X(J)$ such that
\begin{equation*}
    ( R_2 - R_1 ) \cap \overline{C}(J) \,=\,   (R_2\cap\overline{C}) + ( R_1 - R_1 ) \,= \, (R_2\cap\overline{C})  - ( R_1 \cap\overline{C} ) .
\end{equation*}
\item[(b)] A $\We_J$-equivariant isomorphism of partially ordered sets is given by:
\begin{equation}\label{map int}
\begin{array}{ccc}
    [R_1,R_2] &\to & \Fa{R_2-R_1}  \\
       R\quad &\mapsto &  R-R_1
\end{array}
\end{equation}
It also preserves the dimension of the faces.
\end{itemize}
\end{thm}
\Proof To (b):  We have $R_2-R_1= R_2+(R_1-R_1)$, where $R_1-R_1$ is a linear space. If we show
\begin{equation}\label{int Lem vor}
   \{  F\in \Fa{R_2} \mid (R_1-R_1)\cap (R_2 - F) \subseteq F-F \}  = [R_1,R_2] ,
\end{equation}
then we find from Lemma \ref{concap} that (\ref{map int}) is an isomorphism of partially ordered sets.

Let $F\in \Fa{R_2}$ such that $(R_1-R_1)\cap (R_2 - F) \subseteq F-F$. Intersecting both sides with $R_2$ we find $R_1\subseteq F$. Hence
in (\ref{int Lem vor}) the inclusion ``$\subseteq $" holds. If $F\in  [R_1,R_2]$ then trivially $R_1-R_1\subseteq F-F$. Hence in (\ref{int Lem vor}) the inclusion ``$\supseteq $" holds.

The map (\ref{map int}) is $\We_J$-equivariant because of $\We_J\subseteq \We_{\upsilon_*(R_1)}=Z_{\mathcal W}(R_1)$. It is dimension preserving because for $R\in [R_1,R_2] $ the linear space $R-R$ is the linear hull of $R$ and of $R-R_1$. 

To (a): To prepare the proof we first show some formulas. We get
\begin{equation} \label{intfa0}
    Y - R_1  \,=\, Y + (R_1-R_1)\,= \, Y  +\left( ( R_1 \cap\overline{C} )  -  ( R_1 \cap\overline{C} ) \right) =  Y  - ( R_1 \cap\overline{C} ) .
\end{equation}
from Corollary \ref{mainlh}.  Similarly, we find
\begin{equation}
\label{intfa1}
    R_2 - R_1  \,= \, R_2  - ( R_1 \cap\overline{C} ) ,
\end{equation}
and also
\begin{eqnarray} 
   (Y\cap\overline{C}) + ( R_1 - R_1 )   &=  &(Y\cap\overline{C})  - ( R_1 \cap\overline{C} ) ,  \label{intfa2}\\
   (R_2\cap\overline{C}) + ( R_1 - R_1 ) &= &  (R_2\cap\overline{C})  - ( R_1 \cap\overline{C} ) .\label{intfa3}
\end{eqnarray}
Furthermore, we show that the set in (\ref{intfa3}) is the intersection of the sets in (\ref{intfa1}), (\ref{intfa2}), i.e.,
\begin{equation}\label{intfa4}
    ( R_2 - R_1) \cap \left( (Y\cap\overline{C}) + ( R_1 - R_1 ) \right)   \,=\, (R_2\cap\overline{C})  - ( R_1 \cap\overline{C} ) .
\end{equation}
Trivially, the inclusion ``$\supseteq$" holds. To show the reverse inclusion we use the equalities  (\ref{intfa1}), (\ref{intfa2}). Let $\la_2\in R_2$, $\la\in Y\cap\overline{C}$, and $\la_1,\,\la_1'\in R_1\cap\overline{C}$ such that
\begin{equation*}
   \la_2-\la_1=\la-\la_1'.
\end{equation*}
Since $R_1\subseteq R_2$, we get
\begin{equation*}
   R_2\ni \la_2+\la_1' =\la +\la_1 \in \overline{C},
\end{equation*}
which shows that
\begin{equation*}
   \la_2-\la_1=(\la_2+\la_1') -(\la_1+\la_1')  \in (R_2\cap\overline{C}) -(R_1\cap\overline{C}).
\end{equation*}

We set $L:=\upsilon_*(R_1)$, which is a facial subset of $I$ by Theorem \ref{maincsx} (b). The convex cone $Y  - R_1$ is $\We_L$-invariant, and from equation (\ref{intfa0}), Theorem \ref{maincsx} (b), and Proposition \ref{titssub} we find
\begin{equation*} 
     Y-R_1 \,= \, Y  - ( R_1 \cap\overline{C} ) \,= \,Y  - ( R_1 \cap\overline{F}_L )\, \subseteq\,  X  -  \overline{F}_L =  X(L). 
\end{equation*}

We now show that
\begin{equation}\label{zwintb}
    Y - R_1= \We_L \left( (Y\cap\overline{C}) + ( R_1 - R_1 )  \right)  .
\end{equation}
The inclusion ``$\supseteq$" is obvious. $R_1$ is fixed pointwise by $W_L$. Hence for the reverse inclusion it is sufficient to show that $Y$ is contained in the set on the right in (\ref{zwintb}). Let $\la\in Y$. Choose $\mu_0\in R_1\cap F_L$, which is possible by Theorem \ref{maincsx} (b). We have $\mu_0\in F_L=\mb{ri}(\overline{F}_L)$ by Proposition \ref{rihfc}. It is easy to see that
\begin{equation*}
   \overline{C}- \overline{F}_L=\overline{C} +(\overline{F}_L-\overline{F}_L)=\overline{C}+(\overline{F}_L- \R^+_0 \mu_0) = \overline{C} - \R^+_0\mu_0.
\end{equation*}
It follows from Proposition \ref{titssub} that there exist $\sigma\in\We_L$, $\mu\in\overline{C}$, $r\in\R^+_0$ such that $ \la=\sigma(\mu- r\mu_0)$. We obtain
\begin{equation*}
  \overline{C}\ni \mu= \sigma^{-1}\la+ r\mu_0 \in Y+ R_1\subseteq Y ,
\end{equation*}
which shows that $\mu\in Y\cap\overline{C}$. Thus $\la=\sigma (\mu-r\mu_0)$ is contained in the set on the right in (\ref{zwintb}).

By (\ref{intfa2}), by Theorem \ref{maincsx} (b), and Proposition \ref{titssub} we get
\begin{equation*}
          (Y\cap\overline{C}) + ( R_1 - R_1 )= (Y\cap\overline{C})  - (R_1\cap\overline{C}) = (Y\cap\overline{C})  - (R_1\cap\overline{F}_L)  \subseteq \overline{C} - \overline{F}_L =  \overline{C}(L) .
\end{equation*}
Because $\overline{C}(L)$ is a fundamental domain for the $\We_L$-action on $X(L)$, we find from (\ref{zwintb}) that
\begin{equation}\label{intfa5}
   ( Y - R_1 ) \cap \overline{C}(L) =   (Y\cap\overline{C}) + ( R_1 - R_1 ) .
\end{equation}

Since $R_2\in [R_1,Y]$, it follows from (b) that $R_2-R_1$ is a face of $Y - R_1$. The linear hull of $R_2-R_1$ is $R_2-R_2$. Hence $\mb{ri}(R_2 -R_1)\supseteq\mb{ri}(R_2)$, which has nonempty intersection with $\overline{C}\subseteq \overline{C}(L)$, and
\begin{equation*}
    \upsilon^*(R_2-R_1)=\Mklz{i\in L}{\al_i\in R_2-R_2 } = L\cap \upsilon^*(R_2)=J.
\end{equation*}
From Theorem \ref{maincsx} (b) we obtain
\begin{equation*}
   R_2-R_1 = \We_J \left( (R_2 - R_1)\cap \overline{C}(L)   \right) ,
\end{equation*}
where in addition
\begin{equation*}
   (R_2 - R_1)\cap \overline{C}(L)  = (R_2-R_1)\cap (Y-R_1)\cap  \overline{C}(L) =  (R_2\cap\overline{C})  - ( R_1 \cap\overline{C} ) 
\end{equation*}
by (\ref{intfa5}) and (\ref{intfa4}). We conclude from $\overline{C}(L)\subseteq \overline{C}(J)$ that $R_2-R_1$ is a $W_J$-invariant convex subcone of $X(J)$. Since $\overline{C}(J)$ is a fundamental domain for the $\We_J$-action on $X(J)$, we get
\begin{equation*}
      (R_2 - R_1)\cap \overline{C}(J)   =  (R_2 - R_1)\cap \overline{C}(L) .    
\end{equation*}
\qed
%
%
%
%
\subsection{Chains in intervals\label{Chains in intervals}}
%
%

If $Y$ is a Weyl group invariant polyhedral cone obtained by a normal reductive linear algebraic monoid as explained in the introduction, then all maximal chains of $\Fa{Y}$ have the same length. Furthermore, the lower and upper type maps are determined by the type map, i.e., for $S\in\Upsilon $ we have
\begin{equation*}
    \upsilon_*(S)\,=\,\bigcap_{R\in\Upsilon, \, R\subseteq S} \upsilon(R)  \quad \mb{ and }\quad  
    \upsilon^*(S)  \,=\,\bigcap_{R\in\Upsilon,\, R\supseteq S} \upsilon(R)  \, ;
\end{equation*}
see Section 7, in particular, Lemma 7.15 of \cite{R3}.
If $Y$ is an arbitray Coxeter group invariant convex subcone of the Tits cone $X$ which contains the zero, then these properties do not need to hold, even if the convex cone $Y\cap \overline{C}$ is finitely generated. In particular, these properties do not need to hold for the Tits cone $X$ itself, which is easy to see. However, we show in this subsection that these properties can be generalized for certain parts of $\Fa{Y}$, which will be reached in Corollary \ref{ketteninfapart} below.

We first investigate chains of faces of $Y$, where the dimensions of successive faces differ by one.

\begin{thm}\label{fupv} Let $R_1,\, R_2\in\Upsilon$ such that
\begin{equation*}
  R_1\subseteq R_2 \quad \mb{ and }\quad \mb{dim}(R_1)=\mb{dim}(R_2)-1.
\end{equation*}
Then the following equivalent equations hold:
\begin{itemize}
\item[(a)] $\upsilon_*(R_1)\cap \upsilon^*(R_2)=\emptyset$.
\item[(b)] $\upsilon_*(R_2)=\upsilon_*(R_1)\cap\upsilon(R_2)$.
\item[(c)] $\upsilon^*(R_1)=\upsilon(R_1)\cap\upsilon^*(R_2)$.
\item[(d)] $\upsilon(R_1)\cap\upsilon(R_2)=\upsilon^*(R_1)\cup\upsilon_*(R_2)$.
\end{itemize}
\end{thm}
\Proof We first show (a) indirectly. Suppose there exists $i\in \upsilon_*(R_1)\cap \upsilon^*(R_2)$. Clearly, $R_1-R_1\subseteq R_2-R_2$. Since $i\in  \upsilon_*(R_1)$, we have 
\begin{equation*}
   R_1-R_1\subseteq \Mklz{\la\in\h^*}{\la(h_i)=0} .
\end{equation*}
Since $i\in\upsilon^*(R_2)$, we have $\al_i\in R_2-R_2$. From $\al_i(h_i)=2$ we conclude that
\begin{equation*}
  R_1-R_1 \subseteq \Mklz{\la\in R_2-R_2}{\la(h_i)=0}\subsetneqq R_2-R_2.
\end{equation*}
By the dimension condition, the first inclusion is an equality. Intersecting both sides with $Y$ we obtain
\begin{equation}\label{R1R2hi}
  R_1 =\Mklz{\la\in R_2}{\la(h_i)=0}.
\end{equation}
Now we choose an element $\la\in R_2\setminus R_1$. Since $i\in \upsilon^*(R_2)$, it is $\sigma_i \la \in R_2$. Because $R_1$ is a face of $R_2$ we find $\la+\sigma_i\la\in R_2\setminus R_1$. However, $(\la+\sigma_i\la)(h_i)=0$ contradicts (\ref{R1R2hi}).

To show the equivalence of (a) - (d) we make use of $\upsilon_*(R_1)\supseteq \upsilon_*(R_2)$ and $\upsilon^*(R_1)\subseteq \upsilon^*(R_2)$ without mentioning it further. The equivalence of (a) and (b) follows from the equation
\begin{eqnarray*}
 &&  \upsilon_*(R_1) \cap\upsilon(R_2) = \upsilon_*(R_1) \cap \left(\upsilon_*(R_2)\,\dot{\cup}\,\upsilon^*(R_2)\right) =\left( \upsilon_*(R_1)\cap \upsilon_*(R_2)\right)\dot{\cup}\left( \upsilon_*(R_1)\cap\upsilon^*(R_2)\right) \\
 &&=\upsilon_*(R_2)\,\dot{\cup}\left( \upsilon_*(R_1)\cap\upsilon^*(R_2)\right) .  
\end{eqnarray*}
The equivalence of (a) and (c) follows from the equation
\begin{eqnarray*}
 &&  \upsilon(R_1) \cap\upsilon^*(R_2) = \left(\upsilon_*(R_1)\,\dot{\cup}\,\upsilon^*(R_1)\right)\cap  \upsilon^*(R_2) 
         =\left( \upsilon_*(R_1)\cap \upsilon^*(R_2)\right)\dot{\cup}\left( \upsilon^*(R_1)\cap\upsilon^*(R_2)\right) \\
 &&=\left( \upsilon_*(R_1)\cap\upsilon^*(R_2)\right)  \dot{\cup} \,\upsilon^*(R_1).  
\end{eqnarray*}
The equivalence of (a) and (d) follows from the equation
\begin{eqnarray*}
 && \upsilon(R_1)\cap\upsilon(R_2) = \left(  \upsilon_*(R_1)\,\dot{\cup}\,  \upsilon^*(R_1)\right) \cap \left(  \upsilon_*(R_2)\,\dot{\cup}\,  \upsilon^*(R_2)\right)\\
  && =  \left(\upsilon_*(R_1)\cap\upsilon_*(R_2)\right)\,\dot{\cup}\, \left(\upsilon_*(R_1)\cap\upsilon^*(R_2)\right)\,\dot{\cup}\, \left(\upsilon^*(R_1)\cap\upsilon_*(R_2)\right)\,\dot{\cup}\, \left(\upsilon^*(R_1)\cap\upsilon^*(R_2)\right) \\
&& = \upsilon_*(R_2)\,\dot{\cup}\, \left(\upsilon_*(R_1)\cap\upsilon^*(R_2)\right)\,\dot{\cup}\, \emptyset \,\dot{\cup}\, \upsilon^*(R_1).
\end{eqnarray*}
\qed

\begin{cor} Let $R\in\Upsilon$ such that $\mb{dim}(R)=\mb{dim}(Y)-1$. Then $R$ is maximal in $\Upsilon\setminus \{Y\}$ and
\begin{equation*}
  \upsilon_*(R)=\upsilon_*(Y) .
\end{equation*} 
Let $S\in\Upsilon$ such that $\mb{dim}(S)=\mb{dim}(Y\cap(-Y))+1$. Then $S$ is minimal in $\Upsilon\setminus\{ Y\cap(-Y)\}$ and
\begin{equation*}
  \upsilon^*(S)=\upsilon^*(Y\cap(-Y)) .
\end{equation*} 
\end{cor}
\Proof Proposition \ref{dimfct} shows that  $R$ is maximal in $\Upsilon\setminus \{Y\}$, and $S$ is minimal in $\Upsilon\setminus\{ Y\cap(-Y)\}$. From $\upsilon(Y)=I$ and Theorem \ref{fupv} (b) we get $\upsilon_*(R)=\upsilon_*(Y) $. From $\upsilon(Y\cap(-Y))=I$ and Theorem \ref{fupv} (c) we obtain $\upsilon^*(S)=\upsilon^*(Y\cap(-Y)) $.
\qed

For $R\in\Upsilon$ we have $\upsilon(R)= \upsilon^*(R)\,\dot{\cup}\,\upsilon_*(R)$, where $ \upsilon^*(R)$ and $\upsilon_*(R)$ are separated. This generalizes for certain chains as follows:
\begin{cor}\label{chainonefo} Let $R_1\subseteq R_2\subseteq\cdots\subseteq R_m$ be a chain in $\Upsilon$ such that
\begin{equation*}
  \mb{dim}(R_k)=\mb{dim}(R_{k-1}) + 1,\quad k=2,\,3,\,\ldots,\,m.
\end{equation*}
Then we have
\begin{equation*}
  \bigcap_{k=1}^m\upsilon(R_k) =\upsilon^*(R_1)\,\dot{\cup}\,\upsilon_*(R_m).
\end{equation*}
Furthermore, $\upsilon^*(R_1)$ and $\upsilon_*(R_m)$ are separated.
\end{cor}
\Proof We prove the Corollary by induction over $m$. Clearly, it holds for $m=1$. Let $m>1$. By induction we get
\begin{eqnarray*}
  &&\;\bigcap_{k=1}^m \upsilon(R_k) = \bigcap_{k=1}^{m-1} \upsilon(R_k) \cap \upsilon(R_m) = \left(\upsilon^*(R_1)\cup\upsilon_*(R_{m-1})\right)\cap \upsilon(R_m) \\
  &&\qquad\qquad\quad=   \left(\upsilon^*(R_1)\cap\upsilon(R_m)\right)\cup  \left(\upsilon_*(R_{m-1})\cap\upsilon(R_m)\right) .
\end{eqnarray*}
By Proposition \ref{inclspeziell} we have $\upsilon^*(R_1)\subseteq\upsilon^*(R_m)\subseteq \upsilon (R_m)$. Therefore, the term in the first bracket is equal to $\upsilon^*(R_1)$. The term in the second bracket is equal to $\upsilon_*(R_m)$ by Theorem \ref{fupv} (b). Furthermore, $\upsilon^*(R_1)\subseteq\upsilon^*(R_m)$ and $\upsilon_*(R_m)$ are disjoint, and even separated.
\qed

\begin{prop}\label{zentwmaxch} Let $R_1,\,R_2\in\Upsilon$ such that $R_1\subseteq R_2$. If there exists a chain of faces between $R_1$ and $R_2$ such the dimensions of successive faces differ by one, then
\begin{equation*}
     Z_{\mathcal W}([R_1,R_2]) =    \We_{\upsilon^*(R_1)\cup  \upsilon_*(R_2)}  .
\end{equation*}
\end{prop}
\Proof Because of Proposition \ref{stabchain1} we only have to show the inclusion ``$\subseteq$". We find from Theorem \ref{chainconj} and Proposition \ref{dimfct} that there exists a chain  $R_1=S_1\subseteq S_2\subseteq\cdots\subseteq S_m=R_2$ in $\Upsilon$ such that $\mb{dim}(S_{k+1})=\mb{dim}(S_k)+1$ for $k=1,\,\ldots,\,m-1$. By Corollary \ref{chainonefo} we have
\begin{equation*}
  \bigcap_{k=1}^m\upsilon(S_k) =\upsilon^*(R_1) \cup \upsilon_*(R_2).
\end{equation*}
If $\sigma\in  Z_{\mathcal W}([R_1,R_2]) $ then $\sigma S_k=S_k$ for all $k=1,\,\ldots,\,m$. Hence by Theorem \ref{maincsx} (a) we obtain
\begin{equation*}
  \sigma \,\in \, \bigcap_{k=1}^m \We_{\upsilon(S_k)}\, =\, \We_{ \bigcap_{k=1}^m\upsilon(S_k)}\,=\,\We_{\upsilon^*(R_1) \cup \upsilon_*(R_2)}.
\end{equation*}
\qed

To draw further conclusions we restrict the form of the convex cone $Y\cap\overline{C}$ in the next theorem.
\begin{thm}\label{mainfgchain} Let the convex cone $Y\cap\overline{C}$ be finitely generated. Let $R_1,\,R_2\in\Upsilon $ such that $R_1\subseteq R_2$. If 
  \begin{equation*}
     (\upsilon_*(R_1)\cap\upsilon^*(R_2))^0\,=\,\upsilon_*(R_1)\cap\upsilon^*(R_2)
\end{equation*} 
then we have:
\begin{itemize}
\item[(a)] The interval $[R_1,R_2]$ is finite, and $\We_{\upsilon_*(R_1)\cap\upsilon^*(R_2)}$ acts faithfully on $[R_1,R_2]$.
\item[(b)] The dimensions of successive faces of every maximal chain in $[R_1,R_2]$ differ by one. In particular, every maximal chain in $[R_1,R_2]$ has the length $\mb{dim}(R_2)-\mb{dim}(R_1)+1$.
\item[(c)] We have
\begin{equation*}
  \upsilon^*(R_1)\,\dot{\cup}\,\upsilon_*(R_2)\,=\,\bigcap_{R\in\Upsilon\atop R_1\subseteq R\subseteq R_2} \upsilon(R)  \, .
\end{equation*}
\end{itemize}
\end{thm}
\Proof We set $J:=\upsilon_*(R_1)\cap\upsilon^*(R_2)$ and use Theorem \ref{mainint} and its notation for the proof. Since $Y\cap\overline{C}$ is a finitely generated convex cone, also its faces $R_2\cap\overline{C}$, $R_1 \cap\overline{C}$, and
\begin{equation*}
    ( R_2 - R_1 ) \cap \overline{C}(J) \,= \, (R_2\cap\overline{C})  - ( R_1 \cap\overline{C} ) 
\end{equation*}
are finitely generated convex cones. By the assumption $J^0=J$ the group $\We_J$ is finite. Therefore,  
\begin{equation*}
    R_2 - R_1  = \We_J\left( ( R_2 - R_1 ) \cap \overline{C}(J) \right) 
\end{equation*}
is a finitely generated convex cone, which is also polyhedral by Theorem 4.23 of \cite{La}. From Corollary 8.5 and 8.6 of \cite{Bro} it follows that $\Fa{R_2-R_1}$ is finite, and the dimensions of successive faces of every maximal chain in $\Fa{R_2-R_1}$ differ by one. Because
\begin{eqnarray*}
    [R_1,R_2] &\to & \Fa{R_2-R_1}  \\
       R\quad &\mapsto & \quad\;\; R-R_1
\end{eqnarray*}
is an isomorphism of partially ordered sets which preserves the dimension of the faces, we find that $[R_1,R_2]$ is finite and has the property described in (b).

From Proposition \ref{stabchain1} and \ref{zentwmaxch} we conclude that   
\begin{equation*} 
   \We_{\upsilon_*(R_1)\cap\upsilon^*(R_2)}\cong N_{\mathcal W}([R_1,R_2])/Z_{\mathcal W}([R_1,R_2])
\end{equation*} 
acts faithfully on $[R_1,R_2]$.

Every $R\in[R_1,R_2]\cap \Upsilon$ is contained in a maximal chain of $[R_1,R_2]\cap\Upsilon$. It follows easily from Theorem \ref{chainconj} that
a maximal chain of $[R_1,R_2]\cap\Upsilon$ is also a maximal chain in $[R_1,R_2]$. In particular,  the dimensions of successive faces differ by one.
From Corollary \ref{chainonefo} we get
\begin{equation*}
   \bigcap_{R\in\Upsilon\atop R_1\subseteq R\subseteq R_2} \upsilon(R) 
     =\; \bigcap_{\text{max.}\:\text{chains } \Gamma \text{ in }  [R_1,R_2]\cap\Upsilon}\;\; \bigcap_{R \in \Gamma}\upsilon(R)
     =\upsilon^*(R_1)\,\dot{\cup}\,\upsilon_*(R_2) .
\end{equation*}
\qed

For every $R\in\Upsilon $ the set $\upsilon_*(R)$ is facial by Theorem \ref{maincsx} (b). Hence we get by Theorem \ref{fsp} a partition
\begin{equation}\label{Updecsp}
   \Upsilon= \dot{\bigcup_{\Th\text{ special facial} } } \Upsilon(\Th)  
\end{equation}
where some of the parts 
\begin{equation*}
   \Upsilon(\Th):=\Mklz{R\in\Upsilon}{\upsilon_*(R)^\infty=\Th},\quad \Th \mb{ special facial},
\end{equation*}
may be empty.

The faces of the Tits cone $X$ are exposed faces. For every face of the Tits cone $X$ we obtain by intersection with $Y$ an exposed face of $Y$. In particular,
\begin{equation*}
     Y\cap R(\Th),\quad\Th\mb{ special facial},
\end{equation*}
are exposed faces of $Y$. 
\begin{thm} Let $\Th\subseteq I$ be special facial. The following are equivalent:
\begin{itemize}
\item[(i)] $ \Upsilon(\Th) \neq\emptyset$.
\item[(ii)] $Y\cap\mb{\rm ri}(R(\Th))\neq\emptyset$.
\end{itemize}
In this case $Y\cap R(\Th)$ is the biggest element of $\Upsilon(\Th)$, and $\upsilon(Y\cap R(\Th))=\Th\cup\Th^\bot$.
\end{thm}
\Proof ``(i)$\Rightarrow$(ii)": Let $R\in  \Upsilon $ such that $\upsilon_*(R)^\infty =\Th$. We find from Theorem \ref{maincsx} (b) and Corollary \ref{riRTht} that $\emptyset\neq R\cap F_{\upsilon_*(R)}\subseteq \mb{ri}(R(\Th))$.

``(ii)$\Rightarrow$(i)":  We show that $Y\cap R(\Th)\in \Upsilon(\Th)$, and that it satisfies the remaining statements of the theorem. Since $Y\cap R(\Th)\in\Fa{Y}$, there exist $\tau\in\We$, $S\in\Upsilon$ such that $Y\cap R(\Th)=\tau S$. By the definition of $R(\Th)$ and by Theorem \ref{maincsx} (b) we get
\begin{equation*}
  \We_{\Th^\bot}(Y\cap\overline{F}_\Th)  = \tau \We_{\upsilon^*(S)}(S\cap \overline{F}_{\upsilon_*(S)})\quad\mb{ and }\quad  S\cap  F_{\upsilon_*(S)} \neq\emptyset.
\end{equation*}
Furthermore, from (ii) and Theorem \ref{FXf} (a) we find that there exists $\Th_f\subseteq\Th^{\bot}$, $(\Th_f)^0=\Th_f$ such that $\We_{\Th^{\bot}} (Y\cap F_{\Th\cup\Th_f} )\neq\emptyset$.
Comparing the types of the facets with nonempty intersection with $Y$ we get $\Th\subseteq \upsilon_*(S)\subseteq \Th\cup\Th_f$, from which we find $\upsilon_*(S)^\infty=\Th$. We also obtain $\upsilon^*(S)\subseteq \upsilon_*(S)^\bot\subseteq \Th^\bot$ and $\upsilon(S)=\upsilon_*(S)\cup\upsilon^*(S)\subseteq\Th\cup{\Th^\bot}$.

Now let $\la\in S\cap  F_{\upsilon_*(S)}$ Then there exist $\la'\in Y\cap\overline{F}_\Th$, $\sigma\in\We_{\Th^{\bot}}\subseteq N_{\mathcal W}(R(\Th))$ such that $\sigma\la'=\tau\la$. It follows that $\la'=\la$ and $\sigma^{-1}\tau\in\We_{\upsilon_*(S)}\subseteq N_{\mathcal W}(S)$. We get 
\begin{equation*}
  Y \cap R(\Th) = \sigma^{-1}( Y \cap R(\Th) )=\sigma^{-1}\tau S=S .
\end{equation*}
From $\upsilon(R(\Th))=\Th\cup\Th^{\bot}$ we obtain the inclusion  $\Th\cup\Th^{\bot}\subseteq \upsilon(Y\cap R(\Th))$.

It remains to show that $Y\cap R(\Th)$ is the biggest element of $\Upsilon(\Th)$. Let $R\in \Upsilon(\Th)$. Since $\upsilon_*(R)$ and $\upsilon^*(R)$ are separated, we get 
$ \Th=  \upsilon_*(R)^\infty\subseteq \upsilon_*(R)$ and $\upsilon_*(R)^0\cup \upsilon^*(R)\subseteq \Th^\bot$. If we use Theorem \ref{maincsx} (b) to describe $R$ we find
\begin{equation*}
   R= \We_{\upsilon^*(R)}(R\cap\overline{F}_{\upsilon_*(R)}) \subseteq \We_{\Th^\bot}\overline{F}_{\Th} = R(\Th) ,
\end{equation*}
which shows that $R\subseteq Y\cap R(\Th)$.
\qed

The partition (\ref{Updecsp}) of $\Upsilon$ induces a partition
\begin{equation}
    \Fa{Y} = \dot{\bigcup_{\Th \text{ special facial} }}  \Fa{Y}_\Th
\end{equation}
where some of the parts 
\begin{equation*}
   \Fa{Y}_\Th :=\We\,\Upsilon(\Th), \quad \Th \mb{ special facial}, 
\end{equation*}
may be empty. 
From Theorem \ref{mainfgchain} we obtain the following corollary.
\begin{cor}\label{ketteninfapart} Let the convex cone $Y\cap\overline{C}$ be finitely generated, and $\Th\subseteq I$ be special facial. 
\begin{itemize}
\item[(a)] If $R_1,\,R_2\in \Fa{Y}_\Th $ such that $R_1\subseteq R_2$, then every maximal chain between $R_1$ and $R_2$ is contained in $\Fa{Y}_\Th$. The dimensions of its successive faces differ by one. In particular, it has the length $\mb{dim}(R_2)-\mb{dim}(R_1)+1$.
\item[(b)] If $R_1,\,R_2\in \Upsilon(\Th)$ such that $R_1\subseteq R_2$, then
\begin{equation*}
  \upsilon^*(R_1)\,\dot{\cup}\,\upsilon_*(R_2)\,=\,\bigcap_{R\in\Upsilon(\Th)\atop R_1\subseteq R\subseteq R_2} \upsilon(R)  \, .
\end{equation*}
\end{itemize}
\end{cor}
\Proof  Let $R_1,\,R_2\in \Fa{Y}_\Th $ such that $R_1\subseteq R_2$. By Theorem \ref{chainconj} there exist $\tau\in\We$ and uniquely determined $R_1',\,R_2'\in\Upsilon$, $R_1'\subseteq R_2'$, such that $R_1=\tau R_1'$ and $R_2=\tau R_2'$. It follows that
\begin{equation*}
   [R_1,R_2] =\tau [R_1',R_2']  \quad\mb{ with } \quad R_1',\,R_2'\in \Upsilon(\Th).
\end{equation*}
We now show
\begin{equation*}
   [R_1',R_2']\subseteq \Fa{Y}_\Th .
\end{equation*}
Let $R'\in [R_1',R_2']$. From Theorem \ref{chainconj} we find that there exist $\sigma\in\We$ and a uniquely determined $R''\in\Upsilon$
such that $R'=\sigma R''$ and
\begin{equation*}
  R_1'\subseteq R''\subseteq R_2'.
\end{equation*}
We find $\upsilon_*(R_1')\supseteq \upsilon_*(R'')\supseteq \upsilon_*(R_2')$, from which we get 
\begin{equation*}
    \Th = \upsilon_*(R_1')^\infty\supseteq \upsilon_*(R'')^\infty\supseteq \upsilon_*(R_2')^\infty=\Th .
\end{equation*}
We conclude that $R''\in\Upsilon(\Th)$ and $R'=\sigma R''\in\Fa{Y}_\Th$.

Therefore, it is sufficient to show the corollary for $R_1,\,R_2\in \Upsilon(\Th)$ such that $R_1\subseteq R_2$. By definition, $\upsilon_*(R_1)^\infty =\upsilon_*(R_2)^\infty =\Th$. Because of $\Th\cap \upsilon^*(R_2)\subseteq \upsilon_*(R_2)\cap \upsilon^*(R_2)=\emptyset$ we get
\begin{equation*}
  \upsilon_*(R_1)\cap\upsilon^*(R_2) = \left(\Th\cup \upsilon_*(R_1)^0 \right) \cap\upsilon^*(R_2) = \upsilon_*(R_1)^0 \cap\upsilon^*(R_2),
\end{equation*}
which is either empty or has only connected components of finite type. Now the Corollary follows from Theorem \ref{mainfgchain}. 
\qed

\begin{cor} Let the convex cone $Y\cap\overline{C}$ be finitely generated. Let $\Th\subseteq I$ be special facial such that $\Upsilon(\Th)\neq\emptyset$. Then every face in $\Upsilon(\Th)$ is an exposed face of the biggest face $Y\cap R(\Th)\in \Upsilon(\Th)$. In particular, the faces in $ \Fa{Y}_{\upsilon_*(Y)^\infty}$ are exposed faces of $Y$.
\end{cor}
\Proof Let $R\in\Upsilon(\Th)$. Then $R$ is a face of $Y\cap R(\Th)$. As we have seen in the proof of Corollary \ref{ketteninfapart} we have
\begin{equation*}
    \left( \upsilon_*(R)\cap \upsilon^*(Y\cap R(\Th)) \right)^0 =  \upsilon_*(R)\cap \upsilon^*(Y\cap R(\Th)).
\end{equation*}
As in the proof of Theorem \ref{mainfgchain} we find that $(Y\cap R(\Th))-R$ is a finitely generated convex cone. Its faces are exposed. In particular, $R-R$ is an exposed face of $(Y\cap R(\Th))-R$. It follows that $R$ is an exposed face of $Y\cap R(\Th)$.
\qed

%
\subsection{Invariant convex subcones containing the dual imaginary cone \label{dimc as subc}}
%
%
%
%
%
%
The following theorem has been given for free integral linear Coxeter systems as Proposition 2.4 in \cite{Loo}. For general linear Coxeter systems it can be proved similarly. 
\begin{thm}\label{corint} Let $\la\in \overline{C}\cap X^\circ$. Then 
\begin{equation*}
               \text{\rm co}(\We\la)\cap\overline{C} = (\la-\sum_{i\in I}\R^+_0\al_i)\cap\overline{C} .
\end{equation*} 
\end{thm}

Similarly to the definition of the imaginary cone $Z\subseteq \h$ in (\ref{defimcone1}), we define the dual imaginary cone $Z^\vee\subseteq\h^*$ by
\begin{equation*}
    Z^\vee:= \bigcup_{\sigma\in {\cal W}}\sigma K^\vee = (-X)\cap  \bigcap_{\sigma\in {\cal W}} \sigma (\sum_{i\in I}\R^+_0 \al_i ) ,
\end{equation*}
where $K^\vee:= (\sum_{i\in I}\R^+_0 \al_i ) \cap (-\overline{C}) $ is a closed convex cone, even finitely generated.
\begin{cor}\label{dimc} If $\We$ acts faithfully on $Y$, and $Y\cap\overline{C}$ is closed, then we have $-Z^\vee\subseteq Y$.
\end{cor}
\Proof It is sufficient to show $-K^\vee\subseteq Y$. Let $\beta\in K^\vee=(\sum_{i\in I}\R^+_0\al_i )\cap ( -\overline{C})$. Since $\We$ acts faithfully on $Y$ we have $\upsilon_*(Y)= \emptyset $. Hence by Theorem \ref{maincsx} (b) there exists an element $\la\in Y\cap C$. By Theorem \ref{corint} we obtain
\begin{equation*}
     \overline{C}\ni r\la-\beta\in \co{ \We (r\la) } \subseteq Y     \quad\mb{ for all }\quad r\in\R^+,
\end{equation*}
because $Y$ is a $\We$-invariant convex cone. Since $Y\cap\overline{C}$ is closed, we get $-\beta\in Y\cap\overline{C}$ as $r\to 0$.
\qed

\begin{cor} Let $I^{ind}=I $. Then $-Z^\vee$ is the smallest $\We$-invariant convex subcone of $X$ containing the zero, such that $\We$ acts faithfully on $-Z^\vee$, and $(-Z^\vee)\cap\overline{C}$ is closed.
\end{cor}
\Proof We have $(-Z^\vee)\cap\overline{C}=-K^\vee$, which is closed.
Since $I^{ind}=I $, we conclude from Theorem 4.3 (Ind) of \cite{K} that there exists $r\in (\R^+)^I $ such that $Ar\in (\R^-)^I $. 
The group $\We$ acts faithfully on $-Z^\vee$, because $ -\sum_{i\in I} r_i\al_i \in (-\sum_{i\in I}\R^+ \al_i )\cap  C \subseteq - K^\vee \subseteq - Z^\vee$ is fixed only by $1\in\We$. Now the corollary follows from Corollary \ref{dimc}.
\qed
%
%
%
%
%

%
\end{document}